\documentclass[11pt]{article}
\usepackage{a4wide}
\usepackage{amsmath, amsthm, xcolor,bm}
\usepackage{amssymb}
\usepackage{graphicx,soul}
\usepackage[normalem]{ulem}
\usepackage{relsize}

\usepackage{cancel}
\usepackage{soul}
\usepackage{dsfont}
\usepackage{hyperref}

\newtheorem{theorem}{Theorem}
\newtheorem{thm}[theorem]{Theorem}
\newtheorem{proposition}{Proposition}
\newtheorem{mainproblem}{Main Problem}

\newtheorem{conjecture}{Conjecture}

\newtheorem{corollary}{Corollary}
\newtheorem{lemma}{Lemma}

\newtheorem{lemmaDBC}{Lemma DBC (Divergence Borel-Cantelli)  \!\!\!\!}

\newtheorem{lemmaGDBC}{Lemma GDBC (General Divergence Borel-Cantelli)\!}

\theoremstyle{remark}
\newtheorem{rem}{{Remark}}
\renewcommand{\Bbb}[1]{\mathbb{#1}}
\newcommand{\N}{{\Bbb N}}         
\newcommand{\Q}{{\Bbb Q}}         
\newcommand{\R}{{\Bbb R}}         
\newcommand{\Z}{{\Bbb Z}}         

\newcommand{\cSl}{\cS'}
\newcommand{\cSke}{{\cS_k^{\ve}}}
\newcommand{\cSkle}{{\cS_{k_\ell}^{\ve}}}
\newcommand{\cSlm}{\cSl_{\neq m}}
\newcommand{\aaa}{c}

\newcommand{\cA}{{\cal A}}

\newcommand{\cF}{{\cal F}}
\newcommand{\cG}{{\cal G}}

\newcommand{\cI}{{\cal I}}

\newcommand{\cK}{{\cal K}}

\newcommand{\cN}{{\cal N}}

\newcommand{\cS}{{\cal S}}

\newcommand{\ve}{\varepsilon}
\newcommand{\supp}{\operatorname{supp}}

\DeclareMathOperator{\lcm}{lcm}

\DeclareMathOperator{\Prime}{\mathbb{P}}

\DeclareMathOperator{\Log}{\operatorname{Log}}

\begin{document}

\large

\title{\bf Borel--Cantelli, zero-one laws and\\  inhomogeneous Duffin--Schaeffer}

\author{Victor Beresnevich \and Manuel Hauke \and Sanju Velani}

\maketitle
\begin{center}
    \date{{\em To intrinsic  beauty:  Khintchine's Theorem - one hundred  years on!}}
\end{center}

\maketitle

\begin{abstract}
The most versatile version of the classical  divergence Borel--Cantelli lemma shows that for any divergent sequence of events $E_n$ in a probability space satisfying a quasi-independence condition, its corresponding  limsup set $E_\infty$ has positive probability. In particular, it provides a lower bound on the probability of $E_\infty$. In this paper we establish a new version of this classical result which guarantees, under an  additional mild assumption, that the probability of $E_\infty$ is not just positive but is one. Unlike existing optimal results, it is applicable within the setting of  arbitrary probability spaces. We then go onto to consider a range of  applications in number theory and dynamical systems. These include new results on the inhomogeneous Duffin--Schaeffer conjecture.  In particular, we establish alternatives to the classical (homogeneous) zero-one laws of Cassels and Gallagher and use them to resolve the so-called weak Duffin-Schaeffer conjecture for an arbitrary rational inhomogeneous shift. As a bi-product,  we establish the Duffin--Schaeffer conjecture with congruence relations.  The applications to dynamical systems include new characterisations of Borel--Cantelli sequences and new dynamical Borel--Cantelli lemmas, as well as characterising  Khintchine-type sequences for shrinking targets.
\end{abstract}

\bigskip
\newpage 

\tableofcontents

\newpage

\section{Introduction}  \label{introA}

Let $(X,\cA,\mu)$ be a probability space  and
let $\{E_i\}_{i \in \mathbb{N}}$ be a sequence of  measurable subsets (events) of $X$.
Borel--Cantelli lemmas represent a powerful tool in probability theory with wide ranging applications which is used for establishing the probability measure of the limsup set
$$
E_{\infty}:=\limsup_{i \to \infty} E_i:= \bigcap_{t=1}^{\infty}
\bigcup_{i=t}^{\infty} E_i \; ;
$$
that is,  the set of points $x \in X$ such that $ x \in E_i $ for infinitely many $i \in \N$.

Determining the measure of $E_{\infty}$ represents a fundamental problem considered in probability theory -- see for example   \cite[Chp.1 \S4]{MR1324786}   and    \cite[Chp.47]{MR1265493}. The first or convergence Borel--Cantelli Lemma states that $\mu(E_\infty)=0$ whenever $\sum_{i=1}^\infty\mu(E_i)<\infty$, that is to say that the sequence $\{E_i\}_{i\in\N}$ is \emph{convergent}. Naturally, for $\mu(E_\infty)$ to have full, or even just positive measure, it is necessary that $\sum_{i=1}^\infty\mu(E_i)=\infty$, that is to say that the sequence $\{E_i\}_{i\in\N}$ is \emph{divergent}. 
Understanding  conditions under which  $\mu(E_\infty)=1$ for  divergent sequences 
is a fundamental problem that crops up in many areas of mathematics and beyond.
In particular, the Borel--Cantelli lemmas are one of the most fundamental tools used in number theory  (see for example \cite{MR2508636, MR2184760, MR3618787, MR1672558, MR1764799, MR4425845, MR548467} and references within) and ergodic theory and dynamical systems (see for example \cite{MR2944100, MR1826488, MR2327135, MR2669635, MR1309976, MR1724857,MR4417342,MR2276480,MR4572386, MR2214457} and references within). Specifically, divergence Borel--Cantelli lemmas are at the heart of numerous recent advances in metric number theory, including the recent breakthrough proof of the notorious Duffin--Schaeffer Conjecture \cite{MR4859} given by Koukoulopoulos $\&$ Maynard \cite{MR4125453}. The present work is also motivated by such applications, this time to the more general  inhomogeneous form of the Duffin--Schaeffer Conjecture, which remains an open problem.

The classical divergence Borel--Cantelli Lemma requires the mutual independence (in the pure probabilistic sense) of the events $E_i$ in order to conclude that $\mu(E_\infty)=1$.
However, in most applications mutual independence is rarely satisfied, and hence Borell--Cantelli lemmas with weaker hypothesis are required. There are many versions of the classical divergence Borel--Cantelli Lemma requiring various `degrees' of independence such as pairwise independence and pairwise quasi-independence, see \cite{MR4497313} for an overview. The most versatile version of the divergence Borel--Cantelli Lemma, stated as Lemma~DBC below, imposes the condition of  {\em quasi-independence on average}; namely  \eqref{eqn02} in what follows.

\medskip

\begin{lemmaDBC}
 \! Let $(X,\cA,\mu)$ be a probability space  and let $\{E_i\}_{i \in \N} $ be a sequence of subsets in $\cA$.
Suppose that
\begin{equation}\label{eqn01}
\sum_{i=1}^\infty \mu(E_i)=\infty
\end{equation}
and that there exists a constant $C>0$ such that
\begin{equation}\label{eqn02}
\sum_{s,t=1}^Q  \mu(E_s\cap E_t)\le C\left(\sum_{s=1}^Q  \mu(E_s)\right)^2\quad\text{for infinitely many $Q\in\N$\,.}
\end{equation}
Then
\begin{equation}\label{eq3}
\mu(E_{\infty}) \ge C^{-1}\,.
\end{equation}
\end{lemmaDBC}

\bigskip

\noindent  It is known that in various natural settings, condition \eqref{eqn02} is not only sufficient for $\mu(E_\infty)>0$ but also necessary, up to a certain trimming of the sequence $\{E_i\}_{i\in\N}$. For instance, for sequences of balls $\{E_i\}_{i\in\N}$ in a metric space $X$ equipped with a doubling Borel probability measure $\mu$ this was shown in \cite[Theorem~3]{MR4497313}. Furthermore, without additional assumptions on the sequence $\{E_i\}_{i\in\N}$ the lower bound \eqref{eq3} is sharp and cannot be improved. For example, if $\mu$ is Lebesgue measure on $[0,1]$, $E_i=[0,\tfrac12]$ for each $i$, then \eqref{eqn01} and \eqref{eqn02} hold with $C=2$; we also have that $E_\infty=[0,\tfrac12]$ and so $\mu(E_\infty)=1/C$ meaning that \eqref{eq3} is sharp for this example.
Nevertheless, many applications, particularly within number theory and ergodic theory and dynamical systems, require establishing that $\mu(E_\infty)=1$.   Thus in terms of applications alone, the following is a natural and fundamental question  which we aim to  address in this paper.

\medskip

\textbf{Question:} {\em Under what additional assumptions in Lemma~DBC does $\mu(E_{\infty})=1$? }

\medskip

\noindent Of course, if \eqref{eqn02} holds for any $C>1$ then $\mu(E_{\infty}) =1$. However, verifying \eqref{eqn02} for $C$ arbitrarily close to $1$ is often a challenging task, if at all possible.  Furthermore, if \eqref{eqn02} holds with $C=1$ and a sufficiently small error term, then one can establish an asymptotic  quantitative version of the Borel--Cantelli Lemma; that is to say
\begin{equation}   \label{asympBC} \sum_{i=1}^N \mathds{1}_{E_i}\sim \sum_{i=1}^N \mu(E_i)   \qquad  \text{as } \ N\to\infty
\end{equation}
almost surely, see for example \cite{MR1826488, MR1672558}, \cite[\S1.1]{MR4572386} and reference within.
Within the context of dynamical systems, more precisely the setup of shrinking targets, see also  Theorem~\ref{gencountthm} in \S\ref{stgen} below for  a quantitative statement of the Borel--Cantelli Lemma.

In some applications concluding that $\mu(E_\infty)=1$ from $\mu(E_\infty)>0$ is not an issue as one can use other means (such as Kolmogorov's theorem \cite[Theorems 4.5 \& 22.3]{MR1324786} or ergodicity \cite[\S24]{MR1324786}) which demonstrate that the limsup set  $E_{\infty}$ satisfies a \emph{zero-one law}; namely that $\mu(E_\infty)$ is either $0$ or $1$. For instance, in the case of the Duffin--Schaeffer Conjecture (see \S\ref{DSBack}) appropriate zero-one laws, which play an {important} role in the proof of the Duffin--Schaeffer Conjecture \cite{MR4125453}, were established by Cassels \cite{MR36787} and Gallagher \cite{MR133297}. For further details and higher dimensional generalisation of their zero-one laws, see \cite{MR3105329, MR2457266}.
However, in other applications zero-one laws are not available. In particular, zero-one laws are not available within the inhomogeneous context of the Duffin--Schaeffer Conjecture.  A particular  goal of this paper is to address this imbalance and indeed prove a form of conjecture when the inhomogeneous ``shift'' is rational. Surprisingly, even the innocent looking  rational case seems fought with difficulties.

We remark that without the presence of a zero-one law, a full measure statement can be obtained if we are willing to impose additional structure on the probability space $X$, namely if the measure sum diverges locally and  quasi-independence on average holds locally in the presence of an appropriate topological structure on $X$. We refer the reader to \cite{MR4497313} for further details. The approach explored in this paper is different. It generally does not require any additional structure on the probability space $X$.

Prior to presenting our new results we state a more general version of Lemma~DBC.  Its proof is similar to that of Lemma~DBC and is provided in the appendix for completeness.

\begin{lemmaGDBC}
Let $(X, \cA, \mu) $ be a probability space and let $\{E_i\}_{i\in\N}$ be a sequence of subsets in\/ $\cA$. Suppose that there exist constants $C'>0$ and $\aaa>0$ and {a sequence of finite subsets $\cS_k\subset\Z$ such that $\min\cS_k\to+\infty$ as $k\to\infty$, and such that}
\begin{equation}\label{eqn04}
\sum_{i\in\cS_k}\mu(E_i) \ge \aaa
\end{equation}
and
\begin{equation} \label{eqn05}
\sum_{\substack{s<t\\[0.5ex] s,t\in\cS_k}} \mu\big(E_s\cap E_t \big) \ \le \  C'\,  \left(\sum_{i\in\cS_k}\mu(E_i)\right)^2
\end{equation}
{for all sufficiently large $k\in\N$.}
Then
$$
\mu(E_\infty) \ge \frac{1}{2C'+\aaa^{-1}}\,.
$$
\end{lemmaGDBC}

Note that conditions \eqref{eqn01} and \eqref{eqn02} in Lemma~DBC imply conditions \eqref{eqn04} and \eqref{eqn05} in Lemma~GDBC  by taking $\cS_k=\{k,\dots,Q\}$ for a sufficiently large $Q$ depending on $k$ and satisfying \eqref{eqn02}.  It turns out  that Lemma~GDBC implies Lemma~DBC (see the appendix) and thus the two statements are equivalent.   However, Lemma~GDBC is more general in the sense that it allows us to choose the sets $\cS_k$ flexibly. This in-built flexibility is an important feature when it comes to applications.   In the main body of the paper we will establish a  full-measure version of this lemma; namely Theorem~\ref{T2}.

We now state our new contributions to the ``classical'' body of Borel-Cantelli lemmas within  probability theory.   The new statements  are then subsequently  exploited to establish various results in the theory of inhomogeneous Diophantine approximation  (see \S\ref{DSBack}--\ref{PDSlove}) and the dynamical framework of shrinking targets  (see  \S\ref{STlove}).

\section{New results {on Borel--Cantelli}}\label{sec2}

\subsection{Results for arbitrary probability spaces}

The following new results give full-measure generalisations of Lemmas~DBC and GDBC.

\begin{theorem}\label{T1}
Let $(X,\cA,\mu)$ be a probability space and let $\{E_i\}_{i \in \N} $ be a sequence of subsets in $\cA$. Suppose  that the  conditions of Lemma~DBC hold; namely \eqref{eqn01} and \eqref{eqn02} hold for some $C>0$.
In addition,  suppose  that
\begin{itemize}
  \item[]
\begin{itemize}
  \item[{\bf(M1)}] for any $\delta>0$ and any natural numbers $q_1<q_2$ there exists $i_0=i_0(q_1,q_2,\delta)$ such that for all $i\ge i_0$
\begin{equation}\label{eqn03}
\mu\left(A\cap E_i\right)\le (1+\delta)\mu\left(A\right)\mu(E_i)\,,\qquad\text{where }A=\bigcup_{j=q_1}^{q_2}E_j\,.
\end{equation}
\end{itemize}
\end{itemize}
Then $\mu(E_\infty)=1$.
\end{theorem}

\bigskip

\begin{theorem}\label{T2}
Let $(X,\cA,\mu)$ be a probability space and let $\{E_i\}_{i \in \N} $ be a sequence of subsets in $\cA$. Suppose that the  conditions of Lemma~GDBC hold, namely
there exist constants $C'>0$ and $\aaa>0$ and {a sequence of finite subsets $\cS_k\subset\Z$ such that $\min\cS_k\to+\infty$ as $k\to\infty$, and such that} \eqref{eqn04} and \eqref{eqn05} hold {for all sufficiently large $k\in\N$}.
In addition,  suppose that {\bf(M1)} holds. Then
$\mu(E_\infty)=1$.
\end{theorem}

We emphasize that Theorems~\ref{T1} and \ref{T2} do not require any additional structure on $X$.  Current full measure results  require that $X$ is a metric space since they make use of Lebesgue density, see for example \cite[\S8]{MR2184760}, or \cite{MR4497313} and references within.

In order to meet {\bf(M1)} the sets $E_i$ should be ``spread'' well enough. This can be achieved by grouping together `smaller' events to form a new sequence. Below we provide a more general version of Theorem~\ref{T2} that incorporates this idea of `grouping together' by essentially requiring that {\bf(M1)} holds only on average.

\medskip

\begin{theorem}\label{T3}
Let $(X, \cA, \mu) $ be a probability space and let $\{E_i\}_{i \in \N} $ be a sequence of subsets in $\cA$. Suppose that {there exists a constant $C'>0$ such that} for any $\varepsilon>0$ there exist $0<\varepsilon^*<\varepsilon$ and
{a sequence of finite subsets $\cSke\subset\Z$ such that $\min\cSke\to+\infty$ as $k\to\infty$, and} such that
\begin{equation}\label{eq07}
\varepsilon^*\le \sum_{i\in\cSke}\mu(E_i) \le \varepsilon
   \end{equation}
and
\begin{equation} \label{eqn05ve}
{\sum_{\substack{s<t\\[0.5ex] s,t\in\cSke}} \mu\big(E_s\cap E_t \big) \ \le \  C'\, \left(\sum_{i\in\cSke}\mu(E_i)\right)^2}
\end{equation}
{hold for all sufficiently large $k\in\N$}. Furthermore,  suppose that
\begin{itemize}
  \item[]
\begin{itemize}
  \item[{\bf(M2)}]
for any $\delta>0$ and every pair of integers $q_1<q_2$ there exists $k_0=k_0(q_1,q_2,\delta)$ such that for all $k\ge k_0$
\begin{equation}\label{eq08'}
\mu\left(A\cap \bigcup_{i\in\cSke} E_i\right)\le (1+\delta)\mu\left(A\right)\sum_{i\in\cSke}\mu(E_i)\,
\qquad\text{where }\quad
A=\bigcup_{j=q_1}^{q_2}E_j\,.
\end{equation}
\end{itemize}
\end{itemize}
Then $\mu(E_\infty)=1$.
\end{theorem}

\bigskip

\begin{rem}
We  note that {\bf(M2)} is an obvious consequence of the following version of {\bf(M1)} on average and thus is a consequence of {\bf(M1)}:
\begin{itemize}
  \item[]
\begin{itemize}
  \item[{\bf(M2${}'$)}]
for any $\delta>0$ and every pair of integers $q_1<q_2$ there exists $k_0=k_0(q_1,q_2,\delta)$ such that for all $k\ge k_0$
\begin{equation}\label{eq08}
\sum_{i\in\cSke}\mu\left(A\cap E_i\right)\le \sum_{i\in\cSke}(1+\delta)\mu\left(A\right)\mu(E_i)\,,
\qquad\text{where}\quad
A=\bigcup_{j=q_1}^{q_2}E_j\,.
\end{equation}
\end{itemize}
\end{itemize}
\end{rem}

Condition \eqref{eq07} of Theorem~\ref{T3} implies that $\liminf_{i\to\infty}\mu(E_i)=0$. However, this is not particularly restrictive. In any case, the following theorem deals with the complementary situation.

\medskip

\begin{theorem}\label{T4}
Let $(X,\cA,\mu)$ be a probability space  and let $\{E_i\}_{i \in \N} $ be a sequence of subsets in $\cA$. Suppose that there exists $0<c_0<1$ such that
\begin{itemize}
  \item[]
\begin{itemize}
  \item[{\bf(M1-$\cI$)}] for any $\delta>0$ and every pair of integers $q_1<q_2$ there is an infinite subset $\cI\subset\N$ and
      $i_0=i_0(q_1,q_2,\delta)$ such that 
\begin{equation}\label{eq19}
      \inf\{\mu(E_i):{i\in\cI}\}> c_0
\end{equation}
      and \eqref{eqn03} holds for all $i \in \cI$ with $i\ge i_0$.
\end{itemize}
\end{itemize}
Then $\mu(E_\infty)=1$. In particular, $\mu(E_\infty)=1$ if\/ $\limsup\limits_{i\to\infty}\mu(E_i)>0$ and {\bf(M1)} is fulfilled.
\end{theorem}

\medskip

{\begin{rem}[Proof of the `In particular' part assuming the main part of the theorem]\label{rem2}
By $\limsup_{i\to\infty}\mu(E_i)>0$, there exists $\cI\subset\N$ such that $$\inf\{\mu(E_i):{i\in\cI}\}=:2c_0>c_0>0.$$ 
Then {\bf(M1-$\cI$)} follows from {\bf(M1)} for this choice $\cI$ (the same for every $\delta>0$ and $q_1<q_2$).
Therefore, the main part of the theorem implies that $\mu(E_\infty)=1$.
\end{rem}}

\subsection{Results for measure metric spaces}

Now we specialise Theorems~\ref{T1}--\ref{T4} to probability measures $\mu$ supported on a metric space $X$. The goal is to simplify Properties {\bf(M1)} and {\bf(M2)} by replacing $A$ with a ball.

In what follows,  {$\supp \mu $ will denote the support of the measure $\mu$}
and, given $x\in X$ and $r>0$, $B=B(x,r)$ will denote the ball centred at $x$ of radius $r$. {Also, given a real number $a > 0$, we denote by  $a
 B$  the ball $B$ scaled by a factor $a$; i.e.  $aB:= B(x, a r)$.}  Recall, that $\mu $ is said to be \emph{doubling}  if there are constants $ \lambda  \geq 1$ {and $r_0>0$} such that for any $x \in \supp\mu$ and {$0<r<r_0$}
\begin{equation}\label{doub}
\mu (B(x, 2r))  \, \leq \, \lambda \ \mu(B(x,r)) \, .
\end{equation}

\noindent The metric measure space $(X, \cA, \mu, d) $ is also said to be  {\em doubling}\/ if $\mu$
is doubling \cite{MR1800917}. Note that the doubling property is imposed only on the measure of balls centred in $\supp\mu$.  {The following four theorems specialise Theorems~\ref{T1}\;--\;\ref{T4} to metric spaces.}

\begin{theorem}\label{T1MS}
Let $\mu$ be a doubling Borel regular probability measure on a metric space $X$. Let $\{E_i\}_{i\in\N}$ be a sequence of $\mu$-measurable subsets of $X$. Suppose that \eqref{eqn01} and \eqref{eqn02} hold for some $C>0$. In addition,
suppose that
\begin{itemize}
  \item[]
\begin{itemize}
  \item[{\bf(B1)}] for any $\delta>0$ and any closed ball $B$ centred at $\supp\mu$ there exists $i_0 =i_0 (\delta, B) $ such that for all $i  \geq i_0$  for all large enough $i$
\begin{equation}\label{vb89}
\mu\left(B\cap E_i\right)\le (1+\delta)\mu\left(B\right)\mu(E_i)\,.
\end{equation}
\end{itemize}
\end{itemize}
Then $\mu(E_\infty)=1$.
\end{theorem}

\medskip

\begin{theorem}\label{T2MS}
Let $\mu$ be a doubling Borel regular probability measure on a metric space $X$. Let $\{E_i\}_{i\in\N}$ be a sequence of $\mu$-measurable subsets of $X$. Suppose that there exist constants $C'>0$ and $\aaa>0$ and {a sequence of finite subsets $\cS_k\subset\Z$ such that $\min\cS_k\to+\infty$ as $k\to\infty$, and such that} \eqref{eqn04} and  \eqref{eqn05} {hold for all sufficiently large $k\in\N$}. In addition, suppose that  {\bf(B1)} holds. Then
$\mu(E_\infty)=1$.
\end{theorem}

\medskip

\begin{theorem}\label{T3MS}
Let $\mu$ be a doubling Borel regular probability measure on a metric space $X$. Let $\{E_i\}_{i\in\N}$ be a sequence of $\mu$-measurable subsets of $X$. Suppose that {there exists a constant $C'>0$ such that} for any $\varepsilon>0$ there exist $0<\varepsilon^*<\varepsilon$ and {a sequence of finite subsets $\cSke\subset\Z$ such that $\min\cSke\to+\infty$ as $k\to\infty$, and} such that
\eqref{eq07} and \eqref{eqn05ve} {hold for all sufficiently large $k\in\N$.} In addition,  suppose that
\begin{itemize}
  \item[]
\begin{itemize}
  \item[{\bf(B2)}] for any $\delta>0$ and any closed ball $B$ centred at $\supp\mu$ there exists $k_0 =k_0 (\delta, B) $ such that for all $k  \geq k_0$
\begin{equation}\label{vb88}
\mu\left(B\cap \bigcup_{i\in\cSke} E_i\right)\le \sum_{i\in\cSke}(1+\delta)\mu\left(B\right)\mu(E_i)\,  .
\end{equation}
\end{itemize}
\end{itemize}
Then $\mu(E_\infty)=1$.
\end{theorem}

\medskip

\begin{theorem}\label{T4MS}
Let $\mu$ be a doubling Borel regular probability measure on a metric space $X$. Let $\{E_i\}_{i\in\N}$ be a sequence of $\mu$-measurable subsets of $X$.
Suppose that there exists $0<c_0<1$ such that
\begin{itemize}
  \item[]
\begin{itemize}
  \item[{\bf(B1-$\cI$)}] for any $\delta>0$ there is an infinite subset $\cI\subset\N$ satisfying \eqref{eq19} such that for any closed ball $B$ centred at $\supp\mu$ there exists $i_0 =i_0 (\delta, B) $ such that \eqref{vb89} holds for all $i \in \cI$ with $i  \geq i_0$.
\end{itemize}
\end{itemize}
Then $\mu(E_\infty)=1$. In particular, $\mu(E_\infty)=1$ if\/ $\limsup\limits_{i\to\infty}\mu(E_i)>0$ and {\bf(B1)} is fulfilled.
\end{theorem}

\bigskip

\begin{rem}
Within the context of Theorems~\ref{T1MS}\;--\;\ref{T4MS}, the fact that $\mu$ is doubling  is assumed for convenience and aesthetics. This assumption can be  replaced with the weaker but more technical property \textbf{(P)} appearing as the conclusion  of Lemma~\ref{lem3} appearing in \S\ref{sec3.5} below. Indeed, it is this property that is at the heart of establishing the aforementioned theorems rather than the fact that $\mu$ is doubling.  It is easily seen that within the context of Lemma~\ref{lem3},  the assumption  that $\mu$ is doubling is only ever used in the proof to exploit the Vitali covering theorem and it is this that eventually leads to property \textbf{(P)}.
The upshot is that if $\mu$ is Borel regular then
$$
\text{doubling  \quad $\Longrightarrow$  \quad Vitali covering \quad $\Longrightarrow$ \quad property \textbf{(P)}  . }
$$
Thus, as an interim the doubling property can be replaced by the assumption that $(X,\mu)$ is a Vitali space, see \cite[p.6]{MR1800917}. As an example of the latter, if $\mu$ is a Radon measure on $\R^n$, then $(\R^n,\mu)$ is a Vitali space, see \cite[p.7]{MR1800917}.
\end{rem}

\medskip

\begin{rem}
When working in the setting of a metric space, a natural question to ask is:
\begin{itemize}
  \item[]{\em  What is the Hausdorff dimension/measure of $E_{\infty}=\limsup\limits_{i\to\infty}E_i$ if the measure sum $\sum_i\mu(E_i)$ converges and so $\mu\big(E_{\infty}\big)=0$?}
\end{itemize}
Indeed,  the shrinking target problem  associated with measure-preserving dynamical systems  discussed in \S\ref{STlove} below is very much in line with this type of questioning.  Within the context of this section,  we note that for various naturally arising sequences $\{E_i\}_{i\in\N}$ of $\mu$-measurable subsets of $X$ (e.g. sequences of balls) and measures supported on $X$  (e.g. Ahlfors regular measures, see \S\ref{stball}), the full $\mu$-measure statements of  Theorems~\ref{T1MS}\;--\;\ref{T4MS} can be combined with the so-called Mass Transference Principle introduced in \cite{MTP} and developed in numerous subsequent works, to obtain lower bounds for the Hausdorff dimension/measure of $E_\infty$. The complementary upper bounds can usually  be obtained using a Hausdorff measure version of the convergence  Borel--Cantelli Lemma, ultimately resulting in computing the exact Hausdorff dimension/measure of $E_{\infty}$. In particular, it is worth emphasising that the Mass Transference Principle is applicable within the context of the Duffin--Schaeffer conjecture and its inhomogeneous versions considered in \S\ref{DSBack} and \S\ref{DSNewres} below. We refer the reader to \cite{MTP} for details regarding such applications  and to  \cite{MTP-survey}  for a recent survey of developments of the Mass Transference Principle. 
\end{rem}

\section{Proofs of results on Borel--Cantelli}\label{sec3}

\subsection{Proof of Theorem~\ref{T3}}

{The basic idea of our proof originates from the work of Duffin $\&$ Schaeffer \cite{MR4859} and is different to the standard proof of Lemma~DBC, see \cite{MR1672558, MR1764799, MR1265493, MR548467}, and the similarly designed proof of Lemma~GDBC which is provided in the Appendix. Indeed, Lemma~DBC uses the Cauchy-Schwarz inequality, while our approach is based on an `inclusion-exclusion argument'.}

To begin with, suppose for a contradiction that $\mu(E_{\infty})<1$. Let $\delta>0$. Recall that
$$
E_\infty=\bigcap_{q_1=1}^\infty\;\bigcup_{i\ge q_1}E_i\,.
$$
Then, since $\mu$ is a finite measure, by the continuity of $\mu$, there exists $q_1\in\N$ large enough such that
\begin{equation}\label{eq03+}
\mu(E_\infty) \;\le\;  \mu\left(\bigcup_{i=q_1}^{\infty}E_i\right) \;\le\; \mu(E_\infty)+\delta\,.
\end{equation}
By the continuity of $\mu$ again, there exists $q_2>q_1$ large enough such that the set $A=\bigcup_{q_1\le i\le q_2}E_i$ defined within \eqref{eqn03} satisfies
\begin{equation}\label{eq01+}
  \mu(A)>\mu(E_\infty)-\delta\,.
\end{equation}
Note that the choice of $A$ depends on $\delta$ and that
\begin{equation}\label{vb7}
\text{$\mu(A)\to\mu(E_\infty)$ \quad as \quad$\delta\to0$.}
\end{equation}
{Let $C'$ be as in Theorem~\ref{T3} and} choose
\begin{equation}\label{vb008}
0<\varepsilon<\frac{1-\mu(E_\infty)}{4C'}\,.
\end{equation}
Let $\varepsilon^*$, {$\cSke$} and $k_0$ satisfy the conditions of Theorem~\ref{T3}. In particular,
we have that \eqref{eq07}, \eqref{eqn05ve} and \eqref{eq08'} hold for {all sufficiently large} $k\ge k_0$, {and $\min\cSke\to+\infty$ as $k\to\infty$}.

By \eqref{eq07} and the compactness of the interval $[\varepsilon^*,\varepsilon]$, there exists $\varepsilon'\in[\varepsilon^*,\varepsilon]$ and a sequence {$k_\ell\in\N$} such that
\begin{equation}\label{vb007+}
\sum_{i\in\cSkle}\mu(E_i) \to \varepsilon'\qquad\text{as}\qquad \ell\to\infty\,.
\end{equation}
By \eqref{vb008} and the fact that $0<\varepsilon^*\le \varepsilon'\le \varepsilon$, we have that
\begin{equation}\label{vb008+}
\varepsilon'=\frac{1-\mu(E_\infty)}{4\tilde C}\qquad\text{for some \;$\tilde C> C'$}\,.
\end{equation}
Note that the choice of $\tilde C$ is independent of $\delta$ and only depends on $\varepsilon'$ {and the sequence $k_\ell$,} and ultimately on $\ve$. Define
$$
U^\ell:=\bigcup_{i\in\cSkle} E_i\,.
$$
Then
\begin{equation}\label{eq02+}
\mu(U^\ell)\ge \sum_{i\in \cSkle} \mu\left(E_i\right)-
\sum_{\substack{s<t\\[0.5ex] s,t\in\cSkle}} \mu\big(E_s\cap E_t \big)
\;\stackrel{\eqref{eqn05ve}}{\ge} \;
\sum_{i\in \cSkle} \mu\left(E_i\right) -C^+\left(\sum_{i\in \cSkle} \mu\left(E_i\right)\right)^2
\end{equation}
{for all sufficiently large $\ell$}, where $C^+$ is any number bigger or equal to $C'$.
Next, by the definition of $A$ and $U^\ell$ {and the fact that $\min\cSke\to+\infty$ as $k\to\infty$}, we trivially have that
$$
A\cup U^\ell \subset \bigcup_{i=q_1}^{\infty}E_i
$$
for all sufficiently large $\ell$ {(how large depends on $\delta$ as well as $\ve$)} and therefore, by the right hand side of \eqref{eq03+}, we get that
$$
\mu(E_\infty)\ge\mu(A\cup U^\ell)-\delta=\mu(A)+\mu(U^\ell)-\mu(A\cap U^\ell)-\delta\,.
$$
Hence, using  \eqref{eq08'} and  \eqref{eq02+} we obtain that
$$
\mu(E_\infty)\ge \mu(A)+\sum_{i\in \cSkle} \mu\left(E_i\right) -C^+\,\left(\sum_{i\in \cSkle} \mu\left(E_i\right)\right)^2
-(1+\delta)\mu(A)\sum_{i\in \cSkle} \mu\left(E_i\right)
-\delta
$$
\begin{equation}\label{eq06+}
=\mu(A)-\delta+T_\ell b-C^+\,T_\ell^2\hspace*{37ex}
\end{equation}
{for all sufficiently large $\ell$,} where
\begin{equation}\label{vb79}
b:=1-(1+\delta)\mu(A)\qquad\text{and}\qquad T_\ell:=\sum_{i\in \cSkle} \mu\left(E_i\right)\,.
\end{equation}
By \eqref{eq03+} and the assumption that $\mu(E_\infty)<1$, we have that
\begin{equation}\label{vb009}
\frac{1-\mu(E_\infty)}{2}<b<1
\end{equation}
provided that $\delta$ is sufficiently small, in particular if $6\delta<1-\mu(E_\infty)$.
Now let $C^+:=b/2\varepsilon'$. By \eqref{vb008+} and \eqref{vb009}, we have that
\begin{equation}\label{vb011}
C'<\tilde C\le C^+\le C'':=\frac{4\tilde C}{1-\mu(E_\infty)}\,.
\end{equation}
Then, \eqref{eq06+} is applicable {for the above choice of $C^+$ and all sufficiently large $\ell$}. By \eqref{vb007+}, we get that $T_\ell\to \varepsilon'=b/2C^+$ as $\ell\to\infty$. And since the rest of the parameters in \eqref{eq06+} do not depend on $\ell$, on taking the limit as $\ell\to\infty$ we obtain from \eqref{eq06+} that
\begin{align}
\mu(E_\infty) & \ge \mu(A)-\delta+(b/2C^+)b-C^+\,(b/2C^+)^2\nonumber\\[1ex]
& = \mu(A)-\delta+b^2/4C^+\nonumber\\[1ex]
&\stackrel{\eqref{vb79}}{=} \mu(A)-\delta+\frac{1}{4C^+}\big(1-(1+\delta)\mu(A)\big)^2\nonumber\\[1ex]
&\stackrel{\eqref{vb011}}{\ge}~\mu(A)-\delta+\frac{1-\mu(E_\infty)}{16\tilde C}\big(1-(1+\delta)\mu(A)\big)^2\label{vb010}
\end{align}
for all sufficiently small $\delta>0$. Recall that $\tilde C$, {which is defined by \eqref{vb008+},} is not dependent on $\delta$.
Then, on letting $\delta\to0$ and using \eqref{vb7}, we get from \eqref{vb010} that
$$
\mu(E_\infty)\ge \mu(E_\infty)+\frac{1}{16\tilde C}\big(1-\mu(E_\infty)\big)^3\,.
$$
Hence $\big(1-\mu(E_\infty)\big)^3\le0$. However, this is impossible due to our assumption that $0\le \mu(E_\infty)<1$. Therefore, this assumption is false and so we must have that $\mu(E_\infty)=1$.

\subsection{Proof of Theorem~\ref{T4}}

Assume for a contradiction that $\mu(E_\infty)<1$. Let $0<c_0<1$ be as in the statement of the theorem. 
Let $\delta>0$ satisfy
\begin{equation}\label{delta1}
0<\delta<\frac{c_0}6(1-\mu(E_\infty))\,.
\end{equation}
By the continuity of $\mu$, there exist $q_1<q_2$ large enough such that \eqref{eq03+} and \eqref{eq01+} hold. Let $\cI$ be an infinite subset of $\N$ that meats {\bf(M1-$\cI$)}. Its existence is an assumption of the theorem.
Then, using \eqref{eq03+}, \eqref{eq01+} and property {\bf(M1-$\cI$)} we get that for any sufficiently large $i\in\cI$
\begin{align*}
\mu(E_\infty)+\delta&\ge\mu(A\cup E_{i})=\mu(A)+\mu(E_{i})-\mu(A\cap E_{i})\\[1ex]
&\ge\mu(A)+\mu(E_{i})-(1+\delta)\mu\left(A\right)\mu(E_{i})\\[1ex]
&\ge\mu(E_\infty)-\delta+\mu(E_{i})-(1+\delta)\left(\mu(E_\infty)+\delta\right)\mu(E_{i})\,.
\end{align*}
By re-arranging the terms we get that
$$
2\delta  \ge \mu(E_{i})(1-(1+\delta)\left(\mu(E_\infty)+\delta\right))\,,
$$
or equivalently
$$
2\delta +\delta\mu(E_i) \ge
\mu(E_{i})(1+\delta)(1-\mu(E_\infty)-\delta)\,.
$$
Hence, since $0\le \mu(E_i)\le 1$ and $1+\delta>1$, we get that
\begin{equation}\label{vbb}
3\delta\ge \mu(E_{i})(1-\mu(E_\infty)-\delta)\stackrel{\eqref{eq19} }{>} c_0(1-\mu(E_\infty)-\delta)\,.
\end{equation}
By \eqref{delta1}, $0<\delta<\frac{1}2(1-\mu(E_\infty))$, and thus \eqref{vbb} implies that
$3\delta > \frac{c_0}{2}(1-\mu(E_\infty))$, which contradicts \eqref{delta1}. Hence, $\mu(E_\infty)$ must be $1$.

{Recall, that the `In particular' part was already considered in Remark~\ref{rem2}. The proof of Theorem~\ref{T4} is thus complete.}

\subsection{Proof of Theorem~\ref{T2}}

{Since we have already verified Theorems~\ref{T3} and \ref{T4}, it suffices to prove that
\begin{equation}\label{T3T2}
  \text{Theorem~\ref{T3}\;\;\&\;Theorem~\ref{T4}\;\;\;$\Rightarrow$\;\;Theorem~\ref{T2}.}
\end{equation}}We will need the following auxiliary lemma, in which an empty sum is set to be $0$.

\begin{lemma}[Reduction lemma]\label{RedLem}
Let $\cSl\subset\N$ be a non-empty finite subset, $C'>0$ and
\begin{equation}\label{new1}
\sum_{\substack{s<t\\[0.5ex] s,t\in\cSl}} \mu\big(E_s\cap E_t \big) \ \le \  C'\,  \left(\sum_{i\in\cSl}\mu(E_i)\right)^2\,.
   \end{equation}
Then there exists $m\in\cSl$ such that
\begin{equation}\label{eqn09}
\sum_{\substack{s<t\\[0.5ex] s,t\in\cSlm}} \mu\big(E_s\cap E_t \big) \ \le \  C'\,  \left(\sum_{i\in\cSlm}\mu(E_i)\right)^2\,,
   \end{equation}
   where $\cSlm=\cSl\setminus\{m\}$.
\end{lemma}

\begin{proof}
If $\cSl$ has less than 3 elements, then the left hand side of \eqref{eqn09} is zero and there is nothing to prove: $m$ can be chosen to be any element of $\cSl$. Thus we can assume that $\cSl$ has at least $3$ elements.
Without loss of generality we will assume equality in \eqref{new1} as otherwise we can simply choose a smaller value for $C'$. Next, observe the following two trivial equations, which are valid for any $m\in\cS$\;:
\begin{equation*}
\sum_{\substack{s<t\\[0.5ex] s,t\in\cSl}} \mu\big(E_s\cap E_t \big) \ =
\sum_{\substack{s<t\\[0.5ex] s,t\in\cSlm}} \mu\big(E_s\cap E_t \big)+
\sum_{\substack{s\in\cSlm}} \mu\big(E_s\cap E_m \big)\,,
   \end{equation*}
$$
C'\left(\sum_{i\in\cSl}\mu(E_i)\right)^2=
C'\left(\sum_{i\in\cSlm}\mu(E_i)\right)^2+C'\mu(E_m)^2+2C'\mu(E_m)\sum_{i\in\cSlm}\mu(E_i)\,.
$$
Since, by our assumption, the left hand sides of these equations are equal, the required result will follow if for some $m\in\cSl$ we have that
\begin{equation}\label{new3}
\sum_{\substack{s\in\cSlm}} \mu\big(E_s\cap E_m \big) \ge C'\mu(E_m)^2+2C'\mu(E_m)\sum_{i\in\cSlm}\mu(E_i)\,.
\end{equation}
Assume the contrary. Then for all $m\in\cSl$
\begin{equation}\label{eqn12}
\sum_{\substack{s\in\cSlm}} \mu\big(E_s\cap E_m \big) < C'\mu(E_m)^2+2C'\mu(E_m)\sum_{i\in\cSlm}\mu(E_i)\,.
\end{equation}
Once again using the assumption that \eqref{new1} is an equation, we note that
\begin{equation}\label{new2}
\sum_{m\in\cS}\;\sum_{\substack{s\in\cSlm}} \mu\big(E_s\cap E_m \big) = 2\sum_{\substack{s<t\\[0.5ex] s,t\in\cSl}} \mu\big(E_s\cap E_t \big) \ = \  2C'\,  \left(\sum_{i\in\cSl}\mu(E_i)\right)^2\,.
\end{equation}
Then on summing \eqref{eqn12} over $m\in\cSl$ and using \eqref{new2} we get that
$$
2C'\,  \left(\sum_{i\in\cSl}\mu(E_i)\right)^2  \ < C'\sum_{m\in\cS}\mu(E_m)^2+2C'\sum_{m\in\cS}\;\sum_{i\in\cSlm}\mu(E_m)\mu(E_i)
$$
$$
\le 2C'\left(\sum_{m\in\cS}\mu(E_m)^2+\sum_{m\in\cS}\sum_{i\in\cSlm}\mu(E_m)\mu(E_i)\right)\;=\;2C'\,  \left(\sum_{i\in\cSl}\mu(E_i)\right)^2 \,,
$$
which is contradictory. Thus, \eqref{new3} must hold for at least one value of $m$.
\end{proof}

\begin{proof}[Proof of \eqref{T3T2}]
In view of Theorem~\ref{T4}, without loss of generality we will assume that $\mu(E_i)\to0$ as $i\to\infty$. {Let $C'$, $\aaa$ and $\cS_k$ satisfy the conditions of Theorem~\ref{T2}.} Fix any $\ve>0$ and define $\varepsilon^*=\min\{\aaa,\ve/2\}$. Since $\mu(E_i)\to0$ as $i\to\infty$ {and $\min\cS_k\to+\infty$ as $k\to\infty$, we have} that $\mu(E_i)<\varepsilon^*$ for all {$i\in\cS_k$ for all sufficiently large $k$.} Then, {by} \eqref{eqn04} and Lemma~\ref{RedLem}, {there exists a sequence of subsets $\cSke\subset\cS_k$} satisfying both \eqref{eq07} and \eqref{eqn05ve} {for all sufficiently large $k$}. Finally, {\bf(M2)} is satisfied as a consequence of {\bf(M1)}. Thus Theorem~\ref{T3} is applicable and we conclude that $\mu(E_\infty)=1$, as required.
\end{proof}

\subsection{Proof of Theorem~\ref{T1}}

Since we have already verified Theorem~\ref{T2}, it suffices to prove that
\begin{equation}\label{T2T1}
  \text{Theorem~\ref{T2}\;\;\;$\Rightarrow$\;\;Theorem~\ref{T1}.}
\end{equation}
Suppose \eqref{eqn01} and \eqref{eqn02} and {\bf(M1)} are satisfied. We need to identify  {a sequence of subsets $\cS_k\subset\Z$ such that $\min\cS_k\to+\infty$ as $k\to\infty$, and such that \eqref{eqn04} and \eqref{eqn05} hold for all sufficiently large $k\in\N$.} We let $\aaa=1$ and $C'=4C$. By \eqref{eqn01}, for any $k\in\N$ there is a sufficiently large $Q=Q_k$ such that
\begin{equation}\label{eqn01++}
\sum_{i=k}^{Q_k} \mu(E_i)>\max\left\{1,\sum_{i=1}^{k-1} \mu(E_i)\right\}\,.
\end{equation}
Then, by \eqref{eqn02}, there exists a constant $C>0$ such that
\begin{equation}\label{eqn02++}
\sum_{s,t=k}^{Q_k}  \mu(E_s\cap E_t)\le \sum_{s,t=1}^{Q_k}  \mu(E_s\cap E_t)\le C\left(\sum_{s=1}^{Q_k}  \mu(E_s)\right)^2
\stackrel{\eqref{eqn01++}}{\le} C\left(2\sum_{s=k}^{Q_k}  \mu(E_s)\right)^2
\end{equation}
for some increasing sequence of natural numbers $Q_k$.
Let $\cS_k=\{k,\dots,Q_k\}$. Then, {we obviously have that $\min\cS_k\to+\infty$ as $k\to\infty$, while} \eqref{eqn01++} and \eqref{eqn02++} verify \eqref{eqn04} and \eqref{eqn05}, and so Theorem~\ref{T2} is applicable and the proof is complete.

{

\subsection{Preliminaries to the proofs of Theorems~\ref{T1MS} -- \ref{T4MS}}\label{sec3.5}

For the rest of the  section $\mu$ is a doubling Borel regular probability measure on a metric space $X$.
The proofs of Theorems~\ref{T1MS} and \ref{T2MS} will follow the strategy of the proofs of Theorems~\ref{T1} and \ref{T2}. Indeed, Theorem~\ref{T1MS} will follow from Theorems~\ref{T2MS} \& \ref{T4MS}, while Theorem~\ref{T2MS} will follow from Theorem~\ref{T3MS} and the fact that $\textbf{(B1)}\Rightarrow\textbf{(B2)}$.
Theorem~\ref{T1MS} will be deduced from Theorem~\ref{T3} on demonstrating that
\begin{equation}\label{B2M2}
  \hspace*{1ex}\textbf{(B2)}\;\Rightarrow\;  \textbf{(M2)}\,.
\end{equation}
The implication $\textbf{(B1)}\;\Rightarrow\;  \textbf{(M1)}$ is generally not true. However,
the implication
\begin{equation}\label{B1IM1I}
\text{$\textbf{(B1-$\cI$)}\;\Rightarrow\;  \textbf{(M1-$\cI$)}$}
\end{equation}
does hold under assumption \eqref{eq19}, and thus we are able to deduce Theorem~\ref{T4MS} from Theorem~\ref{T4}.

Demonstrating \eqref{B2M2} and \eqref{B1IM1I} rests on approximating (in measure) the sets $A$ defined within \textbf{(M2)} with a finite collection of disjoint balls. Ultimately, this rests on the following auxiliary, and rather standard, statement from measure theory. The  proof is provided for completeness. In what follows, given two subsets $A,B\subset X$, $A\Delta B=(A\setminus B)\cup(B\setminus A)$ denotes their symmetric difference.

\begin{lemma}\label{lem3}
Let $\mu$ be a doubling Borel regular probability measure on a metric space $X$, and $A$ be a $\mu$-measurable subset of $X$. Then
\begin{enumerate}
  \item[{\rm\textbf{(P)}}]
for any $\delta'>0$ there exists a finite collection $\cG'$ of disjoint closed balls centred at $\supp\mu$ such that
\begin{equation}\label{vb55}
\mu\left(A  \, \Delta \bigcup_{B\in\cG'}B\right)\le\delta'\,.
\end{equation}
\end{enumerate}
\end{lemma}

\begin{proof}
Let $A$ be any $\mu$-measurable set and $\delta'>0$. Since $\mu$ is finite and Borel regular, we have that
$$
\mu(A)=\inf\{\mu(U):A\subset U,\;U\text{ is open}\}\,,
$$
see \cite[\S2.2.3]{Federer-69:MR0257325}. Thus, there exists an open set $U\supset A$ such that $\mu(U\setminus A)\le\delta'/2$. Let $\cF$ be the collection of all closed balls $B\subset U$ centred at $\supp\mu$. Since $\mu$ is doubling, by the Vitali covering theorem \cite[Theorem~1.6]{MR1800917}, there is a disjoint subcolection $\cG\subset\cF$ which covers $\mu$-almost all of $U$, that is $\mu\left(U\setminus\bigcup_{B\in\cG}B\right)=0$.
Since $\cG$ is disjoint and $\mu$ is finite, $\cG$ can contain only finitely many balls of measure in $[2^{-k},2^{-k+1}]$ for any $k\in\N$. Since balls in $\cG$ are centred at $\supp\mu$, $\mu(B)>0$ for any $B\in\cG$ and thus $\cG=\bigcup_{k\in\N}\{B\in\cG:\mu(B)\in[2^{-k},2^{-k+1}]\}$ is a countable union of finite collections. Thus $\cG$ is countable and we have that the sum
$$
\sum_{B\in\cG}\mu\left(B\right)=\mu\left(\bigcup_{B\in\cG}B\right)\le1
$$
is convergent. Then there exists a finite subcollection $\cG'\subset \cG$ such that
$$
\mu\left(\bigcup_{B\in\cG\setminus\cG'}B\right)=\sum_{B\in\cG\setminus\cG'}\mu\left(B\right)\le \delta'/2\,.
$$
Recall that $A\subset U$, $\bigcup_{B\in\cG}B\subset U$, $\mu\left(U\setminus\bigcup_{B\in\cG}B\right)=0$ and $\mu(U\setminus A)<\delta'/2$. Then
$$
\mu\left(A\Delta \bigcup_{B\in\cG'}B\right)\le \mu(U\setminus A)+\mu\left(\bigcup_{B\in\cG\setminus\cG'}B\right)\le \delta'/2+\delta'/2=\delta'\,.
$$
\end{proof}

\subsection{Proof of Theorem~\ref{T3MS}}

Recall that we only need to verify \eqref{B2M2}, that is $\textbf{(B2)}\Rightarrow\textbf{(M2)}$. Theorem~\ref{T3MS} will then immediately follow from Theorem~\ref{T3}. With this in mind, fix any $\ve>0$ and let $\ve^*$ and $\cSke$ be as in Theorem~\ref{T3MS}. Fix any $\delta>0$, and any natural numbers $q_1<q_2$. Let $A$ be defined by \eqref{eq08'}. Clearly, $A$ is $\mu$-measurable. Without  loss of generality we can assume that $\mu(A)>0$. Choose $\delta'>0$ small enough so that
\begin{equation}\label{vb701}
(1+\delta')\left(1+\frac{\delta'}{\mu\left(A\right)}\right)+\frac{\delta'} ok{\mu(A)\ve^*}\le 1+\delta\,.
\end{equation}
By Lemma~\ref{lem3},
there exists a finite collection of disjoint closed balls $\cG'$ centred at $\supp\mu$ satisfying \eqref{vb55}.
Observe that for any measurable subsets $B,C\subset X$
\begin{equation}\label{vb703}
\mu(A\cap C)\le \mu(B\cap C)+\mu(A\Delta B).
\end{equation}
Then, on taking $B=\bigcup_{B\in\cG'}B$ and $C=\bigcup_{i\in\cSke} E_i$, for sufficiently large $k$ we get from the above inequality and \eqref{vb55} that
\begin{align*}
\mu\left(A\cap \bigcup_{i\in\cSke} E_i\right) &\le \mu\left(\bigcup_{B\in\cG'}B\cap \bigcup_{i\in\cSke} E_i\right)+\delta'=\hspace*{7ex}\text{(by the disjointness of $\cG'$)}\\[0ex]
&=\sum_{B\in\cG'}\mu\left(B\cap \bigcup_{i\in\cSke} E_i\right)+\delta'\\[0ex]
&\stackrel{{\rm\textbf{(B2)}}}{\le} \sum_{B\in\cG'}(1+\delta')\mu(B)\sum_{i\in\cSke} \mu(E_i)+\delta'\\[0ex]
&= (1+\delta')\mu\left(\bigcup_{B\in\cG'}B\right)\sum_{i\in\cSke} \mu(E_i)+\delta'\\[0ex]
&\stackrel{\eqref{eq07}\&\eqref{vb55}}{\le} (1+\delta')\left(\mu(A)+\delta'\right)\sum_{i\in\cSke} \mu(E_i)+\frac{\delta'}{\mu(A)\ve^*}\;\mu(A)\sum_{i\in\cSke} \mu(E_i)\\[0ex]
&= (1+\delta')\left(1+\frac{\delta'}{\mu\left(A\right)}\right)\mu\left(A\right)\sum_{i\in\cSke} \mu(E_i)+\frac{\delta'}{\mu(A)\ve^*}\;\mu(A)\sum_{i\in\cSke} \mu(E_i)\\[0ex]
&\stackrel{\eqref{vb701}}{\le} (1+\delta)\mu\left(A\right)\sum_{i\in\cSke} \mu(E_i)
\end{align*}
for all sufficiently large $k$. This verifies  {\bf(M2)} and completes the proof.

\subsection{Proof of Theorem~\ref{T4MS}}

Recall that we only need to verify \eqref{B1IM1I}, that is $\textbf{(B1-$\cI$)}\Rightarrow\textbf{(M1-$\cI$)}$. Theorem~\ref{T4MS} will then immediately follow from Theorem~\ref{T4}. With this in mind, let $c_0$ be as in Theorem~\ref{T4MS}. Fix any $\delta>0$, and let $\cI$ be as  \textbf{(B1-$\cI$)}, in particular \eqref{eq19} is satisfied.
Fix any natural numbers $q_1<q_2$. Let $A$ be defined by \eqref{eqn03}, and assume without loss of generality that $\mu(A)>0$. Choose $\delta'>0$ small enough so that
\begin{equation}\label{vb702}
(1+\delta')\left(1+\frac{\delta'}{\mu\left(A\right)}\right)+\frac{\delta'}{\mu(A)c_0}\le 1+\delta\,.
\end{equation}
Since $A$ is $\mu$-measurable, by Lemma~\ref{lem3},
there exists a finite collection of disjoint closed balls $\cG'$ centred at $\supp\mu$ satisfying \eqref{vb55}.
Then, using \eqref{vb703} and the disjointness of $\cG'$, for all sufficiently large $i$ we get that
\begin{align*}
\mu\left(A\cap E_i\right)&\le \mu\left(\bigcup_{B\in\cG'}B\cap E_i\right)+\delta'\\[2ex]
&=\sum_{B\in\cG'}\mu\left(B\cap E_i\right)+\delta'\\[2ex]
&\stackrel{\textbf{(B1-$\cI$)}}{\le} \ \sum_{B\in\cG'}(1+\delta')\mu(B)\mu(E_i)+\delta'\\[2ex]
&=(1+\delta')\mu\left(\sum_{B\in\cG'}B\right)\mu(E_i)+\delta'\\[2ex]
&\stackrel{\eqref{vb55}}{\le}(1+\delta')\left(\mu\left(A\right)+\delta'\right)\mu(E_i)+\delta' \\[2ex]
&\le(1+\delta')\left(1+\frac{\delta'}{\mu\left(A\right)}\right)\mu\left(A\right)\mu(E_i)+\frac{\delta'}{\mu(A)c_0}\mu(A)\mu(E_i)\\[2ex]
&\stackrel{\eqref{vb702}}{\le}(1+\delta)\mu\left(A\right)\mu(E_i)\,.
\end{align*}
This verifies {\bf(M1-$\cI$)} and completes the proof.

\subsection{Proof of Theorems~\ref{T1MS} and \ref{T2MS}}

\begin{proof}[Proof of Theorem~\ref{T2MS}]
In view of Theorem~\ref{T4MS}, without loss of generality we will assume that $\mu(E_i)\to0$ as $i\to\infty$. Let $C'$, $\aaa$ and $\cS_k$ satisfy the conditions of Theorem~\ref{T2MS}. Fix any $\ve>0$ and define $\varepsilon^*=\min\{\aaa,\ve/2\}$. Since $\mu(E_i)\to0$ as $i\to\infty$ and $\min\cS_k\to+\infty$ as $k\to\infty$, we have that $\mu(E_i)<\varepsilon^*$ for all $i\in\cS_k$ for all sufficiently large $k$. Then, by \eqref{eqn04} and Lemma~\ref{RedLem}, there exists a sequence of subsets $\cSke\subset\cS_k$ satisfying both \eqref{eq07} and \eqref{eqn05ve} for all sufficiently large $k$. Finally, {\bf(B2)} is satisfied as a consequence of {\bf(B1)}. Thus Theorem~\ref{T3MS} is applicable and we conclude that $\mu(E_\infty)=1$, as required.
\end{proof}

\begin{proof}[Proof of Theorem~\ref{T1MS}]
Suppose \eqref{eqn01} and \eqref{eqn02} and {\bf(B1)} are satisfied. We will appeal to Theorem~\ref{T2MS}. To this end we need to identify a sequence of subsets $\cS_k\subset\Z$ such that $\min\cS_k\to+\infty$ as $k\to\infty$, and such that \eqref{eqn04} and \eqref{eqn05} hold for all sufficiently large $k\in\N$. We let $\aaa=1$ and $C'=4C$. By \eqref{eqn01}, for any $k\in\N$ there is a sufficiently large $Q=Q_k$ such that \eqref{eqn01++} holds.
Then, by \eqref{eqn02}, there exists a constant $C>0$ such that
\eqref{eqn02++} holds for some increasing sequence of natural numbers $Q_k$.
Let $\cS_k=\{k,\dots,Q_k\}$. Then, we obviously have that $\min\cS_k\to+\infty$ as $k\to\infty$, while \eqref{eqn01++} and \eqref{eqn02++} verify \eqref{eqn04} and \eqref{eqn05}, and so Theorem~\ref{T2MS} is applicable and we conclude that $\mu(E_\infty)=1$, as required.
\end{proof}

}

\section{The Duffin--Schaeffer conjecture}\label{DSBack}

\subsection{Notation}

The notation introduced here will apply throughout Section~\ref{DSBack} and the two subsequent sections. By measure $\mu$ we will mean Lebesgue measure on the real line. Also, in what follows $\ll$ will denote the Vinogradov symbol, that is $x\ll y$ means that $x\le Cy$ for some fixed $C>0$ independent of $x$ and $y$. Furthermore, we write $x \asymp y$ whenever $x \ll y \ll x$. Given a real number $x$, $\|x\|$ will denote the distance of $x$ from the nearest integer. Given two integers $a,b$ we will write $(a,b)$ for the $\gcd(a,b)$, unless we specifically mention that $(a,b)$ means a pair of integers. The same will apply to triples of integers. For every positive integer $q$, let 
$$
\mathbb{Z}_q := \{0,\ldots,q-1\},
$$ 
thus $\Z_q$ is a complete system of residues modulo $q$. Also, we let 
$$
\Z_q^*:=\{a\in\Z:1\le a\le q,\;(a,q)=1\},
$$ 
thus $\Z_q^*\subset\Z_q$ is a reduced system of residues modulo $q$. Given a function $\psi:\N\to\R_{\ge0}$, the support of $\psi$ is defined as $\supp(\psi):=\{q\in\N:\psi(q)\neq0\}$. Given any set $S$, the function  $\mathds{1}_S$ will denote the characteristic/indicator function of $S$. Given a finite set $S$, the number of elements in $S$ will be denoted by either $|S|$ or $\#S$. Also, for an interval $I$ we will write $|I|$ for its length. Finally, 
we define 
$$
\Log(x) := \max\{1,\log(x)\}
$$
where $\log$ denotes the natural logarithm.

\subsection{Background}

To set the scene, let $\gamma\in\R$, $\psi:\N\to\R_{\ge0}$ be a real, positive  function and  define the following sequence of subsets of $[0,1]$ indexed by $q\in\N$:
\begin{align}
E_q=E_q(\gamma,\psi):&=
\left\{x \in [0,1]: \left| x - \frac{a+\gamma}{q} \right| \leq \frac{\psi(q)}{q} \text{ for some }a\in\Z\right\}\nonumber\\[2ex]
&=\left\{x \in [0,1]: \left\lVert x - \frac{a+\gamma}{q} \right\rVert \leq \frac{\psi(q)}{q} \text{ for some }a\in\Z_q\right\}\,.\label{e1}
\end{align}

\noindent For obvious reasons the function $\psi$ is often referred to as an approximating function.  In short it governs the ``rate'' at which the ``shifted  rationals'' $(a+\gamma)/q$   approximate $x$. When $\gamma=0$, we fall into the homogeneous case dealing with rational approximations $a/q$ to $x$. Indeed, in this case $x\in E_q$ if and only if $|x-a/q | <\psi(q)/q$ for some $a\in\Z$. We also note that \eqref{e1} has a natural reformulation in terms of the dynamics of circle rotations -- a fundamental theme in dynamical systems. Indeed, $x\in E_q$ if and only if any point $e^{2\pi \alpha i}$ on the unit circle rotated $q$ times by the angle $2 \pi x$ falls into the $2\pi\psi(q)$-neighborhood of $e^{2\pi (\alpha+\gamma) i}$ (measured by angle). When $\gamma=0$, these $q$ rotations simply return the point to its own $2\pi\psi(q)$-neighborhood.

We will be interested in the set
\begin{equation}\label{e46}
  \limsup_{q\to\infty} E_q\,,
\end{equation}
which consists of $x\in[0,1]$ such that $\left|qx-a-\gamma\right|<\psi(q)$ holds for infinitely many $(q,a)\in\N\times\Z$. Understanding the measure theoretic properties of \eqref{e46} has a long history, staring 100 years ago with Khintchine's work \cite{MR1512207} who proved, in the case $\gamma=0$, that for any $\psi$ such that $q\psi(q)$ is monotonically decreasing\footnote{Nowadays Khintchine's theorem is often stated with the weaker assumption that $\psi$ is non-increasing.},
\begin{equation}\label{Khi}
  \mu\Big(\limsup_{q\to\infty}E_q\Big)=\left\{\begin{array}{cl}
                         0 & \text{if }\sum_{q=1}^\infty\mu(E_q)<\infty\,,\\[2ex]
                         1 & \text{if }\sum_{q=1}^\infty\mu(E_q)=\infty\,.
                       \end{array}\right.
\end{equation}
Note that if $\psi(q)> 1/2$ then $E_q=[0,1]$ and thus $\limsup_{q\to\infty}E_q=[0,1]$ if $\psi(q)> 1/2$ for infinitely many $q\in\N$. Thus, in the context of Khintchine's theorem one is really only interested in $\psi$ bounded by $1/2$. It is readily seen that if $\psi(q)\le1/2$ then 
\begin{equation}\label{KKK}
\mu(E_q)=2\psi(q)
\end{equation}
irrespective of $\gamma$.
We also observe that the main substance of \eqref{Khi} lies within its divergence part, since the convergence part of \eqref{Khi} is a trivial consequence of the first Borel--Cantelli Lemma and does not actually require $\psi$ being monotonic. However, the monotonicity assumption on $\psi$ cannot be dropped altogether from
the divergence part of Khintchine's theorem. This was shown in 1941 by Duffin and Schaeffer \cite{MR4859} who also conjectured an alternative statement, known as the Duffin--Schaeffer conjecture. The conjecture incorporates the additional requirement that $a$ is coprime to $q$ within \eqref{e1}. Thus the approximating rational numbers $a/q$ are written in the lowest terms and never repeated in the sequence of all possible approximating fractions. To introduce the conjecture formally, let us define the following modification of $E_q$:
\begin{align}
E'_q=E'_q(\gamma,\psi)&:=\left\{x \in [0,1]: \left| x - \frac{a+\gamma}{q} \right| \leq \frac{\psi(q)}{q} \text{ for some }a\in\Z,\;(a,q)=1\right\}\nonumber\\[2ex]
&=\left\{x\in[0,1]:\left\lVert x-\frac{a+\gamma}{q}\right\rVert \leq\frac{\psi(q)}{q}\text{ for some $a\in\mathbb{Z}_q^{*}$}\right\}\label{e1+}\,.
\end{align}
Then the \emph{Duffin--Schaeffer conjecture} claims that, assuming $\gamma=0$, for any nonnegative function $\psi$,
\begin{equation}\label{DSStatement}
  \mu\Big(\limsup_{q\to\infty}E'_q\Big)=\left\{\begin{array}{cl}
                         0 & \text{if }\sum_{q=1}^\infty\mu(E'_q)<\infty\,,\\[2ex]
                         1 & \text{if }\sum_{q=1}^\infty\mu(E'_q)=\infty\,.
                       \end{array}\right.
\end{equation}
Again assuming $\psi: \N \to [0,\tfrac{1}{2}]$, this time we have that
\begin{equation}\label{DSDS}
\mu(E'_q)=\frac{2\psi(q)\varphi(q)}{q}\,,
\end{equation}
where $\varphi$ is Euler's totient function, and the divergence sum condition in \eqref{DSStatement} becomes
\begin{equation}\label{Div-Sum}
\sum_{q=1}^\infty \frac{\varphi(q)\psi(q)}{q}=\infty\,.
\end{equation}
In fact, Duffin and Schaeffer \cite{MR4859} only ever considered the case when $\psi(q)\le 1/2$ for all $q$. Note that if $\psi(q)>1/2$, $\mu(E'_q)$ may become much smaller than the right hand side of \eqref{DSDS} and thus \eqref{Div-Sum} may seem  weaker than the divergence sum condition in \eqref{DSStatement}. The truth is that the two divergence sum conditions are equivalent. This follows from a theorem of Pollington and Vaughan \cite[Theorem~2]{pv}, stated below as \eqref{PV1}, and the trivial fact that the convergence case in \eqref{DSStatement} always holds. In short, in the case $\gamma =0$,  Pollington and Vaughan proved that
\begin{equation}\label{PV1}
\forall\;\psi:\N\to\R_{\ge0}\qquad \psi\ge\tfrac12\mathds{1}_{\supp(\psi)}\;\text{ and \eqref{Div-Sum} }\Longrightarrow \mu\Big(\limsup_{q\to\infty}E'_q\Big)=1\,.
\end{equation}

Clearly, $E'_q\subset E_q$ and therefore $\limsup_{q\to\infty}E'_q\subset \limsup_{q\to\infty}E_q$. Thus, the Duffin--Schaeffer conjecture implies the following statement that is referred to as the \emph{weak Duffin--Schaeffer conjecture} \cite[Conjecture~9.1]{RandomFR}\footnote{Apart from introducing the conjecture, Ram\'irez also proved in \cite[Eq.~(22)]{RandomFR} that the weak Duffin-Schaeffer conjecture is equivalent to a famous conjecture of Catlin \cite{Catlin1976}}: assuming $\gamma=0$, we have that
\begin{equation}\label{Catl}
  \mu(\limsup_{q\to\infty}E_q)=1\qquad\text{if}\qquad \sum_{q=1}^\infty \frac{\varphi(q)\psi(q)}{q}=\infty\,.
\end{equation}
The weak Duffin-Schaeffer conjecture provides a generalisation to \eqref{Khi} by modifying the diverging sum rather than modifying the sets $E_q$. The sum in \eqref{Catl} can be interpreted as the Khintchine sum with the weights $\varphi(q)/q$.
Since $\liminf_{q\to\infty}\varphi(q)/q=0$, the weak Duffin-Schaeffer conjecture can be viewed as Khintchine's theorem without monotonicity, but with extra divergence. However, in the case of monotonic $\psi$, it is easily verified, via partial summation that the Duffin-Schaeffer sum is comparable to the Khintchine sum in the sense that
\begin{equation}\label{DS_condition}
\sum_{q=1}^Q \psi(q)\;\ll\;  \sum_{q=1}^Q \frac{\varphi(q)\psi(q)}{q} \;\ll\;\sum_{q=1}^Q \psi(q)
\end{equation}
for all $Q\in\N$.

Quite remarkably, Duffin and Schaeffer \cite{MR4859} demonstrated that it is property \eqref{DS_condition} rather than the monotonicity of $\psi$ that is instrumental in the validity of Khintchine's theorem. Specifically, they proved that, for $\gamma=0$, \eqref{Khi} and indeed \eqref{DSStatement} hold provided that \eqref{DS_condition} is satisfied for infinitely many $Q\in\N$. In the literature this finding of Duffin and Schaeffer is often quoted as the Duffin--Schaeffer theorem. This and other results towards the Duffin--Schaeffer conjecture obtained prior to the turn of the Millennium can be found in Harman's monograph \cite{MR1672558}. The progress was then largely stalled. However, within the last decade or so, rapid advances have been made starting with results on a version of the conjecture with extra divergence conditions, see \cite{MR3784747}, \cite{MR4008524}, \cite{MR3101800}, \cite{MR2915535}, and culminating with the major breakthrough of Koukoulopoulos and Maynard \cite{MR4125453} who resolved the Duffin--Schaeffer conjecture in full. Its asymptotic (quantitative) form was subsequently  obtained in \cite{abh} by Aistleitner {\em et al.}.  The error term associated with the asymptotics has most recently be improved    to essentially the best possible by Koukoulopoulos, Maynard and Yang \cite{KMY24}.

\subsection{Open questions}

Note that much of the progress discussed above was underpinned  by the zero-one laws of Cassels \cite{MR36787} and Gallagher \cite{MR133297}, stating that, for $\gamma=0$, for any nonnegative $\psi$
\begin{equation}\label{CG}
  \mu(\limsup_{q\to\infty} E_q)\in\{0,1\}\qquad\text{and}\qquad \mu(\limsup_{q\to\infty} E'_q)\in\{0,1\}\,.
\end{equation}
Unfortunately, similar results are not available for non-integer $\gamma$, see \cite{MR3606945}.
In fact, despite the aforementioned breathtaking progress, to date little is known regarding inhomogeneous forms of the Duffin--Schaeffer conjecture. The inhomogeneous version of Khintchine's theorem (for monotonic $\psi$) was established by Sz\"usz \cite{MR0095165} who proved \eqref{Khi} for arbitrary fixed $\gamma\in\R$. However, beyond the monotonic case the progress was very limited, and can be found in \cite{CT-Mem}, \cite{MR3988812} and \cite{MR4230541}. Similarly to the homogeneous case, it was shown by Ram\'irez \cite{MR3606945} that the monotonicity condition cannot be dropped altogether from Sz\"usz's inhomogeneous version of Khintchine's theorem. It means that the removal of the monotonicity assumption will require either restrictions on the integers $a$ appearing in \eqref{e1} (we formulate this as the Main Problem below), or imposing an extra divergence condition (see Conjecture~\ref{ICC} below).

\medskip

\begin{mainproblem}[Inhomogeneous Duffin--Schaeffer]\label{PA}
Let $\gamma \in \R$, $\psi:\N\to\R_{\ge0}$ and $\mathcal{I}_q$, $q=1,2,\dots$ denote a sequence of subsets of $\mathbb{Z}_q := \{0,\ldots,q-1\}$ such that
\begin{equation}\label{DivSumI_n}
  \sum_{q=1}^\infty\frac{\lvert \mathcal{I}_q\rvert \psi(q)}{q} =\infty\,.
\end{equation}
Let
\begin{equation}\label{EE}
E_q^{\mathcal{I}} = E_q^{\mathcal{I}}(\gamma,\psi) := \left\{x \in [0,1]: \left\lVert x - \frac{a+\gamma}{q} \right\rVert \leq \frac{\psi(q)}{q} \text{ for some }a \in \mathcal{I}_q\right\}.
\end{equation}
\begin{itemize}
  \item[{\rm(i)}] Under what conditions on $\cI_q$ $($possibly depending on $\gamma$ and $\psi)$ do we have that
\begin{equation}\label{DStP}
\mu(\limsup_{q \to \infty}{E_q^{\mathcal{I}}}) = 1\,?
   \end{equation}
  \item[{\rm(ii)}] Is the above always possible with 
  \begin{equation}\label{I_nLB}
  |\cI_q|\ge c_1\varphi(q) 
  \end{equation}
  for some $c_1>0$ independent of $q$?
  \item[{\rm(iii)}] In particular, what are the $\gamma$'s such that \eqref{DSStatement} holds for $E_q'$ given by \eqref{e1+}?
\end{itemize}
\end{mainproblem}

\medskip

Statement \eqref{DSStatement} was first formally conjectured for $\gamma\neq0$ by Ram\'irez \cite{MR3606945} who coined it the inhomogeneous Duffin--Schaeffer Conjecture. However, the relevance of the coprimality condition $(a,q)=1$ used in the definition of $E_q'$ to inhomogeneous approximations remains unclear. Indeed, when $\gamma\neq0$ this condition does not prevent the situation when the approximating fractions $(a+\gamma)/q$ come too close together or even coincide, for instance, in the case $\gamma\in\Z\setminus\{0\}$.
Note that in the more general setting of (i) and (ii), establishing \eqref{DStP} for some $\gamma\in\R$ and sequence $(\mathcal{I}_q(\gamma))_{q \in \mathbb{N}}$ immediately proves \eqref{DStP} for any integer translate $\gamma + m$ of $\gamma$, where $m \in \mathbb{Z}$, and the corresponding translate of $\cI_q(\gamma)$, that is $\mathcal{I}_q(\gamma+m) := \mathcal{I}_q(\gamma) - m \pmod q$. In particular, establishing (ii) for some $\gamma$ proves it for all integer translates of $\gamma$. However, the sets $\cI_q=\mathbb{Z}_q^{*}$ are not invariant under arbitrary integer translations as above. This means that in order to prove Ram\'irez's conjecture, one cannot restrict to, say $\gamma \in [0,1)$.
As a result,  Ram\'irez's conjecture (part (iii)\;) remains open even for $\gamma \in \mathbb{Z}\setminus \{0\}$. In fact, even weaker versions of Ram\'irez's conjecture based on imposing extra divergence as in \cite{MR3101800}, or condition \eqref{DS_condition} are not known and would still be of interest. In particular, the use of \eqref{DS_condition} leads to the following weaker version of Ram\'irez's conjecture:

\begin{conjecture}[Inhomogeneous Duffin--Schaeffer Theorem]\label{IDST}
Statement~\eqref{DSStatement} holds if additionally \eqref{DS_condition} is satisfied for infinitely many $Q\in\N$.
\end{conjecture}

In this paper we will also investigate the following weaker form of \eqref{DSStatement}, see \cite[Conjecture~1.22]{CT-Mem}:

\begin{conjecture}[Weak inhomogeneous Duffin-Schaeffer conjecture]\label{ICC}
Let $\gamma\in\R$, $\psi:\N\to\R_{\ge0}$ and $E_q$ be defined as in \eqref{e1}. Then \eqref{Catl} holds.
\end{conjecture}

\section{New results on Duffin--Schaeffer}\label{DSNewres}

\subsection{Inhomogeneous Duffin--Schaeffer for large approximating functions}

We begin by establishing a complete analogue of Pollington--Vaughan's result \eqref{PV1} for inhomogeneous approximations. In other words we shall prove \eqref{DSStatement} for functions $\psi$ which are $\ge1/2$ on the support. This will be implied by the following dichotomy statement applicable to the general setup of the Main Problem.

\medskip

\begin{theorem}\label{dichotomy}
Let $\psi:\N\to\R_{\ge0}$ be such that $\lim_{q\to\infty}\psi(q)/q=0$. Fix any $\delta>0$ and
suppose that $\psi(q)\ge\delta$ for any $q\in\N$ such that $\psi(q)\neq0$. Let $\cI_q\subset\Z_q$ be any sequence. Then the following two statements are equivalent:
\begin{enumerate}
  \item[$(\exists)$] \qquad$\exists\;\gamma\in\R$ \quad \eqref{DStP} holds.
  \item[$(\forall)$] \qquad$\forall\;\gamma\in\R$ \quad \eqref{DStP} holds.
\end{enumerate}
\end{theorem}

\medskip

The proof is a simple application of the following lemma which is a slight modification of Cassels' Lemma~9 in \cite{MR36787}, and can be explicitly found as Lemma~1 in \cite{MR2457266}.

\medskip

\begin{lemma}[Cassels' Lemma]\label{CasselsLemma}
Let $I_j$ be a sequence of intervals in $\R$ such that $|I_j|\to0$ as $j\to\infty$. Let $U_j\subset I_j$ be a sequence of Lebesgue measurable sets. Suppose that for some $c>0$ we have that $\mu(U_j)\ge c|I_j|$ for all $j$. Then the sets
$$
\limsup_{j\to\infty} I_j\qquad\text{and}\qquad \limsup_{j\to\infty} U_j
$$
have the same Lebesgue measure.
\end{lemma}

\begin{proof}[Proof of Theorem~\ref{dichotomy}]
We only need to verify that $(\exists)\Rightarrow(\forall)$. Suppose that $\exists$ $\gamma\in\R$ such that \eqref{DStP} holds. This means that for almost every $x\in[0,1]$ there are infinitely many $q\in\N$, $a\in\cI_q$ and $k\in\Z$ satisfying
$$
\left|x-\frac{a+\gamma}{q}-k\right|<\frac{\psi(q)}{q}\,.
$$
Fix any other $\gamma'\in\R$, and define $C=\delta^{-1}|\gamma-\gamma'|+1$. By the triangle inequality and the condition that $\psi\ge\tfrac12\mathds{1}_{\supp(\psi)}$, we get that for any $q$ such that $\psi(q)\neq0$
$$
\left|x-\frac{a+\gamma'}{q}-k\right|\le \left|x-\frac{a+\gamma}{q}-k\right|+\left|\frac{\gamma-\gamma'}{q}\right|<\frac{\psi(q)}{q}+\frac{|\gamma-\gamma'|}{q}\le
\frac{C\psi(q)}{q}\,.
$$
Since $\psi(q)/q=o(1)$, by Lemma~\ref{CasselsLemma}, for almost every $x\in[0,1]$, we also have that
$$
\left|x-\frac{a+\gamma'}{q}-k\right|<\frac{\psi(q)}{q}
$$
holds for infinitely many $q\in\N$, $a\in\cI_q$ and $k\in\N$. This means that \eqref{DStP} holds for $\gamma'$ and completes the proof.
\end{proof}

\medskip

When $\cI_q=\Z_q^*$ we know, by \eqref{PV1}, that \eqref{DStP} holds for $\gamma=0$ for any $\psi$ satisfying the assumptions of Theorem~\ref{dichotomy} and \eqref{Div-Sum}. As a result we have the following statement, which resolves part  (iii) of the Main Problem for large $\psi$.

 \medskip
 
\begin{theorem}[The inhomogeneous Duffin--Schaeffer Conjecture for large $\psi$]\label{DS_large_psi}
Let $\gamma\in\R$, $\delta>0$ and $\psi:\N\to\R_{\ge0}$ satisfy the condition $\psi(q)\ge\delta$ for any $q\in\N$ with $\psi(q)\neq0$. Suppose that \eqref{Div-Sum} holds. Then, with $E'_q$ as in \eqref{e1+}, we have that
$$
\mu\Big(\limsup_{q\to\infty}E'_q\Big)=1\,.
$$
\end{theorem}

\begin{proof}
Without loss of generality we can assume that $\psi(q)=o(q)$. If this was not the case, we would replace $\psi$ with, for instance, $\tilde\psi(q)=\min\{\psi(q),q^{1/2}\}$, which would thus satisfy $\tilde\psi(q)\ge q^{1/2}$ for infinitely many $q$. Indeed, the new function $\tilde\psi$ obviously satisfies the condition $\tilde\psi(q)\ge\delta$ for any $q\in\N$ with $\tilde\psi(q)\neq0$. Furthermore,
$$
\sum_{q=1}^\infty \frac{\varphi(q)\tilde\psi(q)}{q}=\infty\,.
$$
If the latter was not the case, then we would have that $\frac{\varphi(q)\tilde\psi(q)}{q}\to0$ as $q\to\infty$, implying that $\tilde\psi(q)<q/\varphi(q)\ll \log\log q$ for all sufficiently large $q$ and contrary to the fact that  $\tilde\psi(q)\ge q^{1/2}$ for infinitely many $q$. Here we used the well known inequality $\varphi(q)\gg q/\log(\log q)$ valid for all sufficiently large $q$. To complete the proof it remains to apply Lemma~\ref{CasselsLemma} and Theorem~\ref{dichotomy} to \eqref{PV1}.
\end{proof}

\bigskip

Note that in the above proof we only relied on \eqref{PV1}, and did not use the validity of the Duffin--Schaeffer conjecture as established by Koukoulopoulos and Maynard.

\subsection{Borel--Cantelli and zero-one laws for inhomogeneous Duffin--Schaeffer}

Recall that the homogeneous results discussed in \S\ref{DSBack} are proved by
\begin{itemize}
\item establishing that relevant limsup sets have positive measure on verifying the conditions of Lemma~DBC/GDBC, and then
\item establishing that the limsup sets have full measure on applying zero-one laws \eqref{CG}.
\end{itemize}
However, as we already mentioned, there is no analogue of \eqref{CG} in the inhomogeneous case. We will demonstrate that, in fact, there is no need for a zero-one law in the inhomogeneous case, and that verifying conditions of either Lemma~DBC or Lemma~GDBC will already be enough to conclude a full measure statement for the limsup sets in question. To this end we use the results of \S\ref{sec2} and verify that condition \textbf{(B1)} holds in the relevant setup. We begin with the following statement of this type.

\medskip

\begin{thm}\label{new_thm10}
   Let $\gamma \in \R$, $\psi:\N\to[0,\tfrac12]$ and $\mathcal{I}_q \subseteq \mathbb{Z}_q$ be a sequence such that $\tfrac{1}{q}\mathcal{I}_q$ is uniformly distributed modulo $1$, meaning that
    for any $0 \leq x \leq y \leq 1$ we have that
\begin{equation}\label{ud}
\lim_{q \to \infty}\frac{\left\lvert \left\{
a \in \mathcal{I}_q: \frac{a}{q} \in \left[x ,y\right)\right\} \right\rvert}
{\lvert \mathcal{I}_q \rvert } = y-x.
\end{equation}
Let $E_q^{\mathcal{I}}$ be given by \eqref{EE}.
Suppose that there exist constants $C'>0$ and $\aaa>0$ and a sequence of finite subsets $\cS_k\subset\Z$ such that $\min\cS_k\to+\infty$ as $k\to\infty$, satisfying
\begin{equation}\label{enough_mass}
\sum_{q\in\mathcal{S}_k}\frac{\lvert \mathcal{I}_q\rvert \psi(q)}{q} \ge \aaa
   \end{equation}
and
\begin{equation} \label{quasi_indep}
\sum_{\substack{r<q\\[0.5ex] q,r\in\mathcal{S}_k}} \mu\big(E_q^{\mathcal{I}}\cap E_r^{\mathcal{I}} \big) \ \le \  C'\,  \left(\sum_{q\in\mathcal{S}_k}\frac{\lvert \mathcal{I}_q\rvert\psi(q)}{q}\right)^2\,.
   \end{equation}
Then \eqref{DStP} holds.
\end{thm}

\medskip

Note that \eqref{enough_mass} and \eqref{quasi_indep} are exactly conditions \eqref{eqn04} and \eqref{eqn05} within Lemma~GDBC rewritten for the setup of Theorem~\ref{new_thm10}. However, unlike in Lemma~GDBC, here we obtain a full measure conclusion. One can also easily obtain the following version of Theorem~\ref{new_thm10} in the style of Lemma~DBC.

\medskip

\begin{thm}\label{new_thm10v2}
Let $\gamma \in \R$, $\psi:\N\to[0,\tfrac12]$, $\mathcal{I}_q \subseteq \mathbb{Z}_q$ be a sequence such that $\tfrac{1}{q}\mathcal{I}_q$ is uniformly distributed modulo $1$ and
\begin{equation}\label{enough_mass2}
\sum_{q=1}^\infty\frac{\lvert \mathcal{I}_q\rvert \psi(q)}{q} =\infty\,.
   \end{equation}
Let $E_q^{\mathcal{I}}$ be given by \eqref{EE}, and suppose that there exists a constant $C>0$ such that
\begin{equation} \label{quasi_indep2}
\sum_{q,r=1}^Q \mu\big(E_q^{\mathcal{I}}\cap E_r^{\mathcal{I}} \big) \ \le \  C\,  \left(\sum_{q=1}^Q\frac{\lvert \mathcal{I}_q\rvert\psi(q)}{q}\right)^2
   \end{equation}
holds for infinitely many $Q\in\N$. Then \eqref{DStP} holds.
\end{thm}

\medskip

\begin{rem}\label{rem_bounded_psi}
Clearly, Theorem~\ref{new_thm10} and Theorem~\ref{new_thm10v2} can be used to prove Conjecture~\ref{ICC} (the weak inhomogeneous Duffin-Schaeffer conjecture). Indeed, since $\limsup_{q\to \infty}{E^{\cI}_q}\subset\limsup_{q \to \infty}{E_q}$, Conjecture~\ref{ICC} will follow if the conditions of Theorem~\ref{new_thm10} or Theorem~\ref{new_thm10v2} can be met with $\cI_q$ {\rm(}\/possibly depending on $\psi$ and $\gamma${\rm)} such that
\eqref{I_nLB} holds.
In fact, by \cite[Theorem~3]{MR4497313},  for any $\psi:\N\to\R_{\ge0}$ and any $\gamma\in\R$, $\mu(\limsup_{q \to \infty}{E_q}) >0$ if and only if there is a sequence $\mathcal{I}_q \subseteq \mathbb{Z}_q$ satisfying \eqref{enough_mass2} and \eqref{quasi_indep2}. However, deriving a full measure statement would require additional assumption(s) imposed on $\cI_q$. For this purpose in Theorems~\ref{new_thm10} and \ref{new_thm10v2}, we use the rather natural assumption of uniform distribution.
\end{rem}

The following proposition contains a particular choice for the sets $\mathcal{I}_q$ satisfying \eqref{I_nLB}, which in view of Remark~\ref{rem_bounded_psi} represent a candidate for proving the weak inhomogeneous Duffin-Schaeffer conjecture.
Indeed, in \S\ref{sec5.2.2}---\ref{explicit} we present several applications based on this choice of $\cI_q$ which include a full proof of the weak inhomogeneous Duffin-Schaeffer  conjecture for arbitrary rational $\gamma$ and the Duffin--Schaeffer conjecture with congruence conditions imposed on the numerators and denominators of the approximating rational fractions.

\begin{proposition}\label{equidist_prop}
For $q\in\N$ let
\begin{equation}\label{def_Iq}
\mathcal{I}_q = \{a \in \mathbb{Z}_q: (A_q + aB_q, q) = 1\}\,,
\end{equation}
where $A_q,B_q$ are coprime integers such that $B_q>0$. Then both \eqref{ud} and \eqref{I_nLB} are satisfied, and
\begin{equation}\label{eq103}
    |\cI_q|=\varphi(q,B_q)\,,\qquad\text{where } \quad \varphi(q,b):=\varphi(q)\prod_{p \mid (q,b)} \left(1 + \frac{1}{p-1}\right)\,.
\end{equation}
Furthermore, for any $\gamma \in \mathbb{R}$ and $\psi: \N \to [0,\tfrac{1}{2}]$ the set $E_q^{\mathcal{I}}$ given by \eqref{EE} satisfies
\begin{equation}\label{sum_of_measures}
    \mu(E_q^{\mathcal{I}}) = \frac{2\psi(q)|\cI_q|}{q}\,.
\end{equation}
\end{proposition}

Applying Proposition~\ref{equidist_prop} with $A_q=0$, $B_q=1$ to Theorem~\ref{new_thm10} together with Theorem~\ref{DS_large_psi}, to deal with the case $\psi > \tfrac{1}{2}$, gives the following statement. It provides sufficient conditions for establishing \eqref{DSStatement} without  a zero-one law.
\medskip

\begin{corollary}\label{T6}
Let $\gamma\in\R$, $\psi:\N\to\R_{\ge 0}$ and $E'_q=E'_q(\gamma,\psi)$ be given by \eqref{e1+}.
Suppose that there exist constants $C'>0$ and $\aaa>0$ and a sequence of finite subsets $\cS_k\subset\Z$ such that $\min\cS_k\to+\infty$ as $k\to\infty$, satisfying
\begin{equation}\label{eqn04new}
\sum_{q\in\cS_k}\frac{\varphi(q)\psi(q)}{q} \ge \aaa
   \end{equation}
and
\begin{equation} \label{eqn05new}
\sum_{\substack{r<q\\[0.5ex] q,r\in\cS_k}} \mu\big(E'_q\cap E'_r \big) \ \le \  C'\,  \left(\sum_{q\in\cS_k}\frac{\varphi(q)\psi(q)}{q}\right)^2\,.
   \end{equation}
Then
\begin{equation}\label{e65}
  \mu(\limsup_{q\to\infty} E_q')=1\,.
\end{equation}
\end{corollary}

\medskip

\noindent Similarly, applying Proposition~\ref{equidist_prop} with $A_q=0$, $B_q=1$ to Theorem~\ref{new_thm10v2} together with Theorem~\ref{DS_large_psi} gives

\medskip

\begin{corollary}\label{T5+}
Let $\gamma\in\R$, $\psi:\N\to\R_{\ge0}$ and $E'_q=E'_q(\gamma,\psi)$ be given by \eqref{e1+}.
Suppose that \eqref{Div-Sum} holds and there exists $C>0$ such that
$$
\sum_{q,r=1}^Q  \mu(E'_{q}\cap E'_{r})\le C\left(\sum_{q=1}^Q \frac{\varphi(q)\psi(q)}{q}\right)^2
$$
for infinitely many $Q\in\N$. Then \eqref{e65} holds.
\end{corollary}

\medskip

In turn, the following obvious modification of Corollary~\ref{T5+} provides sufficient conditions for establishing Conjecture~\ref{IDST} without a zero-one law.

\begin{corollary}\label{T5}
Let $\gamma\in\R$, $\psi:\N\to \R_{\ge 0}$ and $E'_q=E'_q(\gamma,\psi)$ be given by \eqref{e1+}.
Suppose that \eqref{Div-Sum} holds and there exists $C>0$ such that
$$
\sum_{q,r=1}^Q  \mu(E'_{q}\cap E'_{r})\le C\left(\sum_{q=1}^Q \psi(q)\right)^2
$$
and \eqref{DS_condition} hold simultaneously for infinitely many $Q\in\N$. Then \eqref{e65} holds.
\end{corollary}

\subsection{Inhomogeneous Duffin-Schaeffer for rationals}\label{sec5.2.2}

In this section we specify the sets $\{\mathcal{I}_q\}_{q \in \mathbb{N}}$ appearing in Proposition~\ref{equidist_prop} further with the view to obtaining a version of Theorem~\ref{new_thm10} with \eqref{enough_mass} and \eqref{quasi_indep} checked. In particular, this will allow us to provide a solution to the Main Problem and indeed a full proof of the weak inhomogeneous Duffin-Schaeffer conjecture for \emph{all rational} $\gamma$, see Theorem~\ref{main_thm_rational} below.

\medskip


\begin{theorem}\label{general_DS}
Let $\gamma \in \R$ and $\psi: \N \to \mathbb{R}_{\ge0}$.
Let $\frac{p_k}{q_k}$ denote the $k$-th convergent of the continued fraction expansion of $\gamma$ where in the case of a rational $\gamma = [a_0;a_1,\ldots,a_i]$ we define $\frac{p_k}{q_k} = \frac{p_i}{q_i}=\gamma$ for $k \geq i$.
Suppose that there exists a sequence of finite subsets $\cS_k\subset\Z$ such that $\min\cS_k\to+\infty$ as $k\to\infty$, satisfying
\begin{equation}\label{Sk_assumption1}
 \left\lvert  \gamma - \frac{p_k}{q_k}\right\rvert\le \psi(q)\qquad\text{for $q\in \cS_k$}
\end{equation}
and
    \begin{equation}\label{big_Mass}
\sum_{\substack{q \in \mathcal{S}_k} \\ }\frac{\psi(q)\varphi(q)}{q}
    \geq q_k^8
    \end{equation}
for infinitely many $k$. 
Furthermore, assume that
\begin{equation}\label{large_psi}
\lim_{k \to \infty}\max_{\substack{q \in \mathcal{S}_k\\ {\psi(q)\ge1/2}}} \frac{\psi(q)\varphi(q,q_k)}{q} = 0
\end{equation}
where $\varphi(q,q_k)$ is given by \eqref{eq103} and  {the maximum is set to be zero if $\psi(q)<1/2$ for all $q\in\cS_k$.}
 Then for almost every $x\in[0,1]$, there exist infinitely many $k\in\N$ such that for some $q\in\cS_k$ and $a\in\Z$
\begin{equation}\label{shifted_DS-type}
\lvert qx- a - \gamma \rvert \leq \psi(q)\qquad\text{and}\qquad
(p_k + aq_k,q) = 1\,.
\end{equation}
In particular, there exist infinitely many $q \in \N$ such that
\begin{equation}\label{eq66}
\lVert qx - \gamma \rVert \leq \psi(q)\,.
\end{equation}
\end{theorem}

\begin{rem}
    We note that by \eqref{eq103}
    \begin{equation}
        \varphi(q) \leq \varphi(q,q_k) \ll \varphi(q)\Log \Log (q,q_k)\le \varphi(q)\Log \Log q_k,
    \end{equation}
    and therefore condition \eqref{big_Mass} in Theorem \ref{general_DS} immediately implies that both conditions \eqref{Div-Sum} and \eqref{DivSumI_n} hold, in which $\cI_q$ is taken to be as in \eqref{def_Iq} with $A_q=p_k$ and $B_q=q_k$ for the smallest $k$ satisfying $q\in\cS_k$, whenever such $k$ exists. We emphasize that, with reference to \eqref{def_Iq}, in general  \eqref{DivSumI_n} does not imply \eqref{Div-Sum} since it might be that $B_q \to \infty$ and therefore it is possible that $\limsup_{q \to \infty} \frac{\lvert \mathcal{I}_q\rvert}{\varphi(q)} = \limsup_{q \to \infty} \frac{\varphi(q,B_q)}{\varphi(q)} = \infty$. 
\end{rem}

\medskip

\begin{rem} We note that the choice of $q_k^8$ on the right-hand side of \eqref{big_Mass} can be seen as some form of extra-divergence condition. The value of $q_k^8$ can most probably be decreased replaced with a significantly smaller quantity. However, with the methods used in this article, it seems impossible to replace $q_k^8$ with a constant or, even stronger, with a term tending to $0$.
\end{rem}
\medskip

In the special case of $\gamma \in \mathbb{Q}$, the left hand side of \eqref{Sk_assumption1} vanishes for all sufficiently large $k$ and thus \eqref{Sk_assumption1} holds for all sufficiently large $k$. 
Furthermore, for rational $\gamma$, $\varphi(q,q_k) \asymp_{\gamma} \varphi(q)$.
This enables us to deduce the following consequence of Theorem \ref{general_DS}, which confirms Conjecture~\ref{ICC} in the rational case, and indeed answers parts (i) and (ii) of the Main Problem in the rational case.

\medskip

\begin{theorem}[The case of rational $\gamma$]\label{main_thm_rational}
Let $\psi: \mathbb{N} \to \R_{\ge 0}$ be an arbitrary function and $\gamma = A/B \in \mathbb{Q}$, where $A,B \in \mathbb{Z}$, $B>0$ and $(A,B)=1$.  Then for almost every $x\in[0,1]$ the number of integer pairs $(a,q)$ satisfying
\begin{equation}\label{eq75}
\lvert qx - \gamma - a\rvert \leq \psi(q) \qquad\text{and}\qquad
    (A + aB,q) = 1
\end{equation}
is infinite if \eqref{Div-Sum} holds, and finite otherwise.
Consequently, the weak inhomogeneous Duffin-Schaeffer conjecture holds for all $\gamma\in\Q$. That is, for all $\gamma\in\Q$, and for all $\psi: \N \to \R_{\ge0}$ satisfying \eqref{Div-Sum}, for almost every $x\in[0,1]$ inequality \eqref{eq66}
holds for infinitely many $q\in\N$.
\end{theorem}

\medskip

\begin{rem}
We note that unlike Theorems~\ref{new_thm10}\;--\;\ref{general_DS}, 
Theorem~\ref{main_thm_rational} is obtained without any constrains on the size of $\psi$. Indeed, Theorems \ref{new_thm10} and \ref{new_thm10v2} are only obtained for $\psi: \N \to [0,\tfrac{1}{2}]$. The reason is that the assumption of $\tfrac{1}{q}\mathcal{I}_q$  being uniformly distributed does not in general imply that the elements of $\frac{1}{q}\mathcal{I}_q$ are approximately equally spaced. Consequently, the intervals generated by
$$
\left\lVert x - \frac{a+\gamma}{q} \right\rVert \leq \frac{\psi(q)}{q}
$$
that contribute to $E_q^{\cI}$ may overlap for larger values of $\psi$ for different $a$. This overlapping may prevent property {\bf(B1)} from being true. However, when $\mathcal{I}_q$ is given by \eqref{def_Iq} and $\psi$ does not grow too quickly, we have enough additional information about the spacing between the elements of $\tfrac1q\cI_q$ to deduce {\bf(B1)}. This enables us to prove Theorem \ref{general_DS} for any $\psi$ that satisfies assumption \eqref{large_psi} which is less restrictive than $\psi\le\tfrac12$. In Theorem~\ref{main_thm_rational} all restrictions on $\psi$ are finally removed.
\end{rem}

\subsection{Duffin--Schaeffer with congruence relations}\label{DSCongruence}

Here we derive new results regarding Diophantine approximation with congruence relations. In this setting one is interested in approximations by integers lying in given arithmetic progressions, or in other words satisfying given  congruence relations. In dimension $1$ we deal with  approximations by rational fractions $a/q$  satisfying     \begin{equation}\label{cong_relations}
a \equiv r \!\!  \pmod t,\qquad q \equiv s \!\!  \pmod u\,,
    \end{equation}
where $r,s\in\Z$ and $t,u\in\N$ are fixed. This setting of Diophantine approximation has a surprising rich body of research with results obtained both in dimension $1$, e.g. \cite{Adiceam_2015, Harman_1988} and higher dimensions \cite{MR4312709, MR4153361}. 
In particular, Harman \cite[Corollary of Theorem~2]{Harman_1988} obtained the following Khintchine type result: {\em for any monotonic function $\psi:\N\to[0,\tfrac12]$ such that 
$$
\sum_{q=1}^\infty\psi(q)=\infty
$$ 
for almost every $x\in[0,1]$ there are infinitely many $a,q\in\N$ subject to \eqref{cong_relations} that satisfy 
\begin{equation}\label{cong_relations_inequality}
|qx-a|<\psi(q)\,.
\end{equation}}%
In fact, Harman \cite{Harman_1988} proved a quantitative version of this result which also included the inhomogeneous case. In this paper we remove the monotonicity assumption from the above statement, thus establishing the Duffin--Schaeffer conjecture with congruence relations. The first result (Theorem~\ref{cong_thm1}) deals with the case  a congruence relation is imposed on the numerators only. 
The second result (Theorem~\ref{cong_thm2}) deals with the more general case when congruence relations are imposed on both numerators and denominators.




\begin{theorem}[Duffin--Schaeffer with one congruence relation]\label{cong_thm1}
Let $r,t \in \Z$, $t>0$, $\psi: \N \to \R_{\ge0}$ be arbitrary, and suppose that
$$
    \label{div_assum_new1}
    \sum_{q=1}^\infty \frac{\psi(q)\varphi(q)}{q} = \infty.
$$
Then for almost every $x\in[0,1]$ inequality \eqref{cong_relations_inequality} holds for infinitely many integer pairs $(a,q)$ satisfying $a \equiv r \pmod t$ and $(a,q)\mid(q,r,t)$. In particular, if
\begin{equation}
    \label{div_assum_new}
    \sum_{\substack{q=1\\[0.5ex](q,r,t)=1}}^\infty \frac{\psi(q)\varphi(q)}{q} = \infty
\end{equation}
then for almost every $x\in[0,1]$ inequality \eqref{cong_relations_inequality} holds for infinitely many coprime integer pairs $(a,q)$ satisfying $a \equiv r \pmod t$.
\end{theorem}

\begin{proof}[Proof of Theorem~\ref{cong_thm1} modulo Theorem~\ref{main_thm_rational}]
To begin with, note that the `in particular' part of the statement follows from the main part of the statement applied to the following modified approximating function
$$\tilde\psi(q):=\psi(q)\cdot\mathds{1}_{(q,r,t)=1}(q)\,.
$$
Let $d=(r,t)$, $A=r/d$ and $B=t/d$. Then $B>0$ and $(A,B)=1$ and, by Theorem~\ref{main_thm_rational}, for almost every $x\in[0,1]$ there are  infinitely many integer pairs $(a,q)$ with $q>0$ satisfying \eqref{eq75}. In particular, this holds for almost every $x\in[0,\tfrac1{Bd}]$. Note that the function $x\mapsto dBx$ maps any subset of full Lebesgue measure in $[0,\tfrac1{Bd}]$ onto a subset of full Lebesgue measure in $[0,1]$. Then, on multiplying the inequality in \eqref{eq75} by $Bd$ and letting $a'=d(A+aB)$ and $x'=dBx$ we get that for almost all $x'\in[0,1]$ the inequality
\begin{equation}\label{eq81}
|qx'-a'|<Bd\psi(q)
\end{equation}
holds for infinitely many integer pairs $(a',q)$ with $q>0$ and $(a',q)\mid(q,r,t)$. On renaming, in \eqref{eq81}, $x'$ into $x$ and $a'$ into $a$, and using Cassels' Lemma (Lemma~\ref{CasselsLemma}) to remove the factor $Bd$ from \eqref{eq81}, we obtain the required result.
\end{proof}

The following theorem is an obvious consequence of Theorem~\ref{cong_thm1} obtained by redefining $\psi$ to be zero for $q\not\equiv s\pmod u$. On the other hand, Theorem~\ref{cong_thm1} can be regarded as a special case of Theorem~\ref{cong_thm2}, namely the one with $u=1$.

\begin{theorem}[Duffin--Schaeffer with two congruence relations]\label{cong_thm2}
Let $r,t,s,u \in \N$, $\psi: \N \to \R_{\ge0}$ be arbitrary, and suppose that
\begin{equation}
    \label{div_assum_new3}
    \sum_{q\equiv s \!\!\!\!    \pmod u} \frac{\psi(q)\varphi(q)}{q} = \infty.
\end{equation}
Then for almost every $x\in[0,1]$ inequality \eqref{cong_relations_inequality} holds for infinitely many integer pairs $(a,q)$ satisfying \eqref{cong_relations} and $(a,q)\mid(q,r,t)$. In particular, if
\begin{equation}
    \label{div_assum_new2}
    \sum_{\substack{q\equiv s \!\!\!\! \pmod u\\[0.5ex] (q,r,t)=1}} \frac{\psi(q)\varphi(q)}{q} = \infty
\end{equation}
then for almost every $x\in[0,1]$ inequality \eqref{cong_relations_inequality} holds for infinitely many coprime integer pairs $(a,q)$ satisfying \eqref{cong_relations}.
\end{theorem}

\medskip
\begin{rem}
It would be interesting to resolve the inhomogeneous problems and conjectures stated in \S\ref{DSBack} in the more general setting of approximations with congruence relations. 
\end{rem}

\subsection{Explicit irrational constructions}\label{explicit}

In the remaining results of this section we will make the first step towards understanding the Main Problem and conjectures in the case of irrational $\gamma$. Our immediate goal will be to construct examples of irrational $\gamma$ such that for almost all $x$ \eqref{eq66} holds for infinitely many $q\in\N$ provided \eqref{Div-Sum} or a similar assumption with extra divergence holds. We note that for any fixed given $\psi$, simply the existence of such a $\gamma$ is not an issue. This can be deduced from Cassels' doubly metric results \cite{MR0349591}, see also \cite[Theorem~8]{MR2508636}. Indeed, by Cassels' result we know that for almost all $(x,\gamma)\in\R^2$,  \eqref{eq66} holds for infinitely many $q\in\N$ as long as $\sum_{q=1}^\infty\psi(q)=\infty$. Thus, applying Fubini's theorem to Cassels' result gives the desired existence statement, albeit without the coprimality or any other potentially necessary restrictions. 
Two more downsides of using a doubly metric statement are that $\gamma$ is dependent on $\psi$ and that one cannot prove anything for $\gamma$ lying in a `thin' set such as the sets of badly or very-well approximable numbers, or Liouville numbers. The latter is `particularly thin' since it has zero Hausdorff dimension. 
We will be interested in obtaining results for $\gamma$ representing `thin' subsets of $\R$, in particular the set of  Liouville numbers. 
Recall that an irrational number $\gamma$ is Liouville if for any $w>0$ the inequality
$$
\left\lvert  \gamma - \frac{p}{q}\right\rvert<\frac{1}{q^w}
$$
holds for infinitely many $p,q\in\Z$, $q>0$. Equivalently, $\gamma$ is Liouville if and only if the convergents $p_k/q_k$ of the continued fraction expansion of $\gamma$ satisfy
\begin{equation}\label{Liouvdef}
\limsup_{k\to\infty}\frac{\log a_{k+1}}{\log q_k}=\infty\,.
\end{equation}

Our first result for irrational $\gamma$ provides a Liouville example $\gamma$ for each $\psi$.

\begin{thm}[Explicit construction for a given $\psi$]\label{thm_irrational}
Let $\{k_i\}_{i\in\N}$ be any increasing sequence of natural numbers. Suppose that for every $k\in\cK_0:=\N\setminus\{k_i:i\in\N\}$, a natural number $a_k$ is given. Then, for any function $\psi$ satisfying \eqref{Div-Sum} there are $($explicitly constructible$)$ positive integers $a_{k_i}$, $i=1,2,\dots$ such that  $\gamma=[0;a_1,a_2,\dots]$ is a Liouville number 
and 
for almost every $x\in\R$ there are infinitely many $(q,a)\in\N\times\Z$ such that  
\begin{equation}\label{shifted_DS-type+}
\lvert qx- a - \gamma \rvert \leq \psi(q)\qquad\text{and}\qquad
(p_{k} + aq_{k},q) = 1\quad\text{for some $k\in\N$.}
\end{equation}
In particular, 
for almost every $x\in\R$ there are infinitely many $q \in \N$ satisfying \eqref{eq66}. Furthermore, when $\cK_0$ is finite, $\gamma$ can be chosen optimally in the sense that for any $\tilde\psi:\N\to[0,+\infty)$ that does not satisfy  \eqref{Div-Sum},
for almost every $x\in\R$ there are only finitely many $(q,a)\in\N\times\Z$ such that
\begin{equation}\label{eq87}
\lvert qx- a - \gamma \rvert \leq \tilde\psi(q)\qquad\text{and}\qquad
(p_{k} + aq_{k},q) = 1\quad\text{for some $k\in\N$.}
\end{equation}
\end{thm}

The next result provides a version of the above in which $\gamma$ is independent of $\psi$. The cost is an `extra' divergence assumption that need to be imposed on $\psi$. 

\medskip

\begin{thm}[Uniform $\gamma$ assuming extra divergence]\label{extra_div_cor}
Fix any monotonic function $f$ such that $\lim_{Q\to\infty}f(Q)=\infty$. Then there exists an $($explicitly constructible$)$ Liouville number $\gamma$ depending only on $f$ such that for all $\psi$ satisfying any of the conditions
   \begin{equation}\label{extra_div_cond1}
   \exists \; \ve>0\qquad\liminf_{Q \to \infty}\frac{1}{f(Q)}\sum_{q \leq Q}\frac{\psi(q)\varphi(q)}{q\Log q (\Log \Log q)^{1 + \varepsilon}} \;>0\end{equation} or 
    \[
    \exists \; \ve>0\qquad\liminf_{Q \to \infty}\frac{1}{f(Q)}\sum_{q \leq Q}\frac{\psi(q)}{\Log q (\Log \Log q)^{2 + \varepsilon}} \;>0\,,\]
    we have that for almost every $x\in\R$ there are infinitely many $q \in \N$ satisfying \eqref{eq66}.
\end{thm}
\medskip

Theorem~\ref{extra_div_cor} is if fact a simple corollary of the following more technical and seemingly more restrictive statement.

\begin{thm}[Uniform $\gamma$ assuming extra divergence for special $\psi$]\label{thm_extra_div}
Fix any monotonic function $f$ such that $\lim_{Q\to\infty}f(Q)=\infty$. Then there exists an $($explicitly constructible$)$ Liouville $\gamma$ depending only on $f$ such that for all $\psi$ satisfying the following two conditions
    \[\begin{split}
    &\forall\; q \in \N\;\quad \psi(q) \neq 0\; \Rightarrow\;\psi(q) \geq \frac{1}{q}\,, \\[2ex]
    &\liminf_{Q \to \infty} \frac{1}{f(Q)}\sum_{q \leq Q}\frac{\psi(q)\varphi(q)}{q} > 0,
    \end{split}\]
we have that for almost every $x\in\R$, there are infinitely many $q \in \N$ satisfying \eqref{eq66}.
\end{thm}

\medskip

\begin{proof}[Proof of Theorem~\ref{extra_div_cor} modulo Theorem~\ref{thm_extra_div}]
Since $$\frac{q}{\varphi(q)} \ll \Log \Log q,$$ it suffices to show Theorem~\ref{extra_div_cor} assuming condition
    \eqref{extra_div_cond1}.
But then, since $$\sum_{q \in \mathbb{N}}\frac{\varphi(q)/q }{q\Log q (\Log \Log q)^{1 + \varepsilon}} < \infty,$$ we must have that
\[
\liminf_{Q \to \infty} \frac{1}{f(Q)}\sum\limits_{\substack{q \leq Q\\ \psi(q) \geq 1/q}}\frac{\psi(q)\varphi(q)}{q\Log q (\Log \Log q)^{1 + \varepsilon}} = \infty\,.\]
This means that we can assume without loss of generality that $\psi(q) \geq \frac{1}{q}$ for all $q \in \supp(\psi)$ and applying Theorem~\ref{thm_extra_div} gives the required result.
\end{proof}

\section{Proofs of new results on Duffin--Schaeffer   \label{PDSlove}}

\subsection{Proof of Theorems~\ref{new_thm10} and \ref{new_thm10v2}}\label{sec6.1}

To begin with, note that \eqref{ud} implies that $\lvert \mathcal{I}_q \rvert \to \infty$ as $q \to \infty$. Next, observe that, by \eqref{EE} and the fact that $\psi(q)\le\tfrac12$ and $\gamma\in\R$, we have that
\begin{equation}\label{eqLB}
  \mu(E^\cI_q) = \frac{2\psi(q)|\cI_q|}{q}\,.
\end{equation}
Next, using $\gamma\in\R$ and $\psi(q)\le\tfrac12$, we obtain that
\begin{align}
\nonumber\mu(E_q^{\mathcal{I}} \cap [x,y]) &\leq \frac{2\psi(q)}{q} \;\# \left\{
a \in \mathcal{I}_q: \frac{a}{q} \in \left[x - \frac{\psi(q) +\gamma}{q},y + \frac{\psi(q)-\gamma}{q}\right] \pmod 1
\right\}\\[1ex]
\nonumber&\leq \frac{2\psi(q)}{q}  \left(\# \left\{
a \in \mathcal{I}_q: \frac{a}{q} \in \left[x ,y\right]
\right\} + O(1)\right)\\[2ex]
& \stackrel{\eqref{ud}}{=} \frac{2\psi(q)}{q}  \Big((y-x)|\cI_q|+o(|\cI_q|)+O(1)\Big).\label{eq74}
\end{align}
Similarly, we have that
\[
\mu(E_q^{\mathcal{I}} \cap [x,y]) \geq \frac{2\psi(q)}{q}  \Big((y-x)|\cI_q|-o(|\cI_q|)-O(1)\Big)\,.
\]
This together with \eqref{eqLB} and \eqref{eq74} imply that for any $0\le x<y\le1$
\begin{equation}\label{ud2}
\lim_{q \to \infty} \frac{\mu\left(E_q^{\mathcal{I}} \cap [x,y] \right)}{\mu(E_q^{\mathcal{I}})} =\mu([x,y])\,.
\end{equation}
Note that \eqref{ud2} implies {\bf{(B1)}}. Thus applying Theorem \ref{T2MS} gives Theorems~\ref{new_thm10}, while applying Theorem~\ref{T1MS} gives Theorem~\ref{new_thm10v2}.

\subsection{Further auxiliary statements}

\begin{lemma}\label{key_lemma}

Let $\psi: \N \to [0,\tfrac{1}{2}]$ be an arbitrary function. Let $A,B$ be coprime integers, $B>0$ and for $q\in\N$ let
    \[E_q^* = E_{q,A,B,\psi}^* := \bigcup_{\substack{a \in \mathbb{Z}_{Bq}^{*}\\ a \equiv \,A\!\!\!\! \pmod B}} \left[\frac{a}{Bq} - \frac{\psi(q)}{q}, \frac{a}{Bq} + \frac{\psi(q)}{q}\right]\pmod1.\]
Then for any finite set $\mathcal{S} \subseteq \mathbb{N}$ such that 
\begin{equation}\label{cond001}
\sum_{q \in \mathcal{S}} \frac{\psi(q)\varphi(q)}{q} \geq B^8\,,
\end{equation}
we have that
\begin{equation}\label{concl001}
\sum_{q,r \in \mathcal{S}}
\mu(E_q^* \cap E_r^*) \ll \left(\sum_{q \in \mathcal{S}} \mu(E_q^*)\right)^2\,,
\end{equation}
where the implied constant is absolute and in particular independent from $B$.
\end{lemma}

Note that establishing Theorems~\ref{general_DS} and \ref{main_thm_rational} in full generality will require a treatment of the case $\psi(q) > \tfrac12$, which is not covered by Lemma~\ref{key_lemma}, neither it is covered by Proposition~\ref{equidist_prop} or Theorems~\ref{new_thm10} and \ref{new_thm10v2}. We will deal with this issue by establishing an appropriate  generalisation of \eqref{PV1}. Note that for $\psi(q) > \tfrac12$, formula \eqref{sum_of_measures} is no longer true since the balls (intervals) centred at $(a+\gamma)/q$, that $E_q^{\cI}$ is made of,  need not be disjoint for different $a$ -- see \eqref{EE}. 
This partially overlapping structure also makes it hard to prove property {\bf{(B1)}}. However, assuming the growth of $\psi(q)$ is restricted, we still recover the  asymptotic behaviour seen in the case $\psi\le\tfrac12$.

\begin{proposition}\label{large_psi_measure}
     Let $A_q,B_q \in \mathbb{Z}$, $B_q>0$,
     $(A_q,B_q) =1$,  let $E_q^{\mathcal{I}}$ be as in Proposition \ref{equidist_prop} and $\psi(q)\ge1/2$ for infinitely many $q$.
Furthermore, suppose that
\begin{equation}\label{decay_large_psi}\lim_{\substack{q \to \infty\\  {\psi(q)\ge1/2}}} \frac{\varphi(q,B_q)\psi(q)}{q} = 0.\end{equation}
     Then
\begin{equation}\label{eq_measure_large_psi_}
     \lim_{\substack{q \to \infty\\ {\psi(q)\ge1/2}}}
     \mu(E^\cI_q)\left(\frac{2\psi(q)\varphi(q,B_q)}{q}\right)^{-1} = 1.
     \end{equation}
     with the convergence being uniform in $B_q$.
     Furthermore, for any $0 \leq x < y \leq 1$,
\begin{equation}\label{ud_large_psi}
 \lim_{\substack{q \to \infty\\ {\psi(q)\ge1/2}}} \frac{\mu\left(E_q^{\mathcal{I}} \cap [x,y] \right)}{\mu(E_q^{\mathcal{I}})} =\mu([x,y]).
\end{equation}
\end{proposition}

In view of~\eqref{eq_measure_large_psi_} the following second moment bound can be seen as an analogue of Lemma~\ref{key_lemma} in the case of large $\psi$.

\begin{lemma}\label{key_lemma_large_psi}
For any $\psi: \N \to \R_{\ge 0}$ satisfying $\psi(q) \geq \tfrac{1}{2}\mathds{1}_{[q \in \supp \psi]}$,
    we have the following. Let $A,B$ be coprime integers, $B>0$ and let
$E^*_q := E^*_{q,A,B,\psi}$ be the same as in Lemma~\ref{key_lemma}.
Then for any finite set $\mathcal{S} \subseteq \mathbb{N}$ satisfying \eqref{cond001}, we have that 
\begin{equation}\label{conclem5}
\sum_{q,r \in \mathcal{S}}
\mu(E_q^* \cap E_r^*) \ll \left(\sum_{q \in \mathcal{S}} \frac{\varphi(q,B)\psi(q)}{q}\right)^2\,,
\end{equation}
with the implied constant being absolute and in particular, independent of $A$ and $B$.
\end{lemma}

\medskip

\begin{rem}
We note that, by Proposition~\ref{large_psi_measure}, for relatively `small' $\psi$, namely $\psi$ satisfying \eqref{decay_large_psi},  inequality \eqref{conclem5} becomes the same as \eqref{concl001}.
\end{rem}
While postponing the proof of Propositions~\ref{equidist_prop} and \ref{large_psi_measure} and Lemmas~\ref{key_lemma} and \ref{key_lemma_large_psi}  until \S\ref{proofsprop1and2},\ref{sec5.4}, we now focus on the applications of these statements.

\subsection{Proof of Theorem~\ref{general_DS}}

Since the sets $\cS_k$ are finite and $\min\cS_k\to\infty$ as $k\to\infty$, we can assume without loss of generality that these sets are disjoint. Then for every $q\in\cS:=\bigcup_k\cS_k$ there is a unique $k$ such that $q\in\cS_k$. Then, with the view to using Proposition~\ref{equidist_prop} or Proposition~\ref{large_psi_measure}, we define $A_q=p_{k}$ and $B_q=q_{k}$ if $q\in\cS_k$ for some $k$ and $A_q=0$ and $B_q=1$ otherwise. Note that  for $q\in\cS_k$ the set $\cI_q$, defined within Proposition~\ref{equidist_prop}, is explicitly given by
$$
\mathcal{I}_q = \{a \in \mathbb{Z}_q: (p_k + aq_k, q) = 1\}\,.
$$
Also without loss of generality we will assume that $\psi(q)=0$ if $q\not\in\cS$, otherwise we modify $\psi$ by multiplying it by $\mathds{1}_{\cS}$. Obviously, the modified $\psi$ also satisfies the conditions of Theorem~\ref{general_DS}, namely \eqref{Sk_assumption1}, \eqref{big_Mass} and \eqref{large_psi}. Note that to complete the proof of Theorem~\ref{general_DS} it is enough to show that  $\limsup_{q\to\infty} E^{\cI}_q$ has full measure, where $E_q^{\cI}$ is defined by \eqref{EE} and  can be written explicitly as follows:
    \[E_q^{\cI} =
    \bigcup_{\substack{a \in \mathbb{Z}_{q} \\ (p_k + aq_k,q) = 1}}
    \left[\frac{a+\gamma}{q} - \frac{\psi(q)}{q}, \frac{a+\gamma}{q} + \frac{\psi(q)}{q}\right]\pmod 1.\]

First, we will consider the case when $\psi: \mathbb{N} \to [0,\tfrac{1}{4}]$ for all $q\in\N$. Then, by Theorem~\ref{new_thm10} together with Proposition~\ref{equidist_prop}, it suffices to show that for infinitely many $k$, we have that
\begin{equation}\label{eq93}
    \sum_{q,r \in \mathcal{S}_k} \mu(E_q^{\cI} \cap E_r^{\cI}) \ll \left(\sum_{q \in \mathcal{S}_k} \mu(E_q^{\cI})\right)^2\,.
\end{equation}
Indeed, \eqref{eq93} verifies \eqref{quasi_indep}, \eqref{ud} is verified by Proposition~\ref{equidist_prop}, and \eqref{enough_mass} is a consequence of \eqref{big_Mass} and \eqref{eq103}. Consequently, Theorem~\ref{new_thm10} would imply Theorem~\ref{general_DS} in this case.

Regarding \eqref{eq93}, by \eqref{Sk_assumption1}, we have that for any $q\in\cS_k$ and any $a\in\Z$
$$
\left[\frac{a+\gamma}{q} - \frac{\psi(q)}{q}, \frac{a+\gamma}{q} + \frac{\psi(q)}{q}\right]\subset
\left[\frac{a+\frac{p_k}{q_k}}{q} - \frac{2\psi(q)}{q}, \frac{a+\frac{p_k}{q_k}}{q} + \frac{2\psi(q)}{q}\right]\,.
$$
Therefore, using $(p_k,q_k) = 1$, 
\[\begin{split}E_q^{\cI}
&\subseteq
 \bigcup_{\substack{a \in \mathbb{Z}_{q} \\ (p_k + aq_k,q) = 1}}
    \left[\frac{a+\frac{p_k}{q_k}}{q} - \frac{2\psi(q)}{q}, \frac{a+\frac{p_k}{q_k}}{q} + \frac{2\psi(q)}{q}\right]\pmod 1
    \\&= \bigcup_{\substack{a \in \mathbb{Z}_{q} \\ (p_k + aq_k,qq_k) = 1}}
    \left[\frac{aq_k+p_k}{qq_k} - \frac{2\psi(q)}{q},\frac{aq_k+p_k}{qq_k} + \frac{2\psi(q)}{q}\right]\pmod 1
    \\&= \bigcup_{\substack{a \in \mathbb{Z}_{qq_k}^{*}\\ a \equiv \, p_k\!\!\!\! \pmod {q_k}}} \left[\frac{a}{qq_k} - \frac{2\psi(q)}{q}, \frac{a}{qq_k} + \frac{2\psi(q)}{q}\right]\pmod 1
    \\[2ex]&= E^*_{q,p_k,q_k,2\psi},
    \end{split}\]
    where $E^*_{q,p_k,q_k,2\psi}$ is defined in Lemma \ref{key_lemma}. By the assumption $\psi\le\tfrac{1}{4}$, we have that $2\psi\le\tfrac12$ and therefore $\mu(E^*_{q,p_k,q_k,2\psi}) = 2\mu(E_q^{\cI})$. In fact, the measure of either of this sets is given by \eqref{sum_of_measures}, and not dependent on $\gamma$. Then, Lemma \ref{key_lemma} together with \eqref{big_Mass} imply that
    \begin{equation*}
    \begin{split}
\sum_{q,r \in \mathcal{S}_k} \mu(E_q^{\cI} \cap E_r^{\cI})
&\leq \sum_{q,r \in \mathcal{S}_k} \mu(E^*_{q,p_k,q_k,2\psi} \cap E^*_{r,p_k,q_k,2\psi})
\\&\ll \left(\sum_{q \in \mathcal{S}_k}\mu(E^*_{q,p_k,q_k,2\psi})\right)^2
\\&\ll \left(\sum_{q\in \mathcal{S}_k} \mu(E_q^{\cI})\right)^2,
\end{split}
     \end{equation*}
which concludes the proof in the case $\psi\le\tfrac{1}{4}$. In the remaining case we can assume without loss of generality that
\begin{equation}\label{div_large_psi}\sum_{\substack{q \in \N\\ \psi(q) > 1/4}} \frac{\psi(q)\varphi(q)}{q} = \infty.\end{equation}
By \eqref{large_psi}, we can assume that $\lim_{q \to \infty}\frac{\psi(q)}{q} = 0$ and thus,
 Cassels' Lemma (Lemma \ref{CasselsLemma}) enables us to replace $\psi$ by $2\psi$. Note that \eqref{Sk_assumption1},\eqref{big_Mass} and \eqref{large_psi} still hold for $2\psi$. Thus, to complete the proof, we can simply assume without loss of generality that
$\psi(q) \geq \tfrac{1}{2}\mathds{1}_{[q \in \supp \psi]}$ and \eqref{large_psi}. Then, by using Lemma \ref{key_lemma_large_psi} instead of Lemma \ref{key_lemma}, arguing as above, we get that
\[\sum_{q,r \in \mathcal{S}_k} \mu(E_q^{\cI} \cap E_r^{\cI})
\ll \left(\sum_{q \in \mathcal{S}_k} \frac{\varphi(q,q_k)\psi(q)}{q}\right)^2.
\]
By \eqref{eq_measure_large_psi_} from Proposition \eqref{large_psi_measure} and the assumption of \eqref{large_psi}, we get $\frac{\varphi(q,q_k)\psi(q)}{q} \ll \mu(E_q^{\cI})$ and therefore, \eqref{eq93} follows. Further, \eqref{ud_large_psi} from Proposition~\ref{large_psi_measure} verifies {\bf(B1)}.
 As a result, Theorem~\ref{T2MS} is applicable to the sequence $(E_q^{\cI})_{q \in \N}$ and the result follows.

\subsection{Proof of Theorem~\ref{main_thm_rational}}

We will write $E^{\psi}_q$ for the set $E_q^{\mathcal{I}}$ from \eqref{EE} induced by the condition $(A + aB,qB) = 1$ with $\gamma=A/B$.
If 
$$
\sum_{q \in \N} \frac{\psi(q)\varphi(q)}{q} < \infty
$$ 
then \eqref{large_psi} holds, and we can use \eqref{eq_measure_large_psi_} to conclude that $\sum_{q \in \N} \mu(E_q) < \infty$, so the result in this case follows on applying  the first Borel--Cantelli Lemma.

For the rest of the proof we will assume  \eqref{Div-Sum}.
We claim that it is enough to consider functions $\psi$ satisfying 
\begin{equation}
    \label{decay_psi_phi}
    \lim_{q \to \infty}\frac{\varphi(q)\psi(q)}{q} = 0.
\end{equation}

Indeed, assume for the moment that \eqref{decay_psi_phi} does not hold, then
 on letting 
$\psi_1(q) := \min\{\psi(q),\frac{q}{\varphi(q)}\}$ we  observe that 
\[\sum_{q \in \N} \frac{\psi_1(q)\varphi(q)}{q} = \infty\,.\]
Furthermore, since we assumed that \eqref{decay_psi_phi} does not hold we have that $$\limsup_{q \to \infty} \frac{\psi_1(q)\varphi(q)}{q} = \delta$$ for some $0 < \delta \leq 1$. Therefore, there exists an increasing sequence of integers $(q_n)_{n \in \N}$ such that for all $n \in \N$ we have that
\begin{equation}\label{almost_delta}\frac{\delta}{2} \leq \frac{\psi_1(q_n)\varphi(q_n)}{q_n} \leq 2\delta\,.\end{equation}
Next, we set
\[\psi_2(q) = \begin{cases} \frac{\psi_1(q_n)}{n} &\text{ if } q = q_n \text { for some $n \in \N$},\\
0 &\text{ otherwise. }
\end{cases}\]
By \eqref{almost_delta}, we get that
\[\sum_{q \in \N} \frac{\psi_2(q)\varphi(q)}{q} = \sum_{n \in \N} \frac{\psi_1(q_n)\varphi(q_n)}{q_n} \geq \frac{\delta}{2}\sum_{n \in \N}\frac{1}{n} = \infty\,.
\]
Since $\psi\ge\psi_1\ge\psi_2$ we also have that $\mu(\limsup_{q \to \infty} E_q^{\psi}) \geq \mu(\limsup_{q \to \infty} E_q^{\psi_2})$, and therefore in the proof we could replace $\psi$ with $\psi_2$ which clearly satisfies \eqref{decay_psi_phi}. This justifies the claim about $\psi$ we made above, that we can assume that $\psi$ satisfies \eqref{decay_psi_phi}.

Since $\gamma = \tfrac{A}{B}\in \mathbb{Q}$ and \eqref{Div-Sum} holds, it is a trivial task to find a sequence of subsets $\cS_k \subset \mathbb{Z}$ satisfying  \eqref{Sk_assumption1}, \eqref{big_Mass}  and $\min\cS_k\to\infty$ as $k\to\infty$, with
$\frac{p_k}{q_k} = \frac{A}{B}$. 
Note that since $(q_k)_{k \in \mathbb{N}}$ is absolutely bounded, we see that
$\varphi(q,q_k) \asymp_B \varphi(q)$ and therefore, \eqref{decay_psi_phi} implies \eqref{large_psi}. Hence,  Theorem \ref{general_DS} is applicable which completes the proof.

\subsection{Proof of Theorem~\ref{thm_irrational}}

We can assume without  loss of generality that \eqref{decay_psi_phi} holds, as otherwise we can replace $\psi$ with the functions $\psi_2$ defined in the previous section.  Further, we can assume without loss of generality that 
\begin{equation}\label{LB95}
\psi(q) \geq \frac{1}{q^{2}}\quad \text{if}\quad \psi(q) \neq 0\,.
\end{equation}
Otherwise we can replace $\psi$ with the smaller function $\tilde{\psi}(q) = \psi(q)\mathds{1}_{[\psi(q) \geq q^{-2}]},$
which also satisfies the divergent sum condition \eqref{Div-Sum} for $\sum_{q \in \N} {q^{-2}} < \infty$.

We will define $a_{k_1},a_{k_2},\dots$ inductively, simultaneously 
choosing sets $\mathcal{S}_{k_i-1}$ satisfying the conditions of  Theorem \ref{general_DS}. 
    Recall the following basic inequality from the theory of continued fractions
    \begin{equation}\label{DA_ineq}\left\lvert \gamma - \frac{p_k}{q_k}\right\rvert \leq \frac{1}{q_kq_{k+1}} \leq \frac{1}{a_{k+1}q_k^2}.\end{equation}

For $i = 1$, note that  $a_{k}$ is given for any $k\le k_1-1$ and thus  $q_{k_1-1}$ is well defined. Then, by \eqref{Div-Sum}, we are able to choose $Q_1 = Q_1(\psi)$ such that
\begin{equation}\label{w1}
\sum_{\substack{q \leq Q_1} \\ }\frac{\psi(q)\varphi(q)}{q}
    \geq q^9_{{k_1-1}}\,.
\end{equation}
Let $a_{k_1} \geq Q_1^2$. Then, using \eqref{LB95} and \eqref{DA_ineq}, for any $1\leq q \leq Q_1$ such that 
    $\psi(q) \neq 0$, we get that
\begin{equation}\label{w2}
\left\lvert \gamma - \frac{p_{k_1-1}}{q_{k_1-1}}\right\rvert \leq
    \frac{1}{a_{k_1}q_{k_1-1}^2}\leq \frac{1}{q_{k_1-1}Q_1^2}\leq \frac{1}{q_{k_1-1}q^2} \leq \frac{\psi(q)}{q_{k_1-1}}\,.
\end{equation}
Define $\mathcal{S}_{k_1-1} := [1,Q_1] \cap \supp(\psi)$.

Now suppose $\mathcal{S}_{k_{i'}-1}$, $Q_{i'}$ and $a_{k_{i'}}$ have been chosen for $i'\le i-1$. By the conditions of the theorem, we are also given $a_k$ for $k\le k_{i}-1$, which makes $q_{k_i-1}$ well defined.  
Choose $Q_{i}$ such that
\begin{equation}\label{w3}
\sum_{\substack{Q_{i-1} < q \leq Q_{i}} \\ }\frac{\psi(q)\varphi(q)}{q}
    \geq {q_{k_i-1}^9}.
\end{equation}
Again this is possible thanks to \eqref{Div-Sum}.
    Let $a_{k_i} \ge Q_{i}^2$ and $\mathcal{S}_{k_i-1} :=
    [Q_{i-1},Q_{i}] \cap \supp(\psi)$. Then, using \eqref{LB95} and \eqref{DA_ineq} for any $q \in \mathcal{S}_{k_i-1}$, we have that
\begin{equation}\label{w4}
\left\lvert \gamma - \frac{p_{k_i-1}}{q_{k_i-1}}\right\rvert \leq
    \frac{1}{a_{k_i}q^2_{k_i-1}}\leq
      \frac{1}{Q_{i}^2q_{k_i-1}}\leq\frac{1}{q^2q_{k_i-1}} \leq \frac{\psi(q)}{q_{k_i-1}}\,.
\end{equation}
Observe that \eqref{Liouvdef} can be met by additionally requiring that $\log a_{k_i}\ge i\log q_{k_i-1}$.

Now let $$\tilde{\psi}(q) := \frac{\psi(q)}{q_{k_i-1}}\qquad\text{for }q\in\cS_{k_i-1}\,.$$ 
We claim that 
\eqref{Sk_assumption1},\eqref{big_Mass}, and \eqref{large_psi} are satisfied for $k=k_i-1$, $i \in \mathbb{N}$ when $\psi$ is replaced by $\tilde{\psi}$. Indeed, \eqref{Sk_assumption1} follows from \eqref{w2} and \eqref{w4}, \eqref{big_Mass} follows from \eqref{w1} and \eqref{w3}, whereas \eqref{large_psi} follows from \eqref{decay_psi_phi} and the fact that 
\[\varphi(q,q_{k_i-1}) = \varphi(q)\prod_{p \mid (q,q_{k_i-1})} \left(1 + \frac{1}{p-1}\right) \ll \varphi(q)\Log \Log q_{k_i-1}\,.
\]
Thus, the assumptions of Theorem \ref{general_DS} are satisfied for $\tilde{\psi}$ and 
since $\tilde{\psi}(q) \leq \psi(q)$, the required result regarding \eqref{shifted_DS-type+} follows on applying Theorem \ref{general_DS}.

\medskip

Now we turn to the  furthermore part.
In the above proof, $a_{k_i}$ is chosen to be any integer satisfying 
$$
a_{k_i}\ge\max\{Q_i,e^iq_{k_i-1}\}\,.
$$ 
By Dirichlet's theorem on primes in arithmetic progressions, we can choose $a_{k_i}$ so that $q_{k_i}=a_{k_i}q_{k_i-1}+q_{k_i-2}$ is a prime number. In particular, if $\cK_0$ is finite, we can ensure that $q_k$ is a prime number for all but finitely many $k$. In this case the number of $a\in\Z_q$ satisfying the right hand side of \eqref{eq87} is  $\asymp \varphi(q)$. Then, the furthermore part of Theorem~\ref{thm_irrational} follows as a simple application of the first Borel-Cantelli lemma.

\subsection{Proof of Theorem~\ref{thm_extra_div}}

Since $\|qx-\gamma\|\le\tfrac12$ for any $x\in[0,1]$, any $\gamma\in\R$ and any $q\in\N$, the statement is trivial if $\psi(q)\ge\tfrac12$ for infinitely many $q$. Thus, in  what follows we can assume without loss of generality that $\psi(q)<\tfrac12$ for all $q\in\N$.
Suppose we are given any increasing sequence $\{k_i\}_{i \in \mathbb{N}}$ of natural numbers and any 
collection of natural numbers $a_k$ with $k\in\N\setminus\{k_i:i\in\N\}$.
We construct $\gamma = [0;a_1,a_2,\ldots]$ iteratively by choosing
\[a_{k_{i+1}} \ge  \min\{q \in \N: f(q) \geq q_{k_{i+1}-1}^9\}\quad\text{for}\quad i \in \N,\]
which is well-defined since $f(q) \to \infty$.
Now let $\psi$ be as in Theorem \ref{thm_extra_div}.
We claim that 
$$
\mathcal{S}_{k_{i+1}-1} := [q_{k_{i+1}-1}+1,q_{k_{i+1}}] \cap \supp(\psi)
$$ 
satisfies the assumptions of Theorem \ref{general_DS} for all $i \geq i_0(\psi)$. Indeed, since
$$
\liminf_{Q \to \infty} \frac{1}{f(Q)}\;\sum_{q \leq Q}\frac{\psi(q)\varphi(q)}{q} > 0\,,
$$
there exist $Q_0>0$ and $\delta_0 > 0$ such that for any
$Q \geq Q_0$, we have that
$$
\sum_{q \leq Q}\frac{\psi(q)\varphi(q)}{q} \geq \delta_0 f(Q)\,.
$$
Therefore, for $Y \geq Q_0$ we obtain that
\[
\sum_{X < q \leq Y}\frac{\psi(q)\varphi(q)}{q} \geq
\delta_0 f(Y) - X.
\]
Since $q_{k_{i+1}} \geq a_{k_{i+1}}$, we obtain, for $i$ large enough, that
\[
\sum_{q_{k_{i+1}-1} < q \leq q_{k_{i+1}}}
\frac{\psi(q)\varphi(q)}{q} \geq
\delta_0 f(a_{k_{i+1}}) - q_{k_{i+1}-1}
\geq \delta_0 q_{k_{i+1}-1}^9 - q_{k_{i+1}-1}
\geq  q_{k_{i+1}-1}^{8}
.\]

Since $\psi(q) \geq \frac{1}{q}$ for $q \in \supp(\psi)$, for any
$q \in \mathcal{S}_{k_{i+1}-1}$ we have that
\[
\left\lvert \gamma - \frac{p_{k_{i+1}-1}}{q_{k_{i+1}-1}}\right\rvert \leq \frac{1}{q_{k_{i+1}}} \leq \frac{1}{q} \leq \psi(q).
\]

This verifies the assumptions for Theorem \ref{general_DS} and thus the statement follows, since we can additionally demand that $\log a_{k_i}\ge i\log q_{k_i-1}$ to ensure that $\gamma$ is a Liouville number.

\subsection{Proof of Propositions~\ref{equidist_prop} and \ref{large_psi_measure}}
\label{proofsprop1and2}

\begin{proof}[Proof of Proposition~\ref{equidist_prop}]

In order to show the uniform distribution property, it suffices to verify that for any $0 \leq y \leq 1$, we have
that
\begin{equation}\label{show_ud}
\lim_{q \to \infty}\frac{\left\lvert \left\{
a \in \mathcal{I}_q: \frac{a}{q} \in \left[0 ,y\right)\right\} \right\rvert}
{\lvert \mathcal{I}_q \rvert } = y\,.
\end{equation}
Note that
\[\begin{split}
\# \Big\{a \in \mathcal{I}_q: \frac{a}{q} & < y \Big\}
\;=\; \#\big\{0 \leq a < qy: (\underbrace{A_q + aB_q}_{b},q) = 1\big\}
\\&=  \#\big\{A_q \leq b < A_q+qyB_q: b \equiv
A_q \pmod {B_q},\;\; (b,q) = 1\big\}.
\end{split}
\]

Let $P(q) := \prod\limits_{{p \mid q,\;p \nmid B_q}} p$, where the product is taken over primes $p$. Then, since $A_q$ is coprime to $B_q$, we have that
\[
\left\{\begin{array}{r}
b\equiv A_q \pmod {B_q}\,, \\[1ex]  (b,q) = 1\,,
\end{array}\right.
\;\iff\;
\left\{\begin{array}{l}
b\equiv A_q \pmod {B_q}\,,  \\[1ex]  (b,P(q)) = 1\,.
\end{array}\right.
\]
By a standard inclusion-exclusion argument involving the M\"obius function $\mu$, we can therefore write that
\begin{equation}\label{eq99}
\# \left\{a \in \mathcal{I}_q: \frac{a}{q} < y \right\}
= \sum_{d \mid P(q)}\mu(d)\sum_{\substack{A_q \leq b < qyB_q+ A_q\\ d\mid b}}  \mathds{1}_{[b \equiv A_q \pmod {B_q}]}\,.
\end{equation}
Since $(d,B_q)=1$ for arbitrary values $X < Y$ we have that
\begin{equation}\label{eq100}
\sum_{\substack{X \leq b < Y\\ d\mid b}}  \mathds{1}_{[b \equiv A_q \pmod {B_q}]} = \sum_{\substack{X \leq b < Y\\ b \equiv b_0 \pmod{B_qd}}}  1
\end{equation}
where $b_0$ is the unique solution $\pmod{B_qd}$ to the system 
$$
b \equiv 0 \pmod d,\qquad b \equiv A_q \pmod{B_q}\,,
$$
which exists by the Chinese Remainder Theorem.
Then, clearly
\begin{equation}\label{eq101}
\sum_{\substack{X \leq b < Y\\ b \equiv b_0 \pmod{B_qd}}}  1 = \frac{Y-X}{B_qd} + O(1)\,.
\end{equation}
Furthermore,  the error term $O(1)$ in \eqref{eq101} vanishes if $X,Y\in\Z$ and $(B_qd)\mid(Y-X)$.

Now, combining \eqref{eq100} and \eqref{eq101} together with \eqref{eq99} with $X=A_q$ and $Y=A_q+qyB_q+A_q$ gives
\[\begin{split}
\# \left\{a \in \mathcal{I}_q: \frac{a}{q} < y \right\}
&= \sum_{d \mid P(q)}\mu(d)\left(\sum_{\substack{0 \leq b < qyB_q\\ d\mid b}}  \mathds{1}_{[b \equiv A_q \pmod {B_q}]} + O(1)\right)\\[1ex]
&\stackrel{\text{letting }b=dm}{=} \sum_{d \mid P(q)}\mu(d)\left(\sum_{\substack{0 \leq m < \frac{q}{d}B_qy}}  \mathds{1}_{[md \equiv A_q \pmod {B_q}]} + O(1)\right)\\[1ex]
&= \sum_{d \mid P(q)}\mu(d)\left(\sum_{\substack{0 \leq m < \frac{q}{d}B_qy}}  \mathds{1}_{[m \equiv d^{-1}A_q \pmod {B_q}]} + O(1)\right),
\end{split}
\]
where we used in the last line that $(P(q),B_q) = 1$. Note that when $y=1$ we have that $Y-X=qB_q$ is divisible by $B_qd$ and therefore in this case the error term in the above calculation vanishes.
Further, since $(A_q,B_q) = 1$, the above implies that
\begin{equation*}
\sum_{\substack{0 \leq m < \frac{q}{d}B_qy}}  \mathds{1}_{[m \equiv d^{-1}A_q \pmod {B_q}]}
= \frac{q}{d}y + O(1),
\end{equation*}
which shows that
\begin{align}
\# \left\{
a \in \mathcal{I}_q: \frac{a}{q} < y \right\}
&=  qy \sum_{d \mid P(q)}\frac{\mu(d)}{d} + O(\tau(P(q))\nonumber\\[1ex]
&= y\varphi(q)\prod_{p \mid (q,B_q)}\left(1 + \frac{1}{p-1}\right) + O(\tau(q))\,,\label{number_of_points}
\end{align}
where $\tau$ stands for the number of divisors, and we used the well-known identity $$\sum\limits_{d \mid q}\frac{\mu(d)}{d} = \frac{\varphi(q)}{q} =
\prod\limits_{p \mid q} \left(1 - \frac{1}{p}\right)$$
and the obvious fact that
$$
\prod_{p|P(q)}\left(1-\frac1p\right)\prod_{p|(P(q),B_q)}\left(1-\frac1p\right)=\prod_{p|q}\left(1-\frac1p\right)
=\frac{\varphi(q)}{q}\,.
$$

Note that in the case $y = 1$ the error term in the above calculation vanishes, and we obtain that
\eqref{eq103}. 
This immediately proves \eqref{sum_of_measures}. Furthermore, since
$\lim_{q \to \infty} \frac{\tau(q)}{\varphi(q)} = 0$, \eqref{show_ud} follows
from \eqref{number_of_points}, which finishes the proof.
\end{proof}

\bigskip

\begin{proof}[Proof of Proposition \ref{large_psi_measure}]
We start by establishing \eqref{eq_measure_large_psi_}. First of all, since $E^{\cI}_q$ is the union of $\varphi(q,B_q)$ many intervals of length $2\psi(q)/q$, we immediately get that
$$
\mu(E^{\cI}_q) \leq \frac{2\psi(q)\varphi(q,B_q)}{q}\,.
$$
Thus it suffices to show that under the assumption of \eqref{decay_large_psi}, we have that
\begin{equation}\label{large_psi_lower_bound}
     \lim_{\substack{q \to \infty\\\psi(q)\ge1/2}}
     \mu(E^\cI_q)\left(\frac{2\psi(q)\varphi(q,B_q)}{q}\right)^{-1} \geq 1\,.\end{equation}
For this, we are adapting the ideas of \cite[\S 2]{pv} to various coprimality conditions. Recall that $|\cI_q|=\varphi(q,B_q)$, see \eqref{eq103}. Let $a_1 < a_2 < \ldots < a_{\varphi(q,B_q)}$ be the elements of $\mathcal{I}_q$ written in the ascending order. We also extend them both ways modulo $q$, by defining $a_0 = a_{\varphi(q,B_q)} - q$ and $a_{\varphi(q,B_q)+1} = a_1 + q$. Then we get that
\begin{equation}\label{extract_T}\begin{split}
\mu(E^{\cI}_q)&= \sum_{i = 1}^{\varphi(q,B_q)} q^{-1} \left(\min\left\{\psi(q),\frac{a_{i+1}-a_i}{2}\right\} + \min\left\{\psi(q),\frac{a_{i}-a_{i-1}}{2}\right\}\right)
\\&=  \sum_{i = 1}^{\varphi(q,B_q)} q^{-1} \min\left\{2\psi(q),\frac{a_{i+1}-a_i}{2}\right\}
\\&= 2\frac{\psi(q){\varphi(q,B_q)}}{q} - \frac{T}{q}
\end{split}
\end{equation}
where
\begin{equation}
    \label{T_integral}
T = \sum_{\substack{i = 1\\a_{i+1}-a_i \leq 2\psi(q)}}^{\varphi(q,B_q)} 2\psi(q) - (a_{i+1} - a_i)
= \int_{0}^{2\psi(q)} u(x) \mathrm{d}x
\end{equation}
and $u$ is the step function defined by
\[
u(x) = \sum_{\substack{i = 1\\a_{i+1}-a_i \leq x}}^{\varphi(q,B_q)} 1\,.
\]
Note that
\begin{equation}
\label{u_bound}
u(x) \leq \sum_{h \leq x} F(h,q)
\end{equation}
where 
\[F(h,q) := \#\left\{a \in \Z_q: (A+aB_q,q) = (A+aB_q + hB_q,q) = 1\right\}.
\]
Observe that $F(h,q)$ is a multiplicative function of $q$. Thus to compute $F$ it suffices to consider the case of $q = p^k$ for a prime $p$ and $k\in\N$. One checks straightforwardly that
\begin{equation}
F(h,p^k) = \begin{cases}
    p^{k-1}(p-2) &\text{ if } p \nmid hB_q,\\
    p^{k-1}(p-1) &\text{ if } p \mid h, p \nmid B_q,\\
    p^k &\text{ if } p \mid B_q.
\end{cases}
\end{equation}
Therefore, 
\[\begin{split}
F&(h,q) = q \prod_{\substack{p \nmid hB_q\\p \mid q}} \left(1 - \frac{2}{p}\right)\prod_{\substack{p \mid (h,q)\\p \nmid B_q}}\left(1 - \frac{1}{p}\right) \\
&\leq q \prod_{\substack{p \nmid hB_q\\p \mid q}} \left(1 - \frac{1}{p}\right)^2\prod_{\substack{p \mid (h,q)\\p \nmid B_q}}\left(1 - \frac{1}{p}\right) \\
&= q \left(\frac{\varphi(q)}{q}\right)^2 \prod_{p \mid (B_q,q)}\left(1 - \frac{1}{p}\right)^{-2}\prod_{\substack{p \mid (h,q)\\p \nmid B_q}}\left(1 - \frac{1}{p}\right)^{-1}
\\
&= q \left(\frac{\varphi(q)}{q}\right)^2 \prod_{p \mid (B_q,q)}\left(1 + \frac{1}{p}\right)^2\prod_{\substack{p \mid (h,q)\\p \nmid B_q}}\left(1 + \frac{1}{p}\right)\underbrace{\prod_{p \mid (B_q,q)}\left(1 - \frac{1}{p^2}\right)^{-2}\prod_{\substack{p \mid (h,q)\\p \nmid B_q}}\left(1 - \frac{1}{p^2}\right)^{-1}}_{\le\zeta(2)^2}
\\
&\ll q \left(\frac{\varphi(q)}{q}\right)^2 \prod_{p \mid (B_q,q)}\left(1 + \frac{1}{p}\right)^2\prod_{\substack{p \mid (h,q)\\p \nmid B_q}}\left(1 + \frac{1}{p}\right) 
\\
&\leq q \left(\frac{\varphi(q)}{q}\right)^2 \prod_{p \mid (B_q,q)}\left(1 + \frac{1}{p}\right)^2
\prod_{p \mid (h,q)}\left(1 + \frac{1}{p}\right)\,.
\end{split}
\]
Writing $\mu$ for the M\"obius function and using the fact that $\mu(d)^2=1$ if $d$ is a product of distinct primes or $1$, and $0$ otherwise, we get that 
\begin{align*}
\sum_{1 \leq h \leq x}\;\prod_{p \mid (h,q)}\left(1 + \frac{1}{p}\right) & \;=
\sum_{1 \leq h \leq x}\; \sum_{d \mid (h,q)}\frac{\mu^2(d)}{d}=
\sum_{d \mid q}\frac{\mu^2(d)}{d}\sum_{\substack{1 \leq h \leq x\\d \mid h}}1\\[2ex]
&\le\; x\sum_{d \mid q}\frac{\mu^2(d)}{d^2} \ll x\,.
\end{align*}
Therefore, combining this with the above estimate for $F(h,q)$ gives that
\begin{equation}\label{bound_fhq}\begin{split}
\sum_{1\leq h \leq x}F(h,q) &\ll qx \left(\frac{\varphi(q)}{q}\right)^2 \prod_{p \mid (B_q,q)}\left(1 + \frac{1}{p}\right)^2\,,
\end{split}
\end{equation}
where the implied constant is absolute. Plugging this into \eqref{u_bound} and using \eqref{T_integral} and the assumption of \eqref{decay_large_psi} as well as \eqref{eq103}, we get
\begin{align}
\frac{T}{q} & \ll \psi(q)^2\left(\frac{\varphi(q)}{q}\right)^2 \prod_{p \mid (q,B_q)}\left(1 + \frac{1}{p}\right)^2\nonumber\\[2ex]
& \le\psi(q)^2\left(\frac{\varphi(q)}{q}\right)^2 \prod_{p \mid (q,B_q)}\left(1 + \frac{1}{p-1}\right)^2\nonumber\\[2ex]
&= \left(\frac{\psi(q)\varphi(q,B_q)}{q}\right)^2 = o\left(\frac{\psi(q)\varphi(q,B_q)}{q}\right)\label{bound_on_T}
\end{align}
as $q\to\infty$ and $\psi(q)\ge1/2$. Hence,  \eqref{extract_T} implies \eqref{large_psi_lower_bound}
and so \eqref{eq_measure_large_psi_} follows.

In order to show \eqref{ud_large_psi}, let $0 \leq x<y\leq 1$ be fixed. The upper bound is again straightforward and obtained by the union bound:
\begin{align}\label{upper_M}
\frac{\mu(E^{\cI}_q \cap [x,y])}{\dfrac{2\psi(q)\varphi(q,B_q)}{q}(y-x)} &\;\leq\; \frac{\#\left\{a \in \mathcal{I}_q: a \in \Big[x - \tfrac{\psi(q)}{q}, y + \tfrac{\psi(q)}{q}\Big]\right\}}{\varphi(q,B_q)}
= 1 + o(1),    
\end{align}
where we used \eqref{show_ud} and $\frac{\psi(q)}{q} \to 0$. Thus it remains to prove the matching lower bound. To this end, let $\psi(q)\ge1/2$ and define
\begin{align*}
M :&= M(q,x,y) := \#\left\{a \in \mathcal{I}_q: a \in \Big[x + \tfrac{\psi(q)}{q}, y - \tfrac{\psi(q)}{q}\Big]\right\}.
\end{align*}
We observe that 
\begin{equation}\label{lower_M}\frac{M}{\varphi(q,B_q)(y-x)} = 1 - o(1)\end{equation}
by using the same arguments that led to \eqref{upper_M}.
Using the same reasoning as in \eqref{extract_T}, we also obtain that 
\[\mu(E^{\cI}_q \cap [x,y]) \geq 2M\frac{\psi(q)}{q} - \frac{T}{q},\]
where $T$ is as in \eqref{T_integral}. Using \eqref{bound_on_T}, we conclude that
\begin{align*}
\mu(E^{\cI}_q \cap [x,y]) &\geq  2M\frac{\psi(q)}{q} - o\left(\frac{\psi(q)\varphi(q,B_q)}{q}\right)\,.
\end{align*}
Finally, combining this with \eqref{lower_M}, \eqref{upper_M} and \eqref{eq_measure_large_psi_}, we deduce \eqref{ud_large_psi}.
\end{proof}

\subsection{Proof of Lemmas \ref{key_lemma} and \ref{key_lemma_large_psi}}\label{sec5.4}

\subsubsection{Basic overlap estimates}

\begin{lemma}[Overlap estimate for small $\psi$]\label{overlap}
Let $q \neq r$. Further, let $A,B$ be coprime integers, $B>0$ and let
 \[E_q^* := E_{q,A,B,\psi}^* := \bigcup_{\substack{a \in \mathbb{Z}_{Bq}^{*}\\ a \equiv \,A\!\!\!\! \pmod B}} \left[\frac{a}{Bq} - \frac{\psi(q)}{q}, \frac{a}{Bq} + \frac{\psi(q)}{q}\right].\]
Writing
\begin{equation}\label{vbX}
X =X(q,r) :=  2\max\left\{\frac{\psi(q)}{q},\frac{\psi(r)}{r}\right\} B\lcm(q,r),
\end{equation}
we have that
\begin{equation}
\mu(E_q^* \cap E_r^*) \ll \begin{cases}
    \;\;\;0 &\text{ if } X < 1,\\[4ex]
     B\dfrac{\psi(q)\varphi(q,B)}{q} \dfrac{\psi(r)\varphi(r,B)}{r} \prod\limits_{\substack{p \mid \frac{qr}{(q,r)^2}\\ p > X}}\left(1 + \frac{1}{p}\right) &\text{ if } 1 \leq X \leq B^3,\\[5ex]
\phantom{B}\dfrac{\psi(q)\varphi(q,B)}{q} \dfrac{\psi(r)\varphi(r,B)}{r} \prod\limits_{\substack{p \mid \frac{qr}{(q,r)^2}\\ p > X}}\left(1 + \frac{1}{p}\right) &\text{ if } X > B^3,
\end{cases}\label{est1}
\end{equation}
where all implied constants are absolute.
In particular, if
$0 < \psi(q),\psi(r) < 1/2$ then
\begin{equation}
\mu(E_q^* \cap E_r^*) \ll \begin{cases}
    \;\;\;0 &\text{ if } X < 1,\\[4ex]
     B \mu(E_q^*) \mu( E_r^*)  \prod\limits_{\substack{p \mid \frac{qr}{(q,r)^2}\\ p > X}}\left(1 + \frac{1}{p}\right) &\text{ if } 1 \leq X \leq B^3,\\[5ex]
\phantom{B}\mu(E_q^*) \mu( E_r^*)  \prod\limits_{\substack{p \mid \frac{qr}{(q,r)^2}\\ p > X}}\left(1 + \frac{1}{p}\right) &\text{ if } X > B^3\,.
\end{cases}\label{est2}
\end{equation}
\end{lemma}

\medskip

The following statement is a special case of Lemma~\ref{overlap}.

\begin{lemma}[Overlap estimate for large $\psi$]\label{overlap_large_psi}
    Let $q \neq r$ with $\psi(q),\psi(r) > \frac{1}{2}$ and $A,B,E_q^*,E_r^*,X$ as in Lemma \ref{overlap}. Then we have
    \begin{equation}
    \mu(E_q^* \cap E_r^*) \ll \begin{cases}     B\;\dfrac{\psi(q)\varphi(q,B)}{q}\; \dfrac{\psi(r)\varphi(r,B)}{r} &\text{ if } 0 \leq X \leq B^3,\\[3ex]
\dfrac{\psi(q)\varphi(q,B)}{q}\; \dfrac{\psi(r)\varphi(r,B)}{r}  &\text{ if } X > B^3.
\end{cases}
\end{equation}
    
\end{lemma}

\begin{proof}[Proof of Lemma \ref{overlap_large_psi} modulo Lemma~\ref{overlap}]

Note that $\psi(q),\psi(r) > 1/2$ implies that $X > \frac{(qr)^{1/2}}{(q,r)}$, which also ensures $X \ge 1$. Therefore, using Mertens' Theorem gives that    
\[\prod_{\substack{p \mid \frac{qr}{(q,r)^2}\\ p > X}}\left(1 + \frac{1}{p}\right)
\leq \prod\limits_{\frac{(qr)^{1/2}}{(q,r)} \leq p \leq\frac{qr}{(q,r)^2}}\left(1 + \frac{1}{p}\right) \ll 1,
\]
The required bound then follows from \eqref{est1} immediately.
\end{proof}

In the proof of Lemma \ref{overlap}, we will follow the notation of Pollington and Vaughan \cite{pv} and refine the method by allowing additional congruence conditions $A \pmod B$. To begin with, define
\[
\delta := \min\left(\frac{\psi(q)}{q},\frac{\psi(r)}{r}\right)\,,\qquad
\Delta := \max\left(\frac{\psi(q)}{q},\frac{\psi(r)}{r}\right)\,.
\]
Then, by writing
$$
q = \prod\limits_{p \in \Prime} p^u\qquad\text{and}\qquad r = \prod\limits_{p \in \Prime} p^v\,,
$$
where $u=u(q,p)$ and $v=v(r,p)$ are integer exponents and $\Prime$ is the set of all primes, we then define the following three integers:
\[n = \prod_{\substack{p \in \Prime \\ u \neq v}} p^{\max (u,v)},\quad m = \prod_{\substack{p \in \Prime \\ u \neq v}} p^{\min (u,v)},\quad \ell = \prod_{\substack{p \in \Prime \\ u = v}}p^{\max (u,v)}.\]
Observe that, since $q\neq r$, we have that $n>1$.

Further, note that for fixed $a \in \mathbb{Z}_{Bq},b \in \mathbb{Z}_{Br}$, the intervals
$$
\left[\frac{a}{Bq} - \frac{\psi(q)}{q}, \frac{a}{Bq} + \frac{\psi(q)}{q}\right],\left[\frac{b}{Br} - \frac{\psi(r)}{r}, \frac{b}{Br} + \frac{\psi(r)}{r}\right]
$$
are disjoint unless
$\lvert \frac{a}{Bq} - \frac{b}{Br}\rvert \leq 2\Delta$
and clearly,

\[\mu\left(\left[\frac{a}{Bq} - \frac{\psi(q)}{q}, \frac{a}{Bq} + \frac{\psi(q)}{q}\right] \cap \left[\frac{b}{Br} - \frac{\psi(r)}{r}, \frac{b}{Br} + \frac{\psi(r)}{r}\right]\right) \leq 2\delta.\]
Thus we obtain

\[\mu(E_q^* \cap E_r^*) \leq 2\delta \mathop{\sum_{\substack{a \in \mathbb{Z}_{Bq}^{*}\\ a \equiv A \mod B,}}\sum_{\substack{b \in \mathbb{Z}_{Br}^{*}\\ b \equiv A \mod B}}}_{\lvert \frac{a}{Bq} - \frac{b}{Br} \rvert\leq 2\Delta} 1.\]
Observe that
\[\frac{a}{Bq} - \frac{b}{Br} = \frac{c}{B\lcm(q,r)}\]
for some $c \in \mathbb{Z}$, thus we deduce
\begin{equation}\label{overlap_with_hc}\mu(E_q^* \cap E_r^*) \ll
\delta \sum_{\lvert c \rvert \leq 2\Delta B\lcm(q,r) = X} H(c),
\end{equation}
where
\begin{equation}\label{H(c)}
H(c) := \# \left\{(a,b) \in \mathbb{Z}_{Bq}^{*}\times \mathbb{Z}_{Br}^{*} : \begin{array}{l} 
a \equiv b \equiv A \pmod B,\\[0ex] 
\dfrac{a}{Bq} - \dfrac{b}{Br} = \dfrac{c}{B\lcm(q,r)}
\end{array}
\right\}.
\end{equation}

\medskip

\begin{lemma}\label{Hc_lemma}
Let $H(c),A,B$ and $m,\ell,n$ be as above. Then we have the following:
\begin{itemize}
    \item  If  $c \not\equiv \frac{r-q}{(q,r)}A \pmod B$ or $(c,n) > 1$, then $H(c) = 0$.
    \item In any case,
\begin{equation}\label{number_of_solutions} \hspace*{-5ex}H(c) \leq (q,r) \frac{\varphi(m)}{m}\frac{\varphi(\ell)^{2}}{\ell^2} \prod_{p \mid (B,m)}\left(1 - \frac{1}{p}\right)^{-1} \prod_{\substack{ p \mid (\ell,c)}}\left(1 - \frac{1}{p}\right)^{-1}\prod_{p \mid (B,\ell)}\left(1 - \frac{1}{p}\right)^{-2}.  \end{equation}
\end{itemize}
\end{lemma}

\begin{proof}
Suppose that $(a,b)$ satisfies the conditions of \eqref{H(c)}. Then 
\begin{equation}\label{rewrite_c}
c = \frac{r}{(q,r)}a - \frac{q}{(q,r)}b\,,
\end{equation}
and since $a \equiv b \equiv A \pmod B$, the condition
 $c \equiv \frac{r-q}{(q,r)}A \pmod B$ follows immediately.
Therefore, we must have that $H(c)=0$ if $c \not\equiv \frac{r-q}{(q,r)}A \pmod B$.

Now let $p \mid n$. Then $p^u \| q$ and $p^v \| r$ with $u \neq v$ and we can assume without loss of generality that $u > v$. Since $(a,Bq)=1$ implies $p \nmid a$ and therefore, we have
$p \nmid \frac{r}{(q,r)}a$. On the other hand,
$p \mid \frac{q}{(q,r)}b$ and thus by \eqref{rewrite_c}, we have
$p \nmid c$.
Therefore, we must have that $H(c)=0$ if $(c,n)>1$.

\medskip

In order to prove \eqref{number_of_solutions}, suppose that $(a_0,b_0)$ is a solution to
\begin{equation}\label{congruence_solutions}\frac{a}{Bq} - \frac{b}{Br} = \frac{c}{B\lcm(q,r)}\,, \qquad a \equiv b \equiv A \pmod{B}\,.
\end{equation}

By the Chinese Remainder Theorem, all other solution $(a,b)$ to the left hand side of \eqref{congruence_solutions} modulo $(Bq,Br)$ is given by
\[(a,b)=\left(a_0 + j \frac{q}{(q,r)},b_0 + j \frac{r}{(q,r)}\right), \quad
j \in \mathbb{Z}_{B (q,r)}.\]
The right hand side of \eqref{congruence_solutions} means 
\[a_0 + j \frac{q}{(q,r)} \equiv b_0 + j \frac{r}{(q,r)} \equiv A \pmod B\,.\]
Therefore, since  $a_0 \equiv b_0 \equiv A \mod B$, we get that
\begin{equation}
B\mid j \frac{q}{(q,r)}\qquad\text{ and }\qquad  B\mid j \frac{r}{(q,r)}\,.
\end{equation}
Since $\frac{q}{(q,r)}$ is coprime to $\frac{r}{(q,r)}$, $B \mid j$ follows.
Thus, modulo $(Bq,Br)$, all solutions to \eqref{congruence_solutions} are of the form
\[\left(a_0 + jB \frac{q}{(q,r)},b_0 + jB \frac{r}{(q,r)}\right), \quad
j \in \mathbb{Z}_{(q,r)}.\]

Write $(q,r) = \prod\limits_{p \mid (q,r)} p^w$, where we obviously have that $w=w(p,q,r) = \min(u,v)$.
Then, by the Chinese Remainder Theorem, we have that
\begin{equation}\label{vb129}
H(c) \leq
\prod_{p \mid (q,r)} \# \left\{
j \in \Z_{p^w}: 
\big(a_0 + jB \tfrac{q}{(q,r)},p\big) = 1,\; \big(b_0 + jB \tfrac{r}{(q,r)},p\big) = 1
\right\}.
\end{equation}

To estimate the factors in this estimate for each fixed $p \mid (q,r)$, we distinguish the following three cases. Note that $p\mid(q,r)$ if and only if  either $p\mid m$ or $p\mid\ell$.

\bigskip

\noindent\textbf{Case 1:} If $p \mid B$, then we use the trivial upper bound $p^w$.

\bigskip

\noindent\textbf{Case 2:} If $p \nmid B$ and $p \mid m$, then 
either
$\frac{Bq}{(q,r)}$ or $\frac{Br}{(q,r)}$ is coprime to $p$,
and we obtain at most $\varphi(p^w)$ many relevant solutions $j$ modulo $p^w$. Thus $\varphi(p^w)$ is a bound for the corresponding factor in \eqref{vb129}.

\bigskip

\noindent\textbf{Case 3:} If $p \nmid B$ and $p \mid \ell$, then
both $\frac{Bq}{(q,r)}$ and $\frac{Br}{(q,r)}$ are coprime to $p$.
Thus there are $\varphi(p^w)$ many possibilities mod $p^w$ for the condition $(a_0 + jB \frac{q}{(q,r)},p^w) = 1$ to hold.
Observe that after rearranging, we obtain
\begin{equation}
\label{hc_case3}
\frac{r}{(q,r)}\left(a_0 + jB \frac{q}{(q,r)}\right) - \frac{q}{(q,r)}\left(b_0
 + jB \frac{r}{(q,r)}\right)
\equiv c \pmod p.
\end{equation}
We need to distinguish two sub-cases: $p\mid c$ and $p\nmid c$: 

\bigskip

\noindent\textbf{Case 3a:} If $p \mid c$ then
\[
\big(a_0 + jB \tfrac{q}{(q,r)},p^w\big) = 1 \quad\iff\quad \big(b_0 + jB \tfrac{r}{(q,r)},p^w\big) = 1,
\]
thus we obtain $\varphi(p^w)$ solutions modulo $p^w$. 

\bigskip

\noindent\textbf{Case 3b:} If $p\nmid c$, then 
$$
(a_0 + jB \tfrac{q}{(q,r)},p^w) > 1\quad\Rightarrow\quad (b_0 + jB \tfrac{r}{(q,r)},p^w) = 1
$$
and similarly 
$$
(b_0 + jB \tfrac{r}{(q,r)},p^w) > 1\quad\Rightarrow\quad 
(a_0 + jB \tfrac{q}{(q,r)},p^w) = 1\,.
$$
Therefore, we obtain
$p^w - 2p^{w-1}=
p^w(1 - \frac{2}{p})
$ solutions modulo $p^w$.\\

Combining the above estimates, we obtain
that
\begin{align*}
H(c)
&\leq
\prod_{\substack{p \mid B\\p \mid (q,r)}} p^w \prod_{\substack{p \nmid B \\
p \mid m}} p^w (1 - 1/p) \prod_{\substack{p \nmid B \\
p \mid (\ell,c)}} p^w (1 - 1/p) \prod_{\substack{p \nmid B \\
p \mid \ell\\ p \nmid c}} p^w (1 - 2/p)
\\&\leq (q,r) \frac{\varphi(m)}{m}\frac{\varphi(\ell)^{2}}{\ell^2} \prod_{p \mid (B,m)}\left(1 - \frac{1}{p}\right)^{-1} \prod_{\substack{ p \mid (\ell,B)}}\left(1 - \frac{1}{p}\right)^{-2}
\prod_{\substack{p \nmid B \\ p \mid (\ell,c)}}\left(1 - \frac{1}{p}\right)^{-1}.
\end{align*}
which completes the proof of Lemma~\ref{Hc_lemma}.
\end{proof}

\subsubsection{Proof of Lemma~\ref{overlap}}\label{sec6.8.2}

Applying Lemma \ref{Hc_lemma} to \eqref{overlap_with_hc}, we obtain for $t := \frac{r-q}{(q,r)} A \pmod B$ that

\begin{align}
\nonumber\mu(E_q^* \cap E_r^*) \ll
\delta \times(q,r) \frac{\varphi(m)}{m}\frac{\varphi(\ell)^2}{\ell^2}\!\!\! &\prod_{\substack{p \mid (B,m)}}\left(1 + \frac{1}{p-1}\right)  \prod_{\substack{p \mid (B,\ell)}} \left(1 + \frac{1}{p-1}\right)^2\times \\[2ex]
&\times\sum_{\substack{\lvert c \rvert \leq X\\
(c,n) = 1,\\ c \equiv t \pmod B}} \prod_{\substack{p \mid (\ell,c)\\p\nmid B}}  \left(1 + \frac{1}{p-1}\right).\label{overlap_1}
\end{align}

To treat the sieve estimate, we follow the proof of Pollington and Vaughan \cite{pv}.
If $X < 1$, then the sum is empty, since $n>1$ and so $(c,n)=1$ implies that $c\neq0$, and in particular $|c|\ge X$ in this case. Thus $\mu(E_q^* \cap E_r^*) = 0$ for $X<1$.

If $1 \le X \le B^3$, we simply drop, from \eqref{overlap_1}, all the condition related to $B$. Thus
\[\sum_{\substack{\lvert c \rvert \leq X\\
(c,n) = 1,\\ c \equiv t \pmod B}} \prod_{\substack{p \nmid B \\
p \mid (\ell,c)}}  \left(1 + \frac{1}{p-1}\right) \leq
\sum_{\substack{\lvert c \rvert \leq X\\
(c,n) = 1}}\prod_{p\mid (\ell,c)}  \left(1 + \frac{1}{p-1}\right),\]
with the right-hand-side being identical to the corresponding term in \cite{pv}. Thus, following verbatim we get that

\[\sum_{\substack{\lvert c \rvert \leq X\\
(c,n) = 1}}\prod_{p\mid (\ell,c)}  \left(1 + \frac{1}{p-1}\right)
\ll X  \frac{\varphi(n)}{n} \prod_{\substack{p \mid n\\ p > X}} \left(1 + \frac{1}{p-1}\right).
\]

Plugging this into \eqref{overlap_1} and using \eqref{vbX}, we obtain for $1 \leq X \leq B^3$ that
 \begin{align}
 \nonumber\mu(E_q^* \cap E_r^*) &
\ll B \Delta \lcm(q,r)\delta (q,r) \frac{\varphi(m)}{m}\frac{\varphi(\ell)^2}{\ell^2}\frac{\varphi(n)}{n} \prod_{\substack{p \mid (m,B)}}\left(1 + \frac{1}{p-1}\right) \\
\nonumber&\hspace*{18ex} \times\prod_{\substack{p \mid (\ell,B)}} \left(1 + \frac{1}{p-1}\right)^2
\prod_{\substack{p \mid n\\ p > X}}
\left(1 + \frac{1}{p-1}\right)\\
\nonumber&= B \frac{\varphi(q)\varphi(r)}{qr}\psi(q)\psi(r)\prod_{\substack{p \mid (m,B)}}\left(1 + \frac{1}{p-1}\right) \prod_{\substack{p \mid (\ell,B)}} \left(1 + \frac{1}{p-1}\right)^2\\
\nonumber&\hspace*{42ex}\times \prod_{\substack{p \mid n\\ p > X}}
\left(1 + \frac{1}{p-1}\right)\\
&\leq B\dfrac{\psi(q)\varphi(q,B)}{q} \dfrac{\psi(r)\varphi(r,B)}{r}\prod_{\substack{p \mid n\\ p > X}}
\left(1 + \frac{1}{p-1}\right).\label{estim_small_X}
 \end{align}
In the penultimate equation, we used the facts that
$$
(q,r) \lcm(q,r) = m\ell^2 n = qr\quad\text{and}\quad \varphi(q)\varphi(r) = \varphi(m)\varphi(\ell)^2\varphi(n)\,
$$
while in the last equation of \eqref{estim_small_X} we used \eqref{eq103} and
  \eqref{sum_of_measures}, as well as the equation
 \begin{align}
\nonumber\prod_{\substack{p \mid (m,B)}}\left(1 + \frac{1}{p-1}\right)  & \prod_{\substack{p \mid (\ell,B)}} \left(1 + \frac{1}{p-1}\right)^2  \prod_{\substack{p \mid (n,B)}}\left(1 + \frac{1}{p-1}\right)\\[2ex]
&= \prod_{\substack{p \mid (B,r)}} \left(1 + \frac{1}{p-1}\right)\prod_{\substack{p \mid (B,q)}}\left(1 + \frac{1}{p-1}\right).     \label{B_identity}
  \end{align}
Thus, we are left to establish the bound for the case $X > B^3$.
Note that
\[
 \prod_{\substack{p \mid (\ell,c)\\p\nmid B}}  \left(1 + \frac{1}{p-1}\right)=
 \sum_{\substack{d\mid (\ell,c)\\(d,B) = 1}} \frac{\mu(d)^2}{\varphi(d)}.
\]
Recall that $(\ell,n) = 1$. Then substituting the above expression  into \eqref{overlap_1} and changing the order of summation we obtain that 
\begin{align}\notag
\sum_{\substack{\lvert c \rvert \leq X\\
(c,n) = 1,\\ c \equiv t \pmod B}} \prod_{\substack{p \mid (\ell,c)\\p\nmid B}}  \left(1 + \frac{1}{p-1}\right)
&= \sum_{\substack{\lvert c \rvert \leq X\\
(c,n) = 1,\\ c \equiv t \pmod B}} \sum_{\substack{d\mid (\ell,c)\\(d,B) = 1}} \frac{\mu(d)^2}{\varphi(d)}\\[2ex]
&\stackrel{|c|=db}{=} 2\sum_{\substack{d \leq X \\ d \mid \ell \\ (d,B) = 1}}
\;\sum_{\substack{\lvert b \rvert \leq \frac{X}{d} \\ (b,n) = 1 \\ b \equiv td^{-1} \pmod B}} \frac{\mu(d)^2}{\varphi(d)}
\nonumber\\[2ex]
&\ll \sum_{\substack{d \leq X\\ (d,B) = 1}}d^{-9/10}
\;\sum_{\substack{\lvert b \rvert \leq \frac{X}{d} \\ (b,n) = 1 \\ b \equiv td^{-1} \pmod B}} 1.
\label{exchanging_summation}
\end{align}
If $d > \sqrt{X}$, we will use the trivial bound
\begin{equation}\label{trivial_sieve}
    \sum_{\substack{\lvert b\rvert \leq \frac{X}{d} \\ (b,n) = 1 \\ b \equiv td^{-1} \pmod B}} 1 \ll \frac{X}{d}.
\end{equation}
Now suppose that $d \leq \sqrt{X}$. Then
write $P(n) = \prod\limits_{\substack{p \mid n \\ p\nmid B}} p$, and observe that $(t,n) = 1$ and thus, we have that
\[
\left\{\begin{array}{l}
     (b,n) = 1  \\[0.5ex]
     bd \equiv t \pmod B
\end{array}\right.
 \quad\Longleftrightarrow\quad 
\left\{\begin{array}{l}
(b,P(n)) = 1\\[0.5ex]
bd \equiv t \pmod B.
\end{array}\right.
\]
The direction $\Rightarrow$ is obvious, so assume by contradiction that there exists
 $b,d \in \mathbb{Z}$ and $p \in \Prime$ such that
 $p \mid (b,n), p \nmid P(n), db \equiv t \pmod B$. This implies that $p\mid n$ and $p \mid B$ thus
 $t \equiv db \equiv 0 \pmod p$, which means that $(t,n) > 1$, a contradiction since
 every $p \mid n$ divides exactly one of $\frac{q}{(q,r)}$ and $\frac{r}{(q,r)}$.

Note that, by construction, $(P(n),B)= 1$. Therefore, applying a standard result from sieve theory (see e.g. \cite[Theorem 3.8]{mv}), we get for $d \leq \sqrt{X}$ and $B \leq X^{1/3}$ that

\begin{align}
\nonumber\sum_{\substack{\lvert b\rvert \leq \frac{X}{d} \\ (b,n) = 1 \\ b \equiv td^{-1} \pmod B}}  1
&\ll \frac{X/d}{B}\prod_{\substack{p \mid P(n)\\ p \leq \sqrt{2X/Bd}}}
\left(1 -\frac{1}{p}\right)\\[2ex]
\nonumber&\ll \frac{X/d}{B}\prod_{\substack{p \mid P(n)\\ p \leq \sqrt2 X^{1/6}}}
\left(1 -\frac{1}{p}\right) \\[2ex]
&\ll \frac{X/d}{B}\prod_{\substack{p \mid P(n)\\ p \leq X}}
\left(1 -\frac{1}{p}\right),\label{eq138}
\end{align}
where in the last inequality we used the following form of Mertens' Theorem: 
$$
\prod\limits_{\sqrt2 X^{1/6} \leq p \leq X} \left(1 -\frac{1}{p}\right) \gg 1\,.
$$
Plugging \eqref{trivial_sieve} and \eqref{eq138} into \eqref{exchanging_summation}, we obtain that
\begin{align}
\nonumber\sum_{\substack{d \leq X\\ (d,B) = 1}}d^{-9/10}
\sum_{\substack{\lvert b \rvert \leq \frac{{X}}{d} \\ (b,n) = 1 \\ b \equiv td^{-1} \pmod B}} 1
&\ll
X\sum_{\sqrt{X} < d \leq {X}} \frac{1}{d^{19/10}} +
\frac{X}{B}\prod_{\substack{p \mid P(n)\\ p \leq X}}
\left(1 -\frac{1}{p}\right)\sum_{d \leq {\sqrt{X}}} \frac{1}{d^{19/10}}
\\&\ll X^{11/20}+ \frac{X}{B}\prod_{\substack{p \mid P(n)\\ p \leq X}}
\left(1 -\frac{1}{p}\right)\nonumber\\[2ex]
&\ll  \frac{X}{B}\prod_{\substack{p \mid P(n)\\ p \leq X}}
\left(1 -\frac{1}{p}\right),\label{eq139}
\end{align}
where once again we used  Mertens' Theorem and the inequality $B \leq X^{1/3}$.
Finally, note that
\[\begin{split}\prod_{\substack{p \mid P(n)\\ p \leq X}}
\left(1 -\frac{1}{p}\right)
&= \prod_{\substack{p \mid n\\p \leq X}}\left(1 -\frac{1}{p}\right)
\prod_{\substack{p \mid (n,B)\\p\leq X} } \left(1 -\frac{1}{p}\right)^{-1}
\\[2ex]
&\leq \frac{\varphi(n)}{n} \prod_{\substack{p \mid n\\ p > X}}
\left(1 + \frac{1}{p-1}\right) \prod_{\substack{p \mid (n,B)}}
\left(1 + \frac{1}{p-1}\right).
\end{split}
\]
Putting this into \eqref{eq139}  we get that
$$
\sum_{\substack{d \leq X\\ (d,B) = 1}}d^{-9/10}
\sum_{\substack{\lvert b \rvert \leq \frac{X}{d} \\ (b,n) = 1 \\ b \equiv td^{-1} \pmod B}} 1\ll
\frac{X}{B}\frac{\varphi(n)}{n} \prod_{\substack{p \mid n\\ p > X}}
\left(1 + \frac{1}{p-1}\right) \prod_{\substack{p \mid (n,B)}}
\left(1 + \frac{1}{p-1}\right)\,.
$$
Then, by \eqref{overlap_1},   \eqref{B_identity}, \eqref{exchanging_summation} and \eqref{sum_of_measures}, we obtain, for $X \geq B^3$, that
 \[
 \mu(E_q^* \cap E_r^*) \ll \dfrac{\psi(q)\varphi(q,B)}{q} \dfrac{\psi(r)\varphi(r,B)}{r} \prod_{\substack{p \mid n\\ p > X}}\left(1 + \frac{1}{p-1}\right).
 \]
 Since $p \mid n$ if and only if $ p \mid \frac{qr}{(q,r)^2}$
 and 
 $$\prod\limits_{p \in \Prime} \frac{\left(1 + \frac{1}{p-1}\right)}{\left(1 + \frac{1}{p}\right)} \;=\prod\limits_{p \in \Prime} \left(1+\frac{1}{p^2-1}\right) \le \exp\left(\sum_{p\in\Prime}{\frac{1}{p^2-1}}\right) \ll 1,
 $$
 this finishes the proof of Lemma \ref{overlap}.

\subsubsection{Second moment bounds: proof of Lemmas~\ref{key_lemma} and \ref{key_lemma_large_psi}}

Let \[\Psi(\mathcal{S}) := \sum_{q \in \mathcal{S}} \frac{\psi(q)\varphi(q)}{q}\] and recall that by the assumptions of Lemmas~\ref{key_lemma} and \ref{key_lemma_large_psi}, we have that $\Psi(\mathcal{S}) \geq B^8$. Apart from Lemmas~\ref{overlap} and \ref{overlap_large_psi}, in the proofs given in this section we will make use of the following results from \cite{HSW} which, in turn, improve upon estimates obtained in  \cite{abh,MR4125453}. 

\medskip

\begin{proposition}[{\cite[Proposition 6.2]{HSW}}]\label{few_large_gcd}
    For any $s,\varepsilon > 0$, let 
    \[\mathcal{E}_s := 
    \left\{(q,r) \in \mathcal{S}^2: X(q,r) \leq \frac{B\Psi(\mathcal{S})}{s}\;\right\}.
    \]
    Then we have that
    \[\sum_{(q,r) \in \mathcal{E}_s} \frac{\psi(q)\varphi(q)}{q}\frac{\psi(r)\varphi(r)}{r}
    \ll_{\varepsilon} \frac{\Psi(\mathcal{S})^2}{s^{1-\varepsilon}}.
    \]
\end{proposition}

\bigskip

\begin{proposition}[{\cite[Proposition 6.3]{HSW}}]\label{few_large_gcd_large_error}
    Let $s,t,C \geq 0$ be fixed, $\eta \in (0,1/2)$ and
    \[
    \mathcal{E}_{s,t,C} := 
    \left\{(q,r) \in \mathcal{S}^2: X(q,r) \leq sB \Psi(\mathcal{S}),\; L_t(q,r) \geq C\right\}\,,
    \]
where 
$$
L_t(q,r):=\sum_{\substack{p\mid \frac{qr}{(q,r)^2}\\[0.5ex] p\ge t}}\frac1p\,.
$$
Then we have that
    \[\sum_{(q,r) \in \mathcal{E}_{s,t,C}} \frac{\psi(q)\varphi(q)}{q}\frac{\psi(r)\varphi(r)}{r}
    \ll_{\eta} \Psi(\mathcal{S})^2s^{\frac{1}{2} + \eta}e^{-(1 - \eta)tC}.
    \]
\end{proposition}

\medskip

Note that the above propositions are stated and proved in \cite{HSW} with $\cS=\{1,\dots,N\}$. However, the versions we provide above are easily deduced by redefining $\psi$ to be $\psi\cdot\mathds{1}_{\cS}$ for our choice of $\cS\subset\{1,\dots,N\}$ for a suitably large $N$. When comparing Propositions~\ref{few_large_gcd} and \ref{few_large_gcd_large_error} to their counterparts in \cite{HSW}, note that $X(q,r)$ and $\Psi(\cS)$ are the equivalents to  $D(q,r)$ and $\Psi(N)$ used in  \cite{HSW}, and are related by equations $X(q,r)=2B D(q,r)$ and
$\Psi(N)=2\Psi(\cS)$. Also the parameter $s$ we use in Propositions~\ref{few_large_gcd} and \ref{few_large_gcd_large_error} should be understood as $4$ times the parameter $s$ used in Propositions~6.2 and 6.3 of \cite{HSW} in order to see a perfect match between Propositions~\ref{few_large_gcd} and \ref{few_large_gcd_large_error} and their counterparts in \cite{HSW}.

\medskip

\begin{proof}[Proof of Lemma \ref{key_lemma}]
We start with the case where $\psi: \N \to [0,1/2]$. Note that if $B$ is absolutely bounded,
Lemma \ref{overlap} becomes the estimate
\begin{align*}
\mu(E_q^* \cap E_r^*) &\ll \frac{\psi(q)\varphi(q,B)}{q}\frac{\psi(r)\varphi(r,B)}{r}\prod_{\substack{p \mid \frac{qr}{(q,r)^2}\\
p > \Delta \lcm(q,r)}} \left(1 + \frac{1}{p}\right)\\[2ex]
&\ll_{B} \frac{\psi(q)\varphi(q)}{q}\frac{\psi(r)\varphi(r)}{r}\prod_{\substack{p \mid \frac{qr}{(q,r)^2}\\
p > \Delta \lcm(q,r)}} \left(1 + \frac{1}{p}\right)
\end{align*}
of Pollington and Vaughan \cite{pv}. With this estimate, Koukouloupolos and Maynard \cite{MR4125453} showed that
\[
\sum_{q,r \in \mathcal{S}}\frac{\psi(q)\varphi(q)}{q}\frac{\psi(r)\varphi(r)}{r}\prod_{\substack{p \mid \frac{qr}{(q,r)^2}\\
p > \Delta \lcm(q,r)}} \left(1 + \frac{1}{p}\right) \ll 1
\]
whenever $1 \leq \sum_{q \in \mathcal{S}} \frac{\psi(q)\varphi(q)}{q} \leq 2$, thus the statement follows when $B$ is bounded.

Thus from now one we can assume that $B$ is larger than any prescribed constant.
We will adapt the approach of \cite[p. 217 -- 220]{abh} that is built upon the ideas of \cite{MR4125453}. 

Specifically, we partition $\mathcal{S}^2:=\cS\times\cS$ into the following five subsets:
\begin{align*}
    \mathcal{E}^{(1)} &:= \{(q,r) \in \mathcal{S}^2: q = r\},\\[1ex]
    \mathcal{E}^{(2)} &:= \{(q,r) \in \mathcal{S}^2: X(q,r) \leq \sqrt{\Psi(\mathcal{S})},
L_X(q,r) \leq \Log B\},\\[1ex]
\mathcal{E}^{(3)} &:= \left\{(q,r) \in \mathcal{S}^2: X(q,r) \leq \sqrt{\Psi(\mathcal{S})},
L_X(q,r) > \Log B
\right\},\\[1ex]
\mathcal{E}^{(4)} &:= \left\{(q,r) \in \mathcal{S}^2: X(q,r) \geq \sqrt{\Psi(\mathcal{S})}, L_X(q,r) \leq 10
\right\},\\[1ex] 
\mathcal{E}^{(5)} &:= \left\{(q,r) \in \mathcal{S}^2: X(q,r) \geq \sqrt{\Psi(\mathcal{S})}, L_X(q,r) > 10
\right\}.
\end{align*}

We will frequently use the following estimates that follow from \eqref{sum_of_measures} and the well known facts $\frac{n}{\varphi(n)} \ll \log \log n$ and $\log(1 + x) \leq x$:
\begin{equation}
\label{freq_mertens}\mu(E_q^*) \ll \log \log B\frac{\psi(q)\varphi(q)}{q}, \quad
\prod_{\substack{p \mid \tfrac{qr}{(q,r)^2}\\ p > t}} \left(1 + \frac{1}{p}\right) \leq \exp(L_t(q,r)),\; t \in \R.
\end{equation}
Clearly, the contribution from $\mathcal{E}^{(1)}$ to the left hand side of \eqref{concl001} is negligible since
$$
\sum_{q\in\cS} \mu(E_q^*)\ge \sum_{q\in\cS} \frac{\varphi(q)\psi(q)}{q}\ge B^8\ge1
$$
and therefore
\[\sum_{(q,r) \in \mathcal{E}^{(1)}} \mu(E_q^* \cap E_r^*)
= \sum_{q\in\cS} \mu(E_q^*) \le \left(\sum_{q\in\cS} \mu(E_q^*)\right)^2.
\]

To estimate the contribution from $ \mathcal{E}^{(2)}$, using Lemma \ref{overlap} and \eqref{freq_mertens}, we obtain that
\begin{align*}
\sum_{(q,r) \in \mathcal{E}^{(2)}} \mu(E_q^* \cap E_r^*)
&\ll B^2 \sum_{(q,r) \in \mathcal{E}^{(2)}}
\mu(E_q^*)\mu(E_r^*)
\\& \ll B^2 (\Log \Log B)^2 \sum_{(q,r) \in \mathcal{E}^{(2)}}
\frac{\psi(q)\varphi(q)}{q}\frac{\psi(r)\varphi(r)}{r}.
\end{align*}
Furthermore observe that
$\mathcal{E}^{(2)} \subseteq \mathcal{E}_{s}$ with $s = \frac{\sqrt{\Psi(\mathcal{S})}}{B} \geq B^3$, thus an application of Proposition \ref{few_large_gcd} with $\varepsilon = 1/8$ yields

\begin{align*}
    B^2(\Log \Log B)^2 \sum_{(q,r) \in \mathcal{E}^{(2)}}
\frac{\psi(q)\varphi(q)}{q}\frac{\psi(r)\varphi(r)}{r}
&\ll \Psi(\mathcal{S})^2 \frac{B^{5/2}}{B^{5/2}} \leq \Psi(\mathcal{S})^2\\[2ex]
&\ll \left(\sum_{q \in \mathcal{S}} \mu(E^*_q)\right)^2\,.
\end{align*}

\medskip

Next, we estimate the contribution from $\mathcal{E}^{(3)}$. Writing $f(j) := \exp(\exp(j))$, let $j(q,r)$ be the maximal integer satisfying $L_{f(j)}(q,r) \geq 5$. 
Observe that

\begin{align}
\nonumber 5 > L_{f(j+1)}(q,r) \geq L_X(q,r) - \sum_{p \leq f(j+1)} \frac{1}{p} & \geq \Log B - O\left(\Log\Log f(j+1)\right)\\[1ex] 
&\geq \Log B - O(j)\label{estimate_mertens}\,.
\end{align}
Therefore $j \gg \Log B$ and 
$L_{X}(q,r) = O(j)$. Hence, by Lemma \ref{overlap}, we get that
\begin{align}
\nonumber\sum_{\substack{(q,r) \in \mathcal{E}^{(3)}\\ j(q,r) = j}}
\mu(E_q^* \cap E_r^*) &\ll \sum_{\substack{(q,r) \in \mathcal{E}^{(3)}\\ j(q,r) = j}}
B \mu(E_q^*) \mu( E_r^*)  \exp(L_{X}(q,r))
\\&\leq B(\Log\Log B)^2 \exp(O(j)) \sum_{\substack{(q,r) \in \mathcal{E}^{(3)}\\ j(q,r) = j}} \frac{\psi(q)\varphi(q)}{q}\frac{\psi(r)\varphi(r)}{r}.\label{E2_estim}
\end{align}
Furthermore, we see that

\[\{(q,r) \in \mathcal{E}^{(3)}: j(q,r) = j\} \subseteq \mathcal{E}_{1,f(j),5},\]
so an application of Proposition \ref{few_large_gcd_large_error} with $\eta = \tfrac{1}{4}$  gives

\[\sum_{\substack{(q,r) \in \mathcal{E}^{(3)}\\ j(q,r) = j}} \frac{\psi(q)\varphi(q)}{q}\frac{\psi(r)\varphi(r)}{r}
\ll \frac{\Psi(\mathcal{S})^2}{\exp(2f(j))}.\\
\]
Substituting this into \eqref{E2_estim} and then summing over $j \gg \log B$ yields
\begin{equation}\label{conclusion_(2)}\begin{split}
\sum_{(q,r) \in \mathcal{E}^{(3)}
\mu(E_q^* \cap E_r^*)} &\ll  \Psi(\mathcal{S})^2 B^2 \sum_{j \gg \Log B} \frac{\exp(O(j))}{\exp(2f(j))}
\\&\ll \Psi(\mathcal{S})^2 B^2 \sum_{j \gg \Log B}\frac{1}{\exp(\exp(\exp(j)))}
\\[2ex]
&\ll \Psi(\mathcal{S})^2
\ll \left(\sum_{q \in \mathcal{S}} \mu(E^*_q)\right)^2\,
\end{split}
\end{equation}
as required.

\medskip

Regarding the contribution from $\mathcal{E}^{(4)}$ and $\mathcal{E}^{(5)}$, observe that 
$X \geq \sqrt{\Psi(\mathcal{S})} \geq B^4$ implies that Lemma \ref{overlap} gives 
\begin{equation}\label{good_overlap_bound}\mu(E_q^* \cap E_r^*) \ll \mu(E_q^*) \mu( E_r^*)  \prod\limits_{\substack{p \mid \frac{qr}{(q,r)^2}\\ p > X}}\left(1 + \frac{1}{p}\right).\end{equation}

Therefore, in the case of $\mathcal{E}^{(4)}$ we obtain straightforwardly that
\[\sum_{(q,r) \in \mathcal{E}^{(4)}}
\mu(E_q^* \cap E_r^*) \ll \sum_{\substack{(q,r) \in \cS^2\\ L_{X}(q,r)\le10}}
\mu(E_q^*) \mu( E_r^*)  \exp(L_{X}(q,r)) \ll
\left(\sum_{q = 1}^Q \mu(E_q^*) \right)^2.
\]

Finally, we consider the contribution from $\mathcal{E}^{(5)}$.
As in $\mathcal{E}^{(3)}$, let $j(q,r)$ be the maximal integer that satisfies $L_{f(j)}(q,r) \geq 5$. 
Since $L_X(q,r) \geq 10$, an application of Mertens' theorem in the form 
$$
\sum_{f(j) \leq p \leq f(j+1)} \frac{1}{p} \leq 2
$$ 
for $j$ large enough\footnote{We can assume $j$ large enough since $B^3 \leq X \leq F(\exp(\exp(j+1)))$ and $B$ can be assumed to be large enough.} implies that
 $\sqrt{\Psi(\cS)} \leq X \leq f(j)$, which shows that
$j \geq \Log \Log \sqrt{\Psi(\mathcal{S})}$.
Note that for $j(q,r) = j$, we have that
\[L_X(q,r) \leq \sum_{p \mid \frac{qr}{(q,r)^2}}
\frac{1}{p} \leq L_{f(j+1)} + \sum_{p \leq f(j+1)} \frac{1}{p} = O(j).
\]
We thus obtain by \eqref{good_overlap_bound} that 
\[
\begin{split}
\sum_{\substack{(q,r) \in \mathcal{E}^{(5)}:\\ j(q,r)  = j}}
\mu(E_q^*\cap E_r^*)
&\ll \sum_{\substack{(q,r) \in \mathcal{E}^{(5)}:\\ j(q,r)  = j}}\mu(E_q^*)\mu(E_r^*)
\exp(L_{X}(q,r))
\\[1ex]
&\ll B\sum_{\substack{(q,r) \in \mathcal{E}^{(5)}\\ j(q,r)  = j}}
\frac{\psi(q)\varphi(q)}{q}\frac{\psi(r)\varphi(r)}{r}  \exp(O(j)).
\end{split}
\]
Since $j(q,r) = j$ implies that $X \leq f(j)$
and, by construction, $L_{f(j)} \geq 1$, we can apply Proposition \ref{few_large_gcd_large_error} to \[\{(q,r) \in \mathcal{E}^{(5)}: j(q,r) = j\} \subseteq \mathcal{E}_{1,f(j),5}\,.\]
Hence, by Proposition \ref{few_large_gcd_large_error} with $\eta = \tfrac{1}{4}$ we get, similarly to \eqref{conclusion_(2)}, that
\begin{equation*}\begin{split}\sum_{\substack{(q,r) \in \mathcal{E}^{(5)}}}
 \mu(E_q^*\cap E_r^*)
 &\ll
\Psi(\mathcal{S})^2 B \sum_{j \geq  \Log \Log \sqrt{\Psi(\mathcal{S})}} \frac{\exp(O(j))}{\exp(\tfrac{5}{4}f(j))}\\[1ex]
&\ll
\Psi(\mathcal{S})^2 B\sum_{j \geq \Log \Log \sqrt{\Psi(\mathcal{S})}} \frac{1}{\exp(\exp(\exp(j)))}
 \\[1ex]
 &\ll \Psi(\mathcal{S})^2\ll \left(\sum_{q \in \mathcal{S}} \mu(E^*_q)\right)^2\,
 \end{split}
 \end{equation*}
where we used $\Psi(\mathcal{S}) \geq B^8$. This proves Lemma \ref{key_lemma}.
\end{proof}

\bigskip

\begin{proof}[Proof of Lemma~\ref{key_lemma_large_psi}]

By Lemma~\ref{overlap_large_psi}, we have that
 \begin{equation*}
    \mu(E_q^* \cap E_r^*) \ll \begin{cases}
B\;\dfrac{\psi(q)\varphi(q,B)}{q} \dfrac{\psi(r)\varphi(r,B)}{r} &\text{ if } 0 \leq X \leq B^3,\\[3ex]
\dfrac{\psi(q)\varphi(q,B)}{q} \dfrac{\psi(r)\varphi(r,B)}{r}  &\text{ if } X > B^3.
\end{cases}
\end{equation*}

It suffices to consider an analogue to $\mathcal{E}^{(1)},\mathcal{E}^{(2)}$ which after another application of Proposition~\ref{few_large_gcd} gives the result also for the case of large functions $\psi$ considered in Lemma~\ref{key_lemma_large_psi}.
\end{proof}

\section{Shrinking Targets  \label{STlove} }

Let $(X,\cA,\mu,T)$ be a measure-preserving dynamical system.  Recall, that by definition $\mu$ is a probability measure. Now let $\{A_n \}_{n \in \mathbb{N}} $ be a sequence of subsets  in $\cA$ and let
\begin{eqnarray}
\label{eq35}
W\big(T,\{A_n\}\big) & := & \limsup_{n\to\infty}T^{-n}(A_n)  \nonumber  \\[1ex] & = &  \{x\in   X: T^n(x)\in A_n\ \hbox{ for infinitely many }n\in \mathbb{N}\} \, .
\end{eqnarray}
For obvious reasons the sets $A_n$ can be thought of as targets that the orbit under $T$ of points in $X$ have to hit. The interesting situation is usually, when working within a metric space,  the diameters of $A_n$ tend to zero as $n$ increases. It is thus natural to refer to $ W\big(T,\{A_n\}\big)$ as the corresponding shrinking target set associated with the given dynamical system and target sets. Since $T$ is measure-preserving we have that  $ \mu (T^{-n}(A_n))  =  \mu (A_n)$ for all $n$, and  a straightforward consequence of the (convergence) Borel--Cantelli Lemma is that
\begin{equation} \label{appconv}
\mu\big(W(T,\{A_n\})\big)=0   \qquad {\rm if \ } \qquad \sum_{n=1}^\infty \mu(A_n)  \, < \, \infty  \, .
\end{equation}

\noindent  Now two natural questions arise.   Both fall under the umbrella of the ``shrinking target problem''  formulated in \cite{MR1309976}.
\begin{itemize}
  \item[] \vspace*{-1ex}
  \begin{itemize} \item[~ \textbf{(P1)}] What is the $\mu$-measure of $W(T,\{A_n\})$ if the measure sum in \eqref{appconv} diverges? \\[-2ex]
  \item[~ \textbf{(P2)}] What is the Hausdorff dimension of $W(T,\{A_n\})$  if the measure sum converges and so $\mu\big(W(T,\{A_n\})\big)=0$?
      \end{itemize}
      \vspace*{-1ex}
\end{itemize}

\noindent  To be precise, the target sets $A_n$ in the original formulation in  \cite{MR1309976} are restricted to balls $B_n$.     The more general setup naturally incorporates a larger class of problems. For example, within the context of simultaneous  Diophantine approximation,  it enables us to address problems associated with the weighted (the target sets are rectangular) and multiplicative (the target sets are hyperbola) theories -- see \cite[\S1.2: Remark~4]{MR4572386}.

\bigskip

In this paper we concentrate specifically on the shrinking target  problem \textbf{(P1)}.   With this in mind,  following the terminology from \cite{MR1826488}, we say  that the sequence $\{A_n \}_{n \in \mathbb{N}}$ of target sets  is a {\em Borel--Cantelli sequence} (or simply a {\em BC sequence})     if
$$
\mu \Big(W\big(T,\{A_n\}\big)\Big) =1 \, ;
$$ 
that is, for $\mu$-almost every point $x\in X$ we have that $T^n(x)\in A_n$ holds for infinitely many $n \in \N$. In this case one also says that $\{A_n \}_{n \in \mathbb{N}}$ satisfies the dynamical Borel--Cantelli lemma. Obviously any BC sequence must be divergent.
In the following section \S\ref{stgen}, we investigate conditions under which a sequence of target sets $\{A_n \}_{n \in \mathbb{N}}$ is a BC sequence.
In the subsequent section  \S\ref{stball}, we  consider  the special setup in which the target sets $A_n$ are balls and we obtain an analogue of Khintchine's theorem for shrinking targets.

\subsection{Dynamical Borel--Cantelli for general  shrinking targets \label{stgen}}

BC sequences have been investigated for many decades. It is  known \cite[Proposition~1.6]{MR1826488},  that in any measure preserving systems $(X,\cA,\mu,T)$ with $\mu$  non-trivial  (i.e., there are sets $A \in \cA$ with $0 < \mu(A) < \mu(X)=1$)  there are divergent non-BC sequences $\{A_n\}$. BC sequences have been found under various additional assumptions on the measure preserving systems and/or the sequences $\{A_n\}$ in question. BC sequences have also been used to characterise the properties of dynamical systems. In particular, we have the following:

\begin{proposition}[See {\cite[Proposition~1.5]{MR1826488}}]\label{prop1}
Let $(X,\cA,\mu,T)$ be any measure preserving systems. Then
\begin{itemize}
  \item[{\rm(i)}] $T$ is ergodic $\iff$ every constant sequence $A_n\equiv A$, $\mu(A)>0$, is BC;
  \item[{\rm(ii)}] $T$ is weakly mixing $\iff$ every sequence $\{A_n\}$ that only contains finitely many distinct sets, none of them of measure zero, is BC;
\item[{\rm(iii)}] $T$ is lightly mixing $\iff$ every divergent sequence that only contains finitely many distinct sets, possibly of measure zero, is BC.
\end{itemize}
\end{proposition}

\bigskip

\noindent Recall that the system $(X,\cA,\mu,T)$ is said to be
\begin{itemize}
  \item {\em ergodic} if for any $A\in\cA$ such that $T^{-1}(A)=A$ we have that $\mu(A)=0$ or $1$. It is well know and easily verified that $(X,\cA,\mu,T)$ is ergodic if and only if
\begin{equation}\label{ergodic}
\forall\;A,B\in\cA\;\; \text{with $\mu(A)\mu(B)>0$ }\exists\;n\in\N\;\;\text{such that }\mu(A\cap T^{-n}(B))>0\,.
\end{equation}
  \item  {\em mixing} (or {\em strong mixing}) if for any $A,B\in\cA$ with $\mu(B)>0$
\begin{equation}\label{mixing}
\lim_{n\to\infty}\frac{\mu(A\cap T^{-n}(B))}{\mu(B)}=\mu(A)\,.
\end{equation}
  \item {\em weakly mixing} if
\begin{equation}\label{weaklymixing}
\forall\;A,B\in\cA\;\; \lim_{N\to\infty}\frac{1}{N}\sum_{n=1}^N  \, \left|\mu(A\cap T^{-n}(B))-\mu(A)\mu(B)\right|=0\,.
\end{equation}
\item {\em lightly mixing} if
$$
\liminf_{n\to\infty}\mu(A\cap T^{-n}(B))>0\qquad\text{if}\quad \mu(A)\mu(B)>0\,.
$$
\end{itemize}

\noindent It is also well known \cite[Theorem~1.1]{MR2090763} that
$$
\text{mixing $\Rightarrow$ lightly mixing $\Rightarrow$ weakly mixing $\Rightarrow$ ergodic}
$$
and that all the four notions are different.  Furthermore, it is know that
$(X,\cA,\mu,T) $  is weakly mixing   if and only if for every $A,B\in\cA$ with $\mu(B)>0$ there is a subset $\cN=\cN(A,B)\subset\N$ of zero  density  such that
\begin{equation}\label{wekalymixing2}
\lim_{n\to\infty,\;n\not\in\cN}\frac{\mu(A\cap T^{-n}(B))}{\mu(B)}=\mu(A)\,.
\end{equation}

We note that conditions  \eqref{mixing} and \eqref{wekalymixing2} characterising mixing and weakly mixing systems can be weakened by replacing `$=$' with `$\le$'.
Formally, we have the following

\begin{proposition}
\label{T7}
Let $(X,\cA,\mu,T)$ be a measure preserving system. Then 
\begin{itemize}
    \item[{\rm(i)}] 
the system is
mixing if and only if for any $A,B\in\cA$ with $\mu(B)>0$
\begin{equation}\label{ecq08++}
\limsup_{n\to\infty}\frac{\mu(A\cap T^{-n}(B))}{\mu(B)}\le\mu(A)\,.
\end{equation}
\item[{\rm(ii)}] the system is weakly mixing if and only if for any $A,B\in\cA$ with $\mu(B)>0$ there is a subset $\cN=\cN(A,B)\subset\N$ of zero density such that
\begin{equation}\label{eq41}
\limsup\limits_{\substack{n\to\infty\\ n\not\in\cN}}\frac{\mu(A\cap T^{-n}(B))}{\mu(B)}\le\mu(A).
\end{equation}
\item[{\rm(iii)}] 
the system is ergodic if and only if for any $A,B\in\cA$ with $\mu(B)>0$ there is a subset $\cN=\cN(A,B)\subset\N$ such that $\N\setminus\cN$ is infinite and
\begin{equation}\label{vvv}
\limsup\limits_{\substack{n\to\infty\\ n\not\in\cN}}\frac{\mu(A\cap T^{-n}(B))}{\mu(B)}\le\mu(A)\,,
\end{equation}
or, equivalently, if and only if for any $A,B\in\cA$ with $\mu(B)>0$
\begin{equation}\label{vvv++}
\liminf\limits_{n\to\infty}\frac{\mu(A\cap T^{-n}(B))}{\mu(B)}\le\mu(A)\,.
\end{equation}

\end{itemize}
\end{proposition}

\begin{proof}
In (i) and (ii) the necessity is obvious, so we only need to prove the sufficiency. Regarding (i), 
applying \eqref{ecq08++} with $B$ replaced $X\setminus B$ we find that
$$\limsup_{n\to\infty}\mu(A\setminus T^{-n} (B))   = \limsup_{n\to\infty}  \mu\big(A\cap T^{-n}(X\setminus B) \big) \leq \mu(A)(1 - \mu(B))\,  $$  and so it follows that 
\begin{eqnarray*}
\liminf_{n\to\infty} \mu(A\cap T^{-n} (B))  & = & \mu(A) - \limsup_{n\to\infty} \mu(A\setminus T^{-n} (B))  \\[1ex]
&  \geq  &  \mu(A)\big(1 - (1 - \mu(B))\big) = \mu(A)\,  \mu(B)  \, . 
\end{eqnarray*}
This together  with \eqref{ecq08++} implies \eqref{mixing} as desired. The sufficiency in (ii) is verified using an identical argument with $n$ running through a subsequence of $\N$, namely through the set $\N\setminus\big(\cN(A,B)\cup\cN(A,X\setminus B)\big)$. Observe that $\cN(A,B)\cup\cN(A,X\setminus B)$ will have zero asymptotic density as long as $\cN(A,B)$ and $\cN(A,X\setminus B)$ are of zero asymptotic density.

Now we turn to (iii). Unfortunately the above argument cannot be used as the set
$\cN(A,B)\cup\cN(A,X\setminus B)$ may become finite, or even empty, even if both $\N\setminus\cN(A,B)$ and $\N\setminus\cN(A,X\setminus B)$ are infinite.
First, we prove the sufficiency, and so we assume that \eqref{vvv} holds. By Proposition~\ref{prop1}(i), to demonstrate the ergodicity it is enough to show that any divergent constant sequence $\{A_n\}$ is BC; that is, $A_n=B$ for a fixed set $B \in \cA$. With Theorem~\ref{T4} in mind, let $E_n:=T^{-n}(A_n)=T^{-n}(B)$, $\delta>0$, $q_1<q_2$ be any integers and let
$A=\bigcup_{j=q_1}^{q_2}E_j$.
By \eqref{vvv}, for any $\delta > 0 $, there is a strictly increasing subsequence $n_i$ of integers such that for all large enough $i$
$$
\mu\left(A\cap T^{-n_i}(B)\right)\le (1+\delta)\mu\left(A\right)\mu(B)\,.
$$
Then
$$
\mu\left(A\cap E_{n_i}\right)=\mu\left(A\cap T^{-n_i}(B)\right)\le (1+\delta)\mu\left(A\right)\mu(B)=
(1+\delta)\mu\left(A\right)\mu(E_{n_i})\,.
$$
Also, note that
$$
\inf\{\mu(E_{n_i}):i\in \N\}>c_0:=\tfrac12\mu(B)>0\,.
$$
Thus, the conditions of Theorem~\ref{T4} are satisfied with $\cI=\{n_i:i\in\N\}$, depending on $q_1$ and $q_2$, and we have that $\mu(\limsup E_n)=1$, or in other words, $\{A_n\}$ is BC.

\medskip

Now, we prove the necessity in (iii); that is  to say that ergodicity  implies \eqref{vvv++}.  To start with,   recall that in view of  Birkhoff's Ergodic Theorem, for  any $ B \in \cA$,
\begin{equation}\label{svnoparis}
\lim_{n \to \infty}  \frac1n \sum_{i=0}^{n-1}\mathds{1}_B(T^ix)    =  \mu (B)
\end{equation}
for $\mu$-almost all $ x \in X$.   Now note that for  any $ A \in \cA$ and $i \in \N$, we have that 
$$
\int_{A} \mathds{1}_B(T^ix) d\mu(x) = \int_{A} \mathds{1}_{T^{-i}(B)} (x) d\mu(x)  =  \mu \big( A \cap T^{-i}(B)\big)   \, .
$$
It thus follows,  on integrating both sides of \eqref{svnoparis}, that  
$$
\lim_{n \to \infty}  \frac1n \sum_{i=0}^{n-1} \mu \big( A \cap T^{-i}(B)\big)   = \mu(A)  \mu (B) \, . 
$$
The  upshot of this (see for example the Stolz–Ces\'{a}ro theorem) is that
$$
\liminf\limits_{n\to\infty}  \mu \big( A \cap T^{-n}(B)\big)    \le \mu(A)  \mu (B) \le  \limsup\limits_{n\to\infty}  \mu \big( A \cap T^{-n}(B)\big)  \, . 
$$ 
This left hand side is precisely \eqref{vvv++} and we are done.
\end{proof}

Recall that Proposition~\ref{prop1} provides a range of criteria for sequences to be BC. These sequences assume only a finite number of distinct  sets. Now we will try to understand the conditions sufficient for sequences that may contain infinitely many distinct elements to be BC.
We start with a dynamical version of Theorem~\ref{T4}, which rests on the following uniform version of property \eqref{ecq08++}: 
\begin{itemize}
\item[]
\begin{itemize}
  \item[{\bf(U1)}] for any $A\in\cA$ there exists a positive sequence $\{\ve_n\}$ and $n_0\in\N$, which may depend on $A$, such that $\lim_{n\to\infty}\ve_n=0$ and for all $n\ge n_0$ and any $B\in\cA$ we have that
\begin{equation}\label{uniformlymixing0vb}
\mu(A\cap T^{-n}(B))\le \mu(A)\mu(B)+\ve_n\mu(B)\,.
\end{equation}
\end{itemize}
\end{itemize}

\bigskip 

\noindent We remark that condition \textbf{(U1)} can equivalently be stated in a more compact form as follows: for any $A\in\cA$ 
\begin{equation}\label{vb158}
\limsup_{n\to\infty}\sup_{B\in\cA: \; \mu(B)>0}\frac{\mu(A\cap T^{-n}(B))}{\mu(B)}\le\mu(A)\,.
\end{equation}
Observe that  \eqref{vb158} is a uniform (in terms of $B$) version of \eqref{ecq08++}.  In view of Proposition~\ref{T7}, we know that the condition \eqref{mixing} associated with the definition of mixing and the ``weaker'' looking condition \eqref{ecq08++} are equivalent. Thus,   \textbf{(U1)} can be regarded as defining an appropriate  condition for the notion of  \textit{uniform mixing}. However, the  proof  showing  the equivalence of  \eqref{mixing} and \eqref{ecq08++}  cannot be generalised to the uniform case; that is to show that the``weaker'' looking condition
\eqref{vb158} is equivalent to  the following ``natural'' uniform version of \eqref{mixing}: for any $A\in\cA$
\begin{equation}\label{vb158++}
\lim_{n\to\infty} \ \sup_{B\in\cA: \; \mu(B)>0}\left|\frac{\mu(A\cap T^{-n}(B))}{\mu(B)}-\mu(A)\right|=0 \,.
\end{equation} In short,   uniform mixing defined via  \textbf{(U1)} is potentially  less restrictive than defining it via \eqref{vb158++}.

\bigskip 

\begin{theorem}\label{T9}
Let $(X,\cA,\mu,T)$ be a measure preserving system satisfying\/ {\rm\textbf{(U1)}}.
Then any sequence $\{A_n\}$, such that $\limsup\limits_{n\to\infty}\mu(A_n)>0$, is BC.
\end{theorem}

\begin{proof}
With Theorem~\ref{T4} in mind, let $E_n:=T^{-n}(A_n)$, $\delta>0$, $q_1<q_2$ be any integers and let
$A=\bigcup_{j=q_1}^{q_2}E_j$. By \eqref{uniformlymixing0vb}, {\bf(M1)} is verified. Since $T$ is measure preserving
$$\limsup\limits_{n\to\infty}\mu(E_n)=\limsup\limits_{n\to\infty}\mu(A_n)>0  \, . $$
Then, by Theorem~\ref{T4} it follows that  $\mu(\limsup E_n)=1$.  In other words, $\{A_n\}$ is BC.
\end{proof}

\medskip

\begin{theorem}[Dynamical version of Theorem~\ref{T1}]\label{T10}
Let $(X,\cA,\mu,T)$ be a measure preserving system satisfying\/ {\rm\textbf{(U1)}}.
Then any divergent sequence $\{A_n\}$ such that for some constant $C\ge 1$
\begin{equation}\label{uniformlymixing5vb+}
\sum_{\substack{n,m\in\N\\ n+m\le N}}\mu(A_m\cap T^{-n}(A_{m+n}))\le C\left(\sum_{n=1}^N\mu(A_n)\right)^2
\end{equation}
holds for infinitely $N\in\N$, is a BC-sequence.
\end{theorem}
\begin{proof}
Since $T$ is measure preserving, \eqref{uniformlymixing5vb+} is equivalent to
$$
\sum_{1\le m<n\le N}\mu(T^{-m}(A_m)\cap T^{-n}(A_n))\le C\left(\sum_{n=1}^N\mu(A_n)\right)^2 = C\left(\sum_{n=1}^N\mu(T^{-n}(A_n))\right)^2\,.
$$
Hence, for $N$ sufficiently large we get that
$$
\sum_{1\le m, n\le N}\mu(T^{-m}(A_m)\cap T^{-n}(A_n))\le (2C+1)\left(\sum_{n=1}^N\mu(T^{-n}(A_n))\right)^2\,,
$$
which verifies \eqref{eqn02} with $E_n:=T^{-n}(A_n)$. The fact that $\{A_n\}$ is a divergent sequence and $T$ is measure preserving verifies \eqref{eqn01}. Finally, condition {\bf(U1)} verifies {\bf(M1)}. Thus,  Theorem~\ref{T1} is applicable and we conclude that $\limsup_{n\to\infty}T^{-n}(A_n)$ has full measure, and so $\{A_n\}$ is BC.
\end{proof}

The next result can be viewed as a dynamical version of Theorem~\ref{T3}.

\begin{theorem}[Dynamical version of Theorem~\ref{T3}]   \label{dT3}
Let $(X,\cA,\mu,T)$ be a measure preserving system, and let $\{A_n\}$ be a sequence of measurable subsets.
Suppose that for any $\varepsilon>0$ there exist constants $C'>0$ and $0<\varepsilon^*<\varepsilon$ and strictly increasing sequences $M_k$ and $N_k$ such that for sufficiently large $k$
\begin{equation}\label{eq07++}
\varepsilon^*\le \sum_{M_k\le n\le N_k}\mu(A_n) \le \varepsilon
   \end{equation} \vspace*{-2ex}
and
\begin{equation} \label{eqn05++}
\sum_{n,m\ge M_k,\; m+n\le N_k} \mu\big(A_m\cap T^{-n}(A_{m+n}) \big) \ \le \  C'\,  \left(\sum_{{M_k\le n \le N_k}}\mu(A_n)\right)^2\,.
   \end{equation}
\\
Also suppose that
\begin{itemize}
\item[]
\begin{itemize}
  \item[{\bf(U2)}]
for any $A\in\cA$
\begin{equation}\label{eq08++}
\limsup_{k\to\infty}\;\dfrac{\displaystyle\sum_{M_k\le n \le N_k}\mu\left(A\cap T^{-n}(A_n)\right)}{\displaystyle\sum_{M_k\le n \le N_k}\mu(A_n)}\le\mu(A)\,.
\end{equation}
\end{itemize}
\end{itemize}
Then $\{A_n\}$ is BC.
\end{theorem}

\medskip

\noindent The proof is similar to the manner in which the previous theorem is  deduced from Theorem~\ref{T1}.  This time we use Theorem~\ref{T3} with  $\cS_k=[M_k,N_k]\cap\Z$.  We leave the details to the reader.

\medskip

We bring this section to a close  by obtaining a  quantitative  version of  Theorem~\ref{T10} in line with the discussion centered around \eqref{asympBC} in \S\ref{introA}.  In short, the quantitative dynamical Borel--Cantelli Lemma is a direct consequence of  strengthening  the  ``mixing' condition \textbf{(U1)} by assuming that the associated sequence $\{\ve_n\}$  is summable. Formally,   we say that
\emph{$\mu$ is  {\small{$\Sigma$}}-mixing}  (short for summable-mixing) if  \textbf{(U1)}  holds with the additional conditions that $\ve_n$ is not dependent on $A$ and satisfies
$$\textstyle{\sum_{n=1}^\infty } \ \ve_n   < \infty \, . $$
Note that this is a less restrictive notion of {\small{$\Sigma$}}-mixing  than that introduced in \cite[\S1.1:~Defintion~1]{MR4572386} or that of uniform mixing introduced in \cite{MR2262783}.  The key difference is that in both of these, the one side inequality \eqref{uniformlymixing0vb} is replaced by the  assumption  that
$$
\big| \mu(A\cap T^{-n}(B)) -  \mu(A)\mu(B) \big| \le \ve_n\mu(B)   \,.
$$
Nevertheless, the less restrictive form suffices to establish the following strengthening  of Theorem~\ref{T10},  the conclusion of which is the same as \cite[Theorem~1]{MR4572386} in which the stronger form of {\small{$\Sigma$}}-mixing is assumed.

\begin{theorem}\label{gencountthm}
Let $(X,\cA,\mu,T)$ be a measure preserving system  and suppose that   $\mu$ is  {\small{$\Sigma$}}-mixing.   Let $\{A_n\}$  be a sequence of subsets in $\cA$.
  Then, for any given $\varepsilon>0$, we have that
\begin{equation} \label{def_Psi}\# \big\{ 1\le n \le   N :    T^{n}(x) \in A_n \}=\Phi(N)+O\left(\Phi^{1/2}(N) \ (\log\Phi(N))^{3/2+\varepsilon}\right)
\end{equation}
\noindent for $\mu$-almost all $x\in X$, where
$$\Phi(N):=\textstyle{\sum_{n=1}^N} \  \mu(A_n) \,  .
$$
\end{theorem}

\medskip

\noindent A simple consequence of Theorem~\ref{gencountthm} is that if the
measure sum diverges (that is $\{A_n\}$ is a divergent sequence),  then
$$\lim_{N \to \infty} \# \big\{ 1\le n \le   N :    T^{n}(x) \in A_n \} = \infty  $$ for  $\mu$--almost all $x \in X$   and so it clearly implies Theorem~\ref{T10}.   Regarding the proof of Theorem~\ref{gencountthm}, it is exactly the same as the proof appearing  in  \cite[\S2]{MR4572386}. The only observation to make while working through \cite[\S2]{MR4572386}  is that we  only require the weaker form of {\small{$\Sigma$}}-mixing as defined in this paper.  For the sake of completeness, we mention that {\small{$\Sigma$}}-mixing implies that \eqref{uniformlymixing5vb+} holds for $C$ arbitrarily close to one.

\subsection{Dynamical Borel--Cantelli for moving shrinking balls}\label{stball}

With the general shrinking target setup in mind, suppose in addition that  $X$ is  a metric space and  $\mu$ is a measure such that every Borel subset is $\mu$-measurable.  Furthermore, given a real, positive function $\psi:\N\to\R_{\ge 0}$ and $n \in \N$,   let  $A_n:=B_n=B\big(x_n,\psi(n)\big)$ be  a  ball in $X$  centred at the point $x_n \in X$ and radius $\psi(n)$, and let
$$
W(T, \psi):=W\big(T,\{B_n\}\big)=\limsup_{n\to\infty}T^{-n}\big(B(x_n,\psi(n))\big)
$$
denote the   associated shrinking target set \eqref{eq35}. Suppose for the moment that the balls $B_n$  have a fixed centre and so are not ``moving''. In other words,  we fix some point $x_0 \in  X$ and take $x_n = x_0$ for all $n \in \N$.   Thus, $x\in W(T, \psi)$ means that $T^n(x)\in B(x_0,\psi(n))$ for infinitely many $n \in \N$. When $\psi(n)\to0$ as $n\to\infty$, it means that $T^n(x)$ revisits the shrinking neighbourhood $B(x_0,\psi(n))$ of $x_0$ again and again. If $x_0$ can be taken to be an arbitrary point of $X$, the shrinking target problem of  determining  when  $W(T, \psi)$ has full $\mu$-measure leads to a quantitative refinement of the density of the $T$-orbit of  $\mu$-almost every point $x\in X$.  The following  result, that can be regarded as a dynamical version of Theorem~\ref{T1MS} for balls  or equivalently a balls version of Theorem~\ref{T10}, provides natural  conditions under which the shrinking target set for moving balls  has full measure.

\begin{theorem}\label{T12}
Let $(X,\cA,\mu,T)$ be a measure preserving system on a  metric space $X$ equipped with a doubling Borel regular measure $\mu$, satisfying\/
\begin{itemize}
\item[]
\begin{itemize}
  \item[{\bf(S1)}] for any ball $A\in\cA$ there exists a positive sequence $\{\ve_n\}$ and $n_0\in\N$, which may depend on $A$, such that $\lim_{n\to\infty}\ve_n=0$ and for all $n\ge n_0$ and any ball $B\in\cA$ we have that
\begin{equation}\label{eq45}
\mu(A\cap T^{-n}(B))\le \mu(A)\mu(B)+\ve_n\mu(B)\,.
\end{equation}
\end{itemize}
\end{itemize}
Let $\psi:\N\to\R_{>0}$, $B_n:=B(x_n,\psi(n))$, and suppose that
\begin{itemize}
\item[]
\begin{itemize}
  \item[{\rm(i)}] $\{B_n\}$ is divergent, that is $\sum_{n=1}^\infty\mu(B_n)=\infty$, and  \\ \vspace*{-2ex}
  \item[{\rm(ii)}] for some constant $C>1$ \eqref{uniformlymixing5vb+} holds for infinitely many $N\in\N$.
\end{itemize}
\end{itemize}
Then $\{B_n\}$ is a BC-sequence; that is, $\mu(W(T, \psi))=1$.
\end{theorem}

\noindent Clearly,  property  \textbf{(S1)} imposed in the above theorem  corresponds precisely to  \textbf{(U1)} in which the measurable sets are  restricted to balls.  We now impose further conditions on the measure which enables us to the determine the measure of balls purely in terms of their radii and thus derive from Theorem~\ref{T12} an  analogue of Khintchine's  theorem for moving shrinking targets.

Let $X$ be a metric space equipped with a non-atomic Borel-regular probability measure $\mu$. We will assume for simplicity that $\supp(\mu)=X$. Suppose there
exist constants $ \delta > 0$, $0<a\le 1\le b<\infty$ and $r_0 > 0$ such
that
\begin{equation} \label{MTPmeasure}
 a \, r ^{\delta}  \ \le    \  \mu(B)  \ \le    \   b \, r
^{\delta}
\end{equation}
for any ball $B=B(x,r)$ with $x\in X$ and radius $r\le r_0$. Such a measure is said to be \emph{$\delta$-Ahlfors regular}. It is well known that if $X$ supports a $\delta$-Ahlfors regular measure $\mu$, then
 \begin{equation} \label{ARdim}
\dim_{\rm H} X = \delta  \,
\end{equation}
where $ \dim_{\rm H} X $ denotes the Hausdorff dimension of $X$.
Moreover, we have that $ \mu$ is comparable  to the $\delta$-dimensional Hausdorff measure ${\cal H}^\delta$ and so   \eqref{MTPmeasure} is valid with $\mu$ replaced by  ${\cal H}^{\delta} $ -- see \cite{F,MR1800917,
MAT} for the details.  The following {\em Khintchine-type} (KT) result   is a direct consequence of  the definition  of $\delta$-Ahlfors regular measures,  the convergence statement  \eqref{appconv} and Theorem~\ref{T12}.  It is natural to refer to the corresponding sequence  of target sets as a \emph{Khintchine-type sequence} or simply a \emph{KT sequence}. 

\begin{corollary}
Within the setting of Theorem~\ref{T12}, assume in  addition that $\mu$ is $\delta$-Ahlfors regular on $X$. Then
$$
\mu(W(T, \psi))=\left\{\begin{array}{cc}
                      0 & \text{if }\sum_{n=1}^{\infty} \psi(n)^\delta<\infty\,,\\[2ex]
                      1 & \text{if }\sum_{n=1}^{\infty}\psi(n)^\delta=\infty\,.
                                          \end{array}
\right.
$$
\end{corollary}

\appendix
\section{Proofs of Lemmas~DBC and GDBC}

We refer the reader to \cite{MR4497313} for further historical remarks on Lemma~DBC, which proof is essentially a consequence of the Cauchy-Schwarz inequality and can be found, for example, in \cite{MR1672558, MR1764799, MR1265493, MR548467}.

For completeness and to demonstrate the difference of our approach to the `standard' proofs of Divergence Borel--Cantelli lemmas here we give a full proof of Lemma~GDBC as well as a proof of the implication
$$
\text{Lemma~GDBC $\Longrightarrow$ Lemma~DBC.}
$$
The complementary implication is straightforward and has been discussed after the  statement of Lemma~GDBC.

\begin{proof}[Proof of Lemma~GDBC]
For each $k\in\N$ define
$$
E(k)=\bigcup_{i\in\cS_k}E_i\qquad\text{and}\qquad f_k:=\sum_{i\in\cS_k}\mathbf{1}_{E_i}\,.
$$
Using the Cauchy-Schwarz inequality we get that
\begin{equation}\label{eqn06+}
\left(\int_{E(k)}f_k\right)^2
\le \left(\int_{E(k)}1^2\right)\left(\int_{E(k)}f_k^2\right)=\mu(E(k))\left(\sum_{s,t\in\cS_k} \mu\big(E_s\cap E_t \big)\right)\,.
\end{equation}
Observe that $\int_{E(k)} f_k=\sum_{i\in\cS_k}\mu(E_i)$. Then, rearranging \eqref{eqn06+} and using \eqref{eqn04} and \eqref{eqn05} gives
$$
\mu(E(k))~\ge~ \frac{\displaystyle\left(\sum_{i\in\cS_k}\mu(E_i) \right)^2}{\displaystyle2\sum_{\substack{s<t\\[0.5ex] s,t\in\cS_k}} \mu\big(E_s\cap E_t \big) +\sum_{i\in\cS_k}\mu(E_i) }
~\stackrel{\eqref{eqn05}}{\ge}~  \frac{\displaystyle\left(\sum_{i\in\cS_k}\mu(E_i) \right)^2}{\displaystyle2C'\left(\sum_{i\in\cS_k}\mu(E_i) \right)^2  +\sum_{i\in\cS_k}\mu(E_i) }
$$
$$
\ge  \frac{1}{\displaystyle2C' +\left(\sum_{i\in\cS_k}\mu(E_i)\right)^{-1} } ~\stackrel{\eqref{eqn05}}{\ge}~   \frac{1}{2C'+a^{-1}}\,.\hspace{19ex}
$$
Therefore, $\mu\Big(\bigcup_{\ell=m}^\infty E(k_\ell)\Big)\ge (2C'+a^{-1})^{-1}$ for any $m\in\N$, where $\{k_\ell\}_{\ell\in\N}$ is the sequence of $k\in\N$ satisfying \eqref{eqn04} and \eqref{eqn05}.
Finally note that, since $\cS_k\subset\Z_{\ge k}$, we have that
$$
\limsup_{\ell\to\infty}E(k_\ell)\subset \limsup_{i\in\cS}E_i\subset E_\infty\,,
$$
and therefore, by the monotonicity and continuity of $\mu$, we get that
$$
\mu(E_\infty)\ge \mu\left(\limsup_{i\in\cS}E_i\right)\ge \mu\left(\bigcap_{m=1}^\infty\bigcup_{\ell=m}^\infty E(k_\ell)\right)=\lim_{m\to\infty}
\mu\left(\bigcup_{\ell=m}^\infty E(k_\ell)\right) \ge \frac{1}{2C'+a^{-1}}\,,
$$
as required.
\end{proof}

\begin{proof}[Lemma~GDBC\/\;$\Rightarrow$\;Lemma DBC; the proof]
Note that \eqref{eqn02} implies that
\begin{equation} \label{eqn06}
\sum_{1\le s<t\le Q} \mu\big(E_s\cap E_t \big) \ \le \  \frac{C}{2}\,  \left(\sum_{i=1}^Q\mu(E_i)\right)^2\,.
   \end{equation}
Fix any $a\in\N$. For sufficiently large $k\in\N$ choose, by \eqref{eqn01}, an integer $Q_k>k$ such that
\begin{equation}\label{eqn07+}
\sum_{i=k}^{Q_k} \mu\big(E_i\big)>a \sum_{i=1}^k \mu\big(E_i\big)>a
\end{equation}
and \eqref{eqn02}, and consequently \eqref{eqn06}, is applicable with $Q=Q_k$. Then
$$
\sum_{k\le s<t\le Q_k} \mu\big(E_s\cap E_t \big) \ \stackrel{\eqref{eqn06}}{\le}
\frac{C}{2}\,  \left(\sum_{i=1}^{Q_k}\mu(E_i)\right)^2
\stackrel{\eqref{eqn07+}}{\le} \
\frac C2\,\frac{(a+1)^2}{a^2}\,  \left(\sum_{i=k}^{Q_k}\mu(E_i)\right)^2\,.
$$
Let $\cS_k=\{k,\dots,Q_k\}\subset\Z_{\ge k}$. Then, by Lemma~GDBC,
$$
\mu(E_\infty)\ge \left(C\,\frac{(a+1)^2}{a^2}+\frac1a\right)^{-1}.
$$
Since $a$ can be made arbitrarily large, we get that
$\mu(E_\infty)\ge C^{-1}$, as required.
\end{proof}

\vspace*{6ex}

\noindent{\it Acknowledgments}.
SV  would like to take this opportunity to thank his wonderful colleagues in the number theory group at York;  in particular Evgeniy, Jason and Victor for their unreserved  support, friendship  and putting up with his ramblings for nearly twenty years -- thanks team!  Regarding this paper, thanks to  David (Simmons) and Jason for various discussions regarding shrinking targets -- it's always enjoyable and fruitful talking to you guys!  Finally, good luck to Bridget with life after sixty (still all in one piece)  and to the dynamic duo, Iona and  Ayesha,  for the start of their new adventures -- hopefully in time getting up and over  ``that great big hill of hope''!

\bibliographystyle{abbrv}
\bibliography{BCbiblio}

\begin{thebibliography}{10}

\bibitem{Adiceam_2015}
F.~Adiceam.
\newblock Rational approximation and arithmetic progressions.
\newblock {\em Int. J. Number Theory}, 11(2):451--486, 2015.

\bibitem{MR3784747}
C.~Aistleitner.
\newblock A note on the {D}uffin-{S}chaeffer conjecture with slow divergence.
\newblock {\em Bull. Lond. Math. Soc.}, 46(1):164--168, 2014.

\bibitem{abh}
C.~Aistleitner, B.~Borda, and M.~Hauke.
\newblock On the metric theory of approximations by reduced fractions: A
  quantitative {K}oukoulopoulos–{M}aynard theorem.
\newblock {\em Compos. Math.}, 159(2):207–231, 2023.

\bibitem{MR4008524}
C.~Aistleitner, T.~Lachmann, M.~Munsch, N.~Technau, and A.~Zafeiropoulos.
\newblock The {D}uffin-{S}chaeffer conjecture with extra divergence.
\newblock {\em Adv. Math.}, 356:106808, 11, 2019.

\bibitem{MR4312709}
M.~Alam, A.~Ghosh, and S.~Yu.
\newblock Quantitative {D}iophantine approximation with congruence conditions.
\newblock {\em J. Th\'{e}or. Nombres Bordeaux}, 33(1):261--271, 2021.

\bibitem{MTP-survey}
D.~Allen and E.~Daviaud.
\newblock A survey of recent extensions and generalisations of the mass
  transference principle.
\newblock {\em Preprint. arXiv:2306.15535}, 2023.

\bibitem{MR2508636}
V.~Beresnevich, V.~Bernik, M.~Dodson, and S.~Velani.
\newblock Classical metric {D}iophantine approximation revisited.
\newblock In {\em Analytic number theory}, pages 38--61. Cambridge Univ. Press,
  Cambridge, 2009.

\bibitem{MR2184760}
V.~Beresnevich, D.~Dickinson, and S.~Velani.
\newblock Measure theoretic laws for lim sup sets.
\newblock {\em Mem. Amer. Math. Soc.}, 179(846):x+91, 2006.

\bibitem{MR3101800}
V.~Beresnevich, G.~Harman, A.~Haynes, and S.~Velani.
\newblock The {D}uffin-{S}chaeffer conjecture with extra divergence {II}.
\newblock {\em Math. Z.}, 275(1-2):127--133, 2013.

\bibitem{MR3105329}
V.~Beresnevich, A.~Haynes, and S.~Velani.
\newblock Multiplicative zero-one laws and metric number theory.
\newblock {\em Acta Arith.}, 160(2):101--114, 2013.

\bibitem{MR3618787}
V.~Beresnevich, F.~Ram\'{\i}rez, and S.~Velani.
\newblock Metric {D}iophantine approximation: aspects of recent work.
\newblock In {\em Dynamics and analytic number theory}, volume 437 of {\em
  London Math. Soc. Lecture Note Ser.}, pages 1--95. Cambridge Univ. Press,
  Cambridge, 2016.

\bibitem{MTP}
V.~Beresnevich and S.~Velani.
\newblock A mass transference principle and the {D}uffin-{S}chaeffer conjecture
  for {H}ausdorff measures.
\newblock {\em Ann. of Math. (2)}, 164(3):971--992, 2006.

\bibitem{MR2457266}
V.~Beresnevich and S.~Velani.
\newblock A note on zero-one laws in metrical {D}iophantine approximation.
\newblock {\em Acta Arith.}, 133(4):363--374, 2008.

\bibitem{MR4497313}
V.~Beresnevich and S.~Velani.
\newblock The divergence {B}orel-{C}antelli {L}emma revisited.
\newblock {\em J. Math. Anal. Appl.}, 519(1):Paper No. 126750, 2023.

\bibitem{MR2090763}
V.~Bergelson and A.~Gorodnik.
\newblock Weakly mixing group actions: a brief survey and an example.
\newblock In {\em Modern dynamical systems and applications}, pages 3--25.
  Cambridge Univ. Press, Cambridge, 2004.

\bibitem{MR1324786}
P.~Billingsley.
\newblock {\em Probability and measure}.
\newblock Wiley Series in Probability and Mathematical Statistics. John Wiley
  \& Sons, Inc., New York, third edition, 1995.
\newblock A Wiley-Interscience Publication.

\bibitem{MR2944100}
M.~Boshernitzan and J.~Chaika.
\newblock Borel-{C}antelli sequences.
\newblock {\em J. Anal. Math.}, 117:321--345, 2012.

\bibitem{MR36787}
J.~W.~S. Cassels.
\newblock Some metrical theorems in {D}iophantine approximation. {I}.
\newblock {\em Proc. Cambridge Philos. Soc.}, 46:209--218, 1950.

\bibitem{MR0349591}
J.~W.~S. Cassels.
\newblock {\em An introduction to {D}iophantine approximation}, volume No. 45
  of {\em Cambridge Tracts in Mathematics and Mathematical Physics}.
\newblock Hafner Publishing Co., New York, 1972.
\newblock Facsimile reprint of the 1957 edition.

\bibitem{Catlin1976}
P.~A. Catlin.
\newblock Two problems in metric {Diophantine} approximation. {I}, {II}.
\newblock J. {Number} {Theory} 8, 282-288, 289-297 (1976)., 1976.

\bibitem{MR1826488}
N.~Chernov and D.~Kleinbock.
\newblock Dynamical {B}orel-{C}antelli lemmas for {G}ibbs measures.
\newblock {\em Israel J. Math.}, 122:1--27, 2001.

\bibitem{CT-Mem}
S.~Chow and N.~Technau.
\newblock Littlewood and {D}uffin--{S}chaeffer-type problems in {Diophantine}
  approximation.
\newblock {\em Mem. Amer. Math. Soc.}, To appear, 2023.

\bibitem{MR4859}
R.~J. Duffin and A.~C. Schaeffer.
\newblock Khintchine's problem in metric {D}iophantine approximation.
\newblock {\em Duke Math. J.}, 8:243--255, 1941.

\bibitem{F}
K.~Falconer.
\newblock {\em Fractal geometry}.
\newblock John Wiley \& Sons, Ltd., Chichester, 1990.
\newblock Mathematical foundations and applications.

\bibitem{Federer-69:MR0257325}
H.~Federer.
\newblock {\em Geometric measure theory}.
\newblock Die Grundlehren der mathematischen Wissenschaften, Band 153.
  Springer-Verlag New York Inc., New York, 1969.

\bibitem{MR2262783}
J.~L. Fern\'{a}ndez, M.~V. Meli\'{a}n, and D.~Pestana.
\newblock Quantitative mixing results and inner functions.
\newblock {\em Math. Ann.}, 337(1):233--251, 2007.

\bibitem{MR133297}
P.~Gallagher.
\newblock Approximation by reduced fractions.
\newblock {\em J. Math. Soc. Japan}, 13:342--345, 1961.

\bibitem{MR2327135}
S.~Gou\"{e}zel.
\newblock A {B}orel-{C}antelli lemma for intermittent interval maps.
\newblock {\em Nonlinearity}, 20(6):1491--1497, 2007.

\bibitem{MR2669635}
C.~Gupta, M.~Nicol, and W.~Ott.
\newblock A {B}orel-{C}antelli lemma for nonuniformly expanding dynamical
  systems.
\newblock {\em Nonlinearity}, 23(8):1991--2008, 2010.

\bibitem{Harman_1988}
G.~Harman.
\newblock Metric {D}iophantine approximation with two restricted variables.
  {I}. {T}wo square-free integers, or integers in arithmetic progressions.
\newblock {\em Math. Proc. Cambridge Philos. Soc.}, 103(2):197--206, 1988.

\bibitem{MR1672558}
G.~Harman.
\newblock {\em Metric number theory}, volume~18 of {\em London Mathematical
  Society Monographs. New Series}.
\newblock The Clarendon Press, Oxford University Press, New York, 1998.

\bibitem{MR1764799}
G.~Harman.
\newblock Variants of the second {B}orel-{C}antelli lemma and their
  applications in metric number theory.
\newblock In {\em Number theory}, Trends Math., pages 121--140. Birkh\"{a}user,
  Basel, 2000.

\bibitem{HSW}
M.~Hauke, S.~Vazquez, and A.~Walker.
\newblock Proving the {D}uffin-{S}chaeffer conjecture without {GCD} graphs.
\newblock {\em Preprint. arXiv:2404.15123}, 2024.

\bibitem{MR2915535}
A.~K. Haynes, A.~D. Pollington, and S.~L. Velani.
\newblock The {D}uffin-{S}chaeffer conjecture with extra divergence.
\newblock {\em Math. Ann.}, 353(2):259--273, 2012.

\bibitem{MR1800917}
J.~Heinonen.
\newblock {\em Lectures on analysis on metric spaces}.
\newblock Universitext. Springer-Verlag, New York, 2001.

\bibitem{MR1309976}
R.~Hill and S.~L. Velani.
\newblock The ergodic theory of shrinking targets.
\newblock {\em Invent. Math.}, 119(1):175--198, 1995.

\bibitem{MR1724857}
R.~Hill and S.~L. Velani.
\newblock The shrinking target problem for matrix transformations of tori.
\newblock {\em J. London Math. Soc. (2)}, 60(2):381--398, 1999.

\bibitem{MR4417342}
M.~Hussain, B.~Li, D.~Simmons, and B.~Wang.
\newblock Dynamical {B}orel-{C}antelli lemma for recurrence theory.
\newblock {\em Ergodic Theory Dynam. Systems}, 42(6):1994--2008, 2022.

\bibitem{MR1512207}
A.~Khintchine.
\newblock Einige {S}\"{a}tze \"{u}ber {K}ettenbr\"{u}che, mit {A}nwendungen auf
  die {T}heorie der {D}iophantischen {A}pproximationen.
\newblock {\em Math. Ann.}, 92(1-2):115--125, 1924.

\bibitem{MR2276480}
D.~H. Kim.
\newblock The dynamical {B}orel-{C}antelli lemma for interval maps.
\newblock {\em Discrete Contin. Dyn. Syst.}, 17(4):891--900, 2007.

\bibitem{MR4125453}
D.~Koukoulopoulos and J.~Maynard.
\newblock On the {D}uffin-{S}chaeffer conjecture.
\newblock {\em Ann. of Math. (2)}, 192(1):251--307, 2020.

\bibitem{KMY24}
D.~Koukoulopoulos, J.~Maynard, and D.~Yang.
\newblock An almost sharp quantitative version of the {D}uffin--{S}chaeffer
  conjecture.
\newblock {\em Preprint. arXiv:2404.14628}, 2024.

\bibitem{MR4572386}
B.~Li, L.~Liao, S.~Velani, and E.~Zorin.
\newblock The shrinking target problem for matrix transformations of tori:
  Revisiting the standard problem.
\newblock {\em Advances in Mathematics}, 421:108994, 2023.

\bibitem{MAT}
P.~Mattila.
\newblock {\em Geometry of sets and measures in {E}uclidean spaces}, volume~44
  of {\em Cambridge Studies in Advanced Mathematics}.
\newblock Cambridge University Press, Cambridge, 1995.
\newblock Fractals and rectifiability.

\bibitem{MR2214457}
F.~Maucourant.
\newblock Dynamical {B}orel-{C}antelli lemma for hyperbolic spaces.
\newblock {\em Israel J. Math.}, 152:143--155, 2006.

\bibitem{mv}
H.~L. Montgomery and R.~C. Vaughan.
\newblock {\em Multiplicative number theory I: Classical theory}.
\newblock Cambridge University Press, 2006.

\bibitem{MR4153361}
E.~Nesharim, R.~R\"{u}hr, and R.~Shi.
\newblock Metric {D}iophantine approximation with congruence conditions.
\newblock {\em Int. J. Number Theory}, 16(9):1923--1933, 2020.

\bibitem{pv}
A.~Pollington and R.~C. Vaughan.
\newblock The k ‐dimensional duffin and schaeffer conjecture.
\newblock {\em Mathematika}, 37(2):190–200, 1990.

\bibitem{MR4425845}
A.~D. Pollington, S.~Velani, A.~Zafeiropoulos, and E.~Zorin.
\newblock Inhomogeneous {D}iophantine approximation on {$M_0$}-sets with
  restricted denominators.
\newblock {\em Int. Math. Res. Not. IMRN}, 2022(11):8571--8643, 2022.

\bibitem{MR1265493}
S.~C. Port.
\newblock {\em Theoretical probability for applications}.
\newblock Wiley Series in Probability and Mathematical Statistics: Probability
  and Mathematical Statistics. John Wiley \& Sons, Inc., New York, 1994.
\newblock A Wiley-Interscience Publication.

\bibitem{MR3606945}
F.~A. Ram\'{\i}rez.
\newblock Counterexamples, covering systems, and zero-one laws for
  inhomogeneous approximation.
\newblock {\em Int. J. Number Theory}, 13(3):633--654, 2017.

\bibitem{RandomFR}
F.~A. Ram\'irez.
\newblock Khintchine's theorem with random fractions.
\newblock {\em Mathematika}, 66(1):178--199, 2020.

\bibitem{MR548467}
V.~G. Sprind\v{z}uk.
\newblock {\em Metric theory of {D}iophantine approximations}.
\newblock Scripta Series in Mathematics. V. H. Winston \& Sons, Washington,
  D.C.; John Wiley \& Sons, New York-Toronto, Ont.-London, 1979.
\newblock Translated from the Russian and edited by Richard A. Silverman, With
  a foreword by Donald J. Newman.

\bibitem{MR0095165}
P.~Sz\"usz.
\newblock \"{U}ber die metrische {T}heorie der {D}iophantischen
  {A}pproximation.
\newblock {\em Acta Math. Acad. Sci. Hungar.}, 9:177--193, 1958.

\bibitem{MR3988812}
H.~Yu.
\newblock A {F}ourier-analytic approach to inhomogeneous {D}iophantine
  approximation.
\newblock {\em Acta Arith.}, 190(3):263--292, 2019.

\bibitem{MR4230541}
H.~Yu.
\newblock On the metric theory of inhomogeneous {D}iophantine approximation: an
  {E}rd\h{o}s-{V}aaler type result.
\newblock {\em J. Number Theory}, 224:243--273, 2021.

\end{thebibliography}

\vspace{10ex}

{\small
\noindent Victor Beresnevich:\\
 Department of Mathematics,
 University of York,\\
 Heslington, York, YO10 5DD,
 England\\
 E-mail: {\tt victor.beresnevich@york.ac.uk}

\bigskip

{\small
\noindent Manuel Hauke:\\
 Department of Mathematics,
 University of York,\\
 Heslington, York, YO10 5DD,
 England\\
 E-mail: {\tt manuel.hauke@gmail.com}

\bigskip

\noindent Sanju Velani:\\
 Department of Mathematics,
 University of York,\\
 Heslington, York, YO10 5DD,
 England\\
 E-mail: {\tt sanju.velani@york.ac.uk}

}

\end{document}